\documentclass[11pt]{article}
\usepackage[margin=2.3cm]{geometry} 

\usepackage[utf8]{inputenc}
\usepackage{amsmath}
\usepackage{amsfonts}
\usepackage{amssymb}
\usepackage{tikz-cd} 
\usepackage{graphicx}
\usepackage{amsthm}
\usepackage{mathtools} 
\usepackage{setspace} 
\usepackage{colonequals} 
\usepackage{anyfontsize} 
\usepackage{quiver}
\usepackage{physics}
\usepackage{tikz}
\usepackage{mathdots}
\usepackage{yhmath}
\usepackage{dynkin-diagrams}
\usepackage{cancel}
\usepackage{array}
\usepackage{multirow}
\usepackage{amssymb}
\usepackage{tabularx}
\usepackage{extarrows}
\usepackage{booktabs}
\usetikzlibrary{fadings}
\usetikzlibrary{patterns}
\usetikzlibrary{shadows.blur}
\usetikzlibrary{shapes}
\usepackage{color}

\usepackage{microtype} 
\usepackage[nottoc]{tocbibind} 
\usepackage[none]{hyphenat} 
\usepackage{subcaption} 

\usepackage{titlesec}
\titleformat*{\section}{\Large\bfseries\sffamily} 
\titleformat*{\subsection}{\bfseries\sffamily}
\titleformat*{\subsubsection}{\bfseries\sffamily}

\usepackage{hyperref}
\usepackage{setspace}

\usepackage{hyperref}
\hypersetup{
    colorlinks=true,
    allcolors=blue
}

\usepackage{lmodern}

\usepackage{lscape}
\usepackage{graphicx}
\usepackage{amssymb}
\usepackage{mathtools}
\usepackage{amsmath}
\usepackage{multirow}
\usepackage{amsthm}
\usepackage{epstopdf}
\usepackage{youngtab}
\usepackage{bm}
\usepackage{tikz}
\usetikzlibrary{cd}
\usetikzlibrary{calc}
\usetikzlibrary{arrows}
\usetikzlibrary{arrows.meta}
\usetikzlibrary{shapes.arrows}
\usetikzlibrary{decorations.pathmorphing}
\usepackage{url}
\usepackage{array}
\begingroup
    \makeatletter
    \@for\theoremstyle:=definition,remark,plain\do{%
        \expandafter\g@addto@macro\csname th@\theoremstyle\endcsname{%
            \addtolength\thm@preskip\parskip
            }%
        }
\endgroup

\usepackage{caption} 

\allowdisplaybreaks

\newcommand{\dash}{\dashrightarrow}
\newcommand{\nin}{\noindent}

\newtheorem{thm}{Theorem}[section]

\newtheorem{propdefn}[thm]{Proposition--Definition}
\newtheorem{lem}[thm]{Lemma}
\newtheorem{cor}[thm]{Corollary}

\theoremstyle{definition}
\newtheorem{defn}[thm]{Definition}
\newtheorem{ex}[thm]{Example}
\newtheorem{rmk}[thm]{Remark}

\numberwithin{equation}{section}

\newcommand{\p}{\mathbb{P}}
\newcommand{\q}{\mathbb{Q}}
\newcommand{\C}{\mathbb{C}}
\newcommand{\F}{\mathbb{F}}
\newcommand{\Aa}{\mathbb{A}}
\newcommand{\RR}{\mathbb{R}}

\newcommand{\ZZ}{\mathbb{Z}}

\newcommand{\oo}{\mathcal{O}}

\newcommand{\cA}{\mathcal{A}}

\newcommand{\cE}{\mathcal{E}}

\DeclareMathOperator{\Bl}{Bl}

\DeclareMathOperator{\Spec}{Spec}

\DeclareMathOperator{\Exc}{Exc}

\DeclareMathOperator{\Bir}{Bir}

\DeclareMathOperator{\codim}{codim}

\DeclareMathOperator{\vp}{vp}

\DeclareMathOperator{\divv}{div}

\DeclareMathOperator{\cf}{coeff}
\DeclareMathOperator{\coreg}{coreg}

\DeclareMathOperator{\Sing}{Sing}

\DeclareMathOperator{\wt}{wt}
\DeclareMathOperator{\rk}{rank}
\DeclareMathOperator{\cone}{Cone}

\title{\textbf{\textsf{Birational geometry of log Calabi-Yau pairs $(\p^3,D)$ of  coregularity 2}}}
\author{Eduardo Alves da Silva}
\date{} 

\begin{document}

\maketitle

\renewcommand*\abstractname{\textsf{Abstract}}

\begin{abstract}
This paper aims to study the birational geometry of log Calabi-Yau pairs $(\p^3, D)$ of coregularity 2, where in this case $D$ is an irreducible normal quartic surface with canonical singularities. We completely classify which toric weighted blowups of a point will initiate a volume-preserving Sarkisov link starting with this pair. Depending on the type of singularity, our results point out that some of these weights do not work generically for a general member of the corresponding parameter space of quartics. 
\end{abstract}




\noindent \textbf{2020 Mathematics Subject Classification}: Primary 14E05, 14E07; Secondary  14E30, 14J17, 14J30

\noindent \textbf{Key words}: Sarkisov program, Calabi-Yau pairs, Cremona maps, volume-preserving maps

{
  \hypersetup{linkcolor=black}

\renewcommand*\contentsname{\textsf{Contents}}

\setcounter{tocdepth}{1}

\renewcommand{\baselinestretch}{0}

\textsf{\tableofcontents}

\renewcommand{\baselinestretch}{1.0}\normalsize
  
}

\section{Introduction}

The study of Calabi-Yau pairs has been an active research area in complex algebraic geometry. In part, this is because they can be seen as distinguished minimal models of the classical minimal model program (MMP) or its log version. Moreover, maximal log Calabi-Yau pairs (see Definition \ref{defn coreg}) have notable properties predicted from mirror symmetry \cite{hk1}. One important tool in the study of Calabi-Yau pairs is a relatively new version of the Sarkisov program \cite{cor1,hm} for volume-preserving maps between Mori fibered Calabi-Yau pairs obtained by Corti \& Kaloghiros \cite{ck}. See Theorem \ref{thm vp sark prog}.

The Sarkisov program asserts that any birational map between Mori fibered spaces can be written as a composition of a finite sequence of elementary maps, called Sarkisov links. This is very useful for studying Mori fibered spaces, which are outcomes of the MMP coming from uniruled varieties. The result by Corti \& Kaloghiros \cite{ck} can be interpreted as a generalization of this theorem with some additional structures.

In some cases, it might be possible to give an explicit description for all possible volume-preserving Sarkisov links between given Calabi-Yau pairs, which becomes very relevant for classification purposes. This was the main motivation of this paper in the case of log Calabi-Yau pairs $(\p^3, D)$ of coregularity 2.

Concerning the geometry of a log Calabi-Yau pair $(X, D_X)$, the most important discrete volume-preserving invariant is the coregularity that varies between 0 and $\dim X$. See Definition \ref{defn coreg}. The case of $\coreg(\p^3,D) = 2$ occurs if and only if $D$ is an irreducible normal quartic surface with canonical singularities. See Lemma \ref{coreg 2}.

In \cite{gue}, Duarte Guerreiro studied Sarkisov links initiated by the toric weighted blowup of a point in $\p^3$ or $\p^4$ using variation of GIT, and gave a complete classification of them with a description of the whole Sarkisov link. In the work \cite{acm}, Araujo, Corti \& Massarenti considered irreducible normal quartic surfaces with single canonical singularities of types $A_1$ and $A_2$ and solved the same problem in the volume-preserving context.

The main result of this paper extends the classification of terminal contractions given in \cite{acm}, contemplating all types of surface canonical singularities, the so-called Du Val singularities, which can be corresponded with simple-laced Dynkin diagrams of type ADE. In our context of $\coreg(\p^3,D) = 2$, we have the following:

\begin{thm}[See Theorem \ref{thm vp weights}]
     Let $(\p^3,D)$ be a log Calabi-Yau pair of coregularity 2 and $\pi \colon (X,D_X) \rightarrow (\p^3,D)$ be a volume-preserving toric $(1,a,b)$-weighted blowup of a torus invariant point. Then this point is necessarily a singularity of $D$ and,  up to permutation, the only possibilities for the weights, depending on the type of singularities, are listed in the following Table 
    \ref{table vp weights I}.

\end{thm}


\begin{table}[htp]
\begin{center}
\begin{tabular}{|c|c|}
\hline
type of singularity & volume-preserving weights \\
\hline
$A_1$ & (1,1,1) \\
$A_2$ & (1,1,1), (1,1,2), \\ 
$A_3$ & (1,1,1), (1,1,2), \textcolor{red}{(1,1,3)} \\
$A_4$ & (1,1,1), (1,1,2), \textcolor{red}{(1,1,3)}, \textcolor{purple}{(1,2,3)}\\
$A_5$ & (1,1,1), (1,1,2), \textcolor{red}{(1,1,3)}, \textcolor{purple}{(1,2,3)}\\
$A_6 $ & (1,1,1), (1,1,2), \textcolor{red}{(1,1,3)}, \textcolor{purple}{(1,2,3)}, \textcolor{cyan}{(1,2,5)}, \textcolor{green}{(1,3,4)}\\
$A_7 $ & (1,1,1), (1,1,2), \textcolor{red}{(1,1,3)}, \textcolor{purple}{(1,2,3)}, \textcolor{cyan}{(1,2,5)}, \textcolor{green}{(1,3,4)}, \textcolor{brown}{(1,3,5)}\\
\hline
$D_4$ & (1,1,1), (1,1,2)\\
$D_5$ & (1,1,1), (1,1,2), \textcolor{purple}{(1,2,3)}\\
$D_6$ & (1,1,1), (1,1,2), \textcolor{purple}{(1,2,3)}\\
$D_7$ & (1,1,1), (1,1,2), \textcolor{purple}{(1,2,3)}, \textcolor{green}{(1,3,4)}, \textcolor{brown}{(1,3,5)}\\
$D_8$ & (1,1,1), (1,1,2), \textcolor{purple}{(1,2,3)}, \textcolor{green}{(1,3,4)}, \textcolor{brown}{(1,3,5)}\\
$D_9$ & (1,1,1), (1,1,2), \textcolor{purple}{(1,2,3)}, \textcolor{green}{(1,3,4)}, \textcolor{brown}{(1,3,5)}, \textcolor{orange}{(1,4,5)}\\
$D_{10}$ & (1,1,1), (1,1,2), \textcolor{purple}{(1,2,3)}, \textcolor{green}{(1,3,4)}, \textcolor{brown}{(1,3,5)}, \textcolor{orange}{(1,4,5)}\\
\hline
$E_6$ & (1,1,1), (1,1,2), \textcolor{purple}{(1,2,3)}, \textcolor{green}{(1,3,4)}, \textcolor{brown}{(1,3,5)} \\
$E_7$ & (1,1,1), (1,1,2), \textcolor{purple}{(1,2,3)}, \textcolor{green}{(1,3,4)}, \textcolor{brown}{(1,3,5)} \\
$E_8$ & (1,1,1), (1,1,2), \textcolor{purple}{(1,2,3)}, \textcolor{green}{(1,3,4)}, \textcolor{brown}{(1,3,5)}, \textcolor{orange}{(1,4,5)}\\
\hline
\end{tabular}
\caption{Table summarizing volume-preserving weights, up to permutation.} 
\label{table vp weights I}
\end{center}
\end{table}

The following result is the volume-preserving version of \cite[Theorem 1.1]{gue} for the case where $\coreg(\p^3,D)=2$. The toric description of the weighted blowup allows us to encompass all types of strict canonical singularities.

\begin{cor}[See Corollary \ref{cor vp weights sark}]
     Let $(\p^3,D)$ be a log Calabi-Yau pair of coregularity 2 and $\pi \colon (X,D_X) \rightarrow (\p^3,D)$ be a volume-preserving toric $(1,a,b)$-weighted blowup of a torus invariant point. Then this point is necessarily a singularity of $D$ and,  up to permutation, the only possibilities for the weights initiating a volume-preserving Sarkisov link, depending on the type of singularities, are listed in the following Table \ref{table vp weights II}.

\end{cor}

\newpage

\begin{table}[htp]
\begin{center}
\begin{tabular}{|c|c|}
\hline
type of singularity & volume-preserving weights \\
\hline
$A_1$ & (1,1,1) \\
$A_2$ & (1,1,1), (1,1,2) \\ 
$A_3$ & (1,1,1), (1,1,2) \\
$A_4$ & (1,1,1), (1,1,2),  \textcolor{purple}{(1,2,3)}\\
$A_5$ & (1,1,1), (1,1,2), \textcolor{purple}{(1,2,3)}\\
$A_{\ge 6} $ & (1,1,1), (1,1,2), \textcolor{purple}{(1,2,3)}, \textcolor{cyan}{(1,2,5)}\\
\hline
$D_4$ & (1,1,1), (1,1,2)\\
$D_{\ge 5}$ & (1,1,1), (1,1,2), \textcolor{purple}{(1,2,3)}\\
\hline
$E_6$ & (1,1,1), (1,1,2), \textcolor{purple}{(1,2,3)} \\
$E_7$ & (1,1,1), (1,1,2), \textcolor{purple}{(1,2,3)} \\
$E_8$ & (1,1,1), (1,1,2), \textcolor{purple}{(1,2,3)} \\
\hline
\end{tabular}
\caption{Table summarizing volume-preserving weights initiating Sarkisov links, up to permutation.} 
\label{table vp weights II}
\end{center}
\end{table}

\begin{rmk}[Important]
   The weights displayed in colors in Tables \ref{table vp weights I} \& \ref{table vp weights II} are not volume-preserving for a generic quartic having that type of singularity. By this, we mean that some closed conditions exist on the coefficients of the equation of $D$ in order for these weights to satisfy the volume-preserving property. These conditions will be stated later on in this work together with criteria to detect such singularities. See Subsection \ref{criteria A_n}. 
\end{rmk}

\begin{rmk}\label{rmk (1,1,3)} While the weights $(1,1,3)$ yield a volume-preserving toric divisorial contraction for some $D'$ with a singularity of type $A_{\ge 3}$, it will not generate a Sarkisov link between terminal Mori fibered Calabi-Yau pairs. See \cite[Lemma 3.2]{gue} or notice that after computing the whole link through variation of GIT, we will get a codomain with worse than terminal singularities, which is not allowed in the Sarkisov program \cite{cor1,hm}.
For the complete description of the Sarkisov links initiated by the remaining weights, see \cite[Table 1]{gue}.

\end{rmk}

Our results can be regarded as a first step in the explicit classification of log Calabi-Yau pairs $(\p^3, D)$ of coregularity 2, up to volume-preserving equivalence, providing advancement of a technology designed to explicitly handle volume-preserving birational maps of threefold Mori fibered Calabi-Yau pairs.

\paragraph{Outline of the proofs of the Theorem \ref{thm vp weights} \& Corollary \ref{cor vp weights sark}.} The volume-preserving property of the toric $(1,a,b)$-weighted blowup is equivalent to the vanishing of $a(E,\p^3,D)$, the discrepancy of the corresponding exceptional divisor $E$ with respect to $(\p^3,D)$. To compute this discrepancy, we used a toric description of the weighted blowup. Our strategy was to realize the divisorial valuation $\nu_E$ associated to $E$ after a sequence of ordinary blowups at points or nonsingular curves such that the valuation on $\p^3$ corresponding to the last exceptional divisor coincides with the valuation associated to $E$ on $\p^3$.

The different types of canonical singularities can be expressed by conditions involving the coefficients of the equation of $D$ and the same holds for the weights $(1,a,b)$ to satisfy the volume-preserving property. By determining all these conditions and comparing them, we obtained the desired weights.

\subsection{Structure of the paper}

Throughout this paper, our ground field will be $\C$, or more generally, any algebraically closed field of characteristic zero. 

In Section \ref{CY pairs and Sark} we will give an overview of Calabi-Yau pairs and the Sarkisov program, introduce some natural classes of singularities of pairs, and give a short compendium on singularities of quartic surfaces. In Section \ref{vp Sark links}, we will approach the context of the problem involving Theorem \ref{thm vp weights} \& Corollary \ref{cor vp weights sark}. In Section \ref{deter vp weights}, with all the results of previous sections, we will finally determine the desired volume-preserving weights.

\subsection{Acknowledgements}

I would like to thank my PhD advisor Carolina Araujo for guiding me in the process of construction of this paper and for her tremendous patience in our several discussions. I am also grateful to Tiago Duarte Guerreiro, a mathematical partner whom I had the opportunity to meet at GAeL XXVIII, for many important clarifications about the theme. I also thank CNPq (Conselho Nacional de Desenvolvimento Científico e Tecnológico) for the financial support during my studies at IMPA, FAPERJ (Fundação de Amparo à Pesquisa do Estado do Rio de Janeiro) for supporting my attendance at events that were very important for meeting people to discuss this work, and the support of the grant ERC StG Saphidir No. 101076412 during my initial postdoctoral period at the University of Paris-Saclay.

\section{Overview about Calabi-Yau pairs and the Sarkisov program}\label{CY pairs and Sark}

In this section, we give an overview of log Calabi-Yau pairs and the volume-preserving Sarkisov program.

\begin{defn}
A \textit{log Calabi-Yau pair} or simply \textit{Calabi-Yau pair} is a log canonical pair $(X,D)$ consisting of a normal projective variety $X$ and a reduced Weil divisor $D$ on $X$ such that $K_X + D \sim 0$. This condition implies the existence of a top degree rational differential form $\omega = \omega_{X,D} \in \Omega_X^n$, unique up to nonzero scaling, such that $D + \divv(\omega)=0$. By abuse of language, we call this differential the \textit{volume form} of $(X,D)$.
\end{defn}

Since $K_X + D \sim 0$, we have that $K_X + D$ is readily Cartier, and hence all the discrepancies with respect to the pair $(X,D)$ are integer numbers.

We now introduce some classes of pairs taking into account the singularities of the ambient varieties and the divisors.

\begin{defn}\label{defn pairs}
    We say that a pair $(X,D)$ is $(t,c)$, respectively, $(t,lc)$, if $X$ has terminal singularities and the pair $(X,D)$ has canonical, respectively log canonical singularities. We say that a pair $(X,D)$ is $\q$-factorial if $X$ is $\q$-factorial.  
\end{defn}

If $(X,D)$ is (t,lc), then $a(E,X,D_X) \leq 0$ implies $a(E,X,D_X) = -1$ or $0$.

\begin{propdefn}[cf. \cite{km} Lemma 2.30]\label{prop crep bir morp}
    A proper birational morphism $f \colon (Z,D_Z) \rightarrow (X,D_X)$ is called crepant if $f_*D_Z=D_X$ and $f^*(K_X + D_X) \sim K_Z+ D_Z$. The term ``crepant'' (coined by Reid) refers to the fact that for every divisor $E$ over $X$, one has $a(E,Z,D_Z)=a(E,X,D_X)$. 
\end{propdefn}

\begin{defn}\label{defn crep bir map}
    A birational map of pairs $\phi \colon (X,D_X) \dash (Y,D_Y)$ is called \textit{crepant} if it admits a resolution 
    \begin{center}
\begin{tikzcd}
	& {(Z,D_Z)} \\
	{(X,D_X)} && {(Y,D_Y)}
	\arrow["\phi", dashed, from=2-1, to=2-3]
	\arrow["p"', from=1-2, to=2-1]
	\arrow["q", from=1-2, to=2-3]
\end{tikzcd}
    \end{center}

\noindent in such a way that $p$ and $q$ are crepant birational morphisms.

\end{defn}

This definition is equivalent to asking that $a(E,X,D_X)=a(E,Y,D_Y)$ for every valuation $E$ of $K(X) \simeq K(Y)$.

For Calabi-Yau pairs, the notion of \textit{crepant birational equivalence} becomes \textit{volume-preserving equivalence}, since $\phi^*\omega_{Y,D_Y} = \lambda \omega_{X,D_X}$, for some $\lambda \in \C^*$. In this case, we call such $\phi$ a \textit{volume-preserving} map.

\begin{rmk}\label{comp of vp is vp}  
    As an immediate consequence of the definition, a composition of volume-preserving maps is volume-preserving. So the set of volume-preserving self-maps of a given Calabi-Yau pair $(X,D)$ forms a group, denoted by $\Bir^{\vp}(X,D)$. In particular, this group is a subgroup of $\Bir(X)$.    
\end{rmk}

\begin{defn}
    A \textit{Mori fibered Calabi-Yau pair} is a $\q$-factorial (t,lc) Calabi-Yau pair $(X,D)$ with a Mori fibered space structure on $X$, that is, a morphism $f \colon X \rightarrow S$, to a lower dimensional variety $S$, such that $f_*\oo_X=\oo_S$, $-K_X$ is $f$-ample, and $\rho(X/S)=\rho(X) - \rho(S)=1$. 
\end{defn}

\paragraph{Sarkisov program.} It is shown in \cite{cor1,hm} that any birational map between Mori fibered spaces is a composition of \textit{Sarkisov links}.

We recall now the definition of the 4 types of Sarkisov links and Mori divisorial contractions to fix terminology regarding the latter. 

A \textit{Mori divisorial contraction} is an extremal divisorial contraction $f\colon  Z \rightarrow X$ from a $\q$-factorial terminal variety $Z$ of an extremal ray $R$ with $K_Z \cdot R < 0$. In particular, $X$ also has $\q$-factorial terminal singularities and $\rho(Z/X)=1$.

A \textit{Sarkisov link} is a birational map between Mori fibered spaces 
\[\begin{tikzcd}
	X & {X'} \\
	S & {S'}
	\arrow[dashed, from=1-1, to=1-2]
	\arrow[from=1-1, to=2-1]
	\arrow[from=1-2, to=2-2]
\end{tikzcd}\]
that fits into one of the following commutative diagrams. Observe that in this description, we dispose the varieties of same Picard rank at the same height.

\begin{enumerate}
    \item A \textit{Sarkisov link of type I} fits into
\begin{center}
\begin{tikzcd}
	& Z & {X'} \\
	X && {S'} \\
	S
	\arrow[from=1-2, to=2-1]
	\arrow[dashed, from=1-2, to=1-3]
	\arrow[from=2-1, to=3-1]
	\arrow[from=2-3, to=3-1]
	\arrow[from=1-3, to=2-3]
\end{tikzcd}
\end{center}

\nin where $Z \rightarrow X$ is a Mori divisorial contraction, and $Z \dash X'$ is a sequence of Mori flips, flops and antiflips. Notice that $\rho(S'/S)=1$.

\item A \textit{Sarkisov link of type II} fits into
\begin{center}
\begin{tikzcd}
	& Z & {Z'} \\
	X &&& {X'} \\
	S &&& {S'}
	\arrow[from=2-1, to=3-1]
	\arrow[from=1-2, to=2-1]
	\arrow[dashed, from=1-2, to=1-3]
	\arrow[from=1-3, to=2-4]
	\arrow[from=2-4, to=3-4]
	\arrow[Rightarrow, no head, from=3-1, to=3-4]
\end{tikzcd}
\end{center}

\nin where $Z \rightarrow X$ and $Z' \rightarrow X'$ are Mori divisorial contractions, and $Z \dash X'$ is a sequence of Mori flips, flops and antiflips. Notice here that $S=S'$.

    \item A \textit{Sarkisov link of type III} fits into
\begin{center}
\begin{tikzcd}
	X & Z \\
	S && {X'} \\
	&& {S'}
	\arrow[from=1-2, to=2-3]
	\arrow[from=1-1, to=2-1]
	\arrow[from=2-3, to=3-3]
	\arrow[from=2-1, to=3-3]
	\arrow[dashed, from=1-1, to=1-2]
\end{tikzcd}
\end{center}

\nin where $X \dash Z$ is a sequence of Mori flips, flops and antiflips, and $Z \rightarrow X'$ is a Mori divisorial contraction. Notice that $\rho(S/S')=1$.

    \item A \textit{Sarkisov link of type IV} fits into
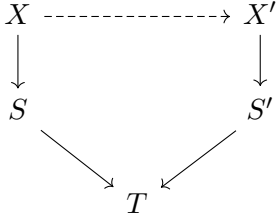
\begin{center}
\begin{tikzcd}
	X && {X'} \\
	S && {S'} \\
	& T
	\arrow[dashed, from=1-1, to=1-3]
	\arrow[from=1-1, to=2-1]
	\arrow[from=1-3, to=2-3]
	\arrow[from=2-1, to=3-2]
	\arrow[from=2-3, to=3-2]
\end{tikzcd}
\end{center}

\nin where $X \dash X'$ is a sequence of Mori flips, flops and antiflips, and $S \rightarrow T$ and $S' \rightarrow T$ are Mori contractions. Hence we have that $\rho(S/T)=\rho(S'/T)=1$.
\end{enumerate}

\begin{defn}\label{defn vp Sark links}
    A \textit{volume-preserving Sarkisov link} is a Sarkisov link as previously described with the following additional data and property: there exist divisors $D_X$ on $X$, $D_{X'}$ on $X'$, $D_Z$ on $Z$, and $D_{Z'}$ on $Z'$, making $(X,D_X),(X',D_{X'}),(Z,D_Z)$ and $(Z',D_{Z'})$ (t,lc) Calabi-Yau pairs, and all the divisorial contractions and extractions, Mori flips, flops and antiflips that constitute the Sarkisov link are volume-preserving for these Calabi-Yau pairs.
\end{defn}

\begin{rmk}\label{all in a vp sark is vp}
    A Sarkisov link of any type is volume-preserving if and only if the corresponding divisorial contractions and extractions, and isomorphisms in codimension 1 involved (when it is the case) are volume-preserving.
\end{rmk}

We state now the result of Corti \& Kaloghiros, which holds in all dimensions.

\begin{thm}[cf. \cite{ck} Theorem 1.1]\label{thm vp sark prog}
    Any volume-preserving map between Mori fibered Calabi-Yau pairs is a composition of volume-preserving Sarkisov links.
\end{thm}
\subsection{The coregularity}

From the perspective of the Minimal Model Program, being log canonical is the minimal requirement about how bad the singularities are. There exists a class of singularities, namely, \textit{divisorial log terminal}, which is slightly better than log canonical.

\begin{defn}[cf. \cite{km} Definition 2.37]\label{def dlt}
Let $(X,D)$ be a pair where $D = \sum a_iD_i $ is a sum of distinct prime divisors with $0 \leq a_i \leq 1$. Assume that $K_X+D$ is $\q$-Cartier. Such a pair is \textit{divisorial log terminal} (dlt for short) if and only if there exists a closed subset $Z \subset X$ with the following properties:
\begin{enumerate}
    \item $X \setminus Z$ is nonsingular and $D|_{X \setminus Z}$ is a simple normal crossing divisor.
    
    \item If $f \colon Y \rightarrow X$ is birational and $E \subset Y$ is an irreducible divisor such that its center in $X$\footnote{The \textit{center} of a divisor $E$ over $X$ is the closure of the set-theoretic image $f(E) \subset X$.} is contained in $Z$, then $a(E,X,D) > -1$.
\end{enumerate}
\end{defn}

In our context of log Calabi-Yau geometry, we can always find a \textit{dlt modification} for any Calabi-Yau pair \cite[Theorem 1.7]{ck}. More precisely, given a Calabi-Yau pair $(X,D_X)$, there always exists a volume-preserving morphism $f \colon (Y,D_Y) \rightarrow (X,D_X)$ where $(Y,D_Y)$ is a $\q$-factorial dlt Calabi-Yau pair and $Y$ has at worst terminal singularities.

For dlt Calabi-Yau pairs, the log canonical centers are well understood via the \textit{strata} of the boundary divisor, which is the collection of irreducible components of all possible finite intersections between the divisors appearing in its support. An element of the strata is called \textit{stratum}.

\begin{thm}[cf. \cite{kol} Theorem 4.16]\label{lc centers dlt}
Let $(X,D)$ be a dlt pair (not necessarily Calabi-Yau) and consider $D_1,\ldots,D_r$ the irreducible divisors appearing in $D$ with coefficient 1 (in the Calabi-Yau case this is automatically satisfied by all components of $D$). Set $I \coloneqq \{1,\ldots,r \}$.
\begin{enumerate}
    \item The log canonical centers of $(X,D)$ are exactly the irreducible components of $D_J \coloneqq \bigcap_{j \in J} D_j$ for any subset $J \subset I$. For $J=\emptyset$, we define $D_{\emptyset} \coloneqq X$.

    \item Every log canonical center is normal and has pure codimension $|J|$.
\end{enumerate}
\end{thm}

Since volume-preserving maps between Calabi-Yau pairs preserve discrepancies, any concept expressed in terms of them will be a volume-preserving invariant. In particular, these maps send log canonical centers onto log canonical centers. Using this fact, it is possible to show that the dimension of a minimal log canonical center (with respect to the inclusion) on a dlt modification of a Calabi-Yau pair $(X,D_X)$ is a volume-preserving invariant.

\begin{defn}\label{defn coreg}
    The \textit{coregularity} $\coreg(X,D_X)$ is defined to be the dimension of a minimal log canonical center in a dlt modification $f \colon (Y,D_Y) \rightarrow (X,D_X)$. A Calabi-Yau pair $(X,D_X)$ is called \textit{maximal} if $\coreg(X,D_X)=0$.
\end{defn}

Thus we have $0 \leq \coreg(X,D_X) \leq \dim X$ for any Calabi-Yau pair $(X,D_X)$. The pairs with ambient varieties with Kawamata log terminal singularities and maximum coregularity $\coreg(X,D_X) = \dim X$ are necessarily of the form $(X,0)$. By definition of Calabi-Yau pair, this implies that $X$ is a Calabi-Yau variety. 

It is also possible to define the notion of coregularity in terms of the dimension of the ambient variety and the corresponding \textit{dual complex} $\mathcal{D}(X,D_X)$ of the pair $(X,D_X)$. In a broad sense, this object is a CW-complex that encodes the geometry of the log canonical centers. We have $\coreg(X,D_X) = \dim X - \dim \mathcal{D}(X,D_X) -1$. Notice that the case of minimum coregularity corresponds to the case where $\dim \mathcal{D}(X,D_X)$ is maximal, which is one justification for the terminology ``maximal pair''. We refer the reader to \cite{duc,kx,mor1,mor2} for more details.

The classification problem of Calabi-Yau pairs up to volume-preserving equivalence can be refined through the coregularity. In the recent work \cite{duc}, Ducat classified all pairs of the form $(\p^3,D)$ with coregularity less than or equal to one, up to volume-preserving equivalence. For those with minimum possible coregularity, he showed that they are all volume-preserving equivalent to toric pairs $(T, D_T)$. 

The following lemma explains the context of the present paper:

\begin{lem}\label{coreg 2}
    Let $(\p^3,D)$ be a Calabi-Yau pair. One has $\coreg(\p^3,D) = 2$ if and only if $D$ is an irreducible normal quartic surface with at worst canonical singularities
\end{lem}

\begin{proof}
    ($\Rightarrow$) Consider $f \colon (X,D_X) \rightarrow (\p^3,D)$ a dlt modification of $(\p^3,D)$. By Theorem \ref{lc centers dlt} combined with \cite[Theorem 1.7]{ck} we have that $X$ is terminal and $D_X$ is irreducible and normal. Since $(X,D_X)$ is a dlt pair, by \cite[Theorem 2.44]{km} there exists a log resolution $g \colon X' \rightarrow X$ of $(X,D_X)$ such that $a(E_i,X,D_X) > -1$ for every exceptional divisor $E_i \subset X'$. 
    
    Discrepancies with respect to Calabi-Yau pairs are always integer numbers because the log canonical divisor is Cartier. This implies that $a(E_i,X,D_X) \geq 0$ for every exceptional divisor $E_i \subset X'$, and hence $(X,D_X)$ is a (t,c) pair. \cite[Lemma 2.18]{acm} ensures that canonicity is preserved under volume-preserving maps, which guarantees that $(\p^3,D)$ is also (t,c). Making use of the adjunction formula, then $D$ has canonical singularities and, in particular, it is normal. By \cite[Proposition 2.6]{acm}, one has the irreducibility of $D$.\\

\nin ($\Leftarrow$) Let us show that $(\p^3,D)$ is already a dlt Calabi-Yau pair. Making use of the adjunction formula, we have that $(\p^3,D)$ is (t,c). Taking $Z = \Sing(D)$ in \cite[Definition 2.37]{km} or by \cite[Theorem 2.44]{km}, we have that $(\p^3,D)$ is dlt. Then $\coreg(\p^3,D) = 2$ is immediate by Theorem \ref{lc centers dlt}. In particular, notice that  $ \Sing(D) = \emptyset$ implies that the pair $(\p^3,D)$ is terminal.    
\end{proof}

As pointed out by Ducat, the case of $\coreg(\p^3,D) = 2$ is the hardest one in terms of an explicit classification of volume-preserving equivalence classes. To illustrate this fact, Oguiso \cite{og} exhibited two nonsingular isomorphic quartic surfaces $D, D' \subset \p^3$ such that there is no $\varphi \in \Bir(\p^3)$ mapping $D$ birationally onto $D'$. By \cite[Proposition 2.6]{acm} any map having this property would automatically be volume-preserving for the pairs $(\p^3,D)$ and $(\p^3,D')$.

\subsection{A short compendium on singularities of quartic surfaces}

To understand the birational geometry of Calabi-Yau pairs $(\p^3, D)$ of coregularity 2,  Lemma \ref{coreg 2} indicates that we must study canonical singularities realizable by irreducible normal quartic surfaces in $\p^3$. In this subsection and the following parts of this paper, we will refer to such varieties as simply quartic surfaces.

It is well known that canonical surface singularities are precisely nonsingular points together with Du Val surface singularities. See for instance \cite{rei3}. It can be proved that a Du Val surface singularity $P \in S$ is analytically isomorphic to a surface singularity $0 \in \{f=0\} \subset \Aa^3_{(x,y,z)}$, where $f$ is one of the following equations together with the corresponding resolution graph:
\begin{itemize}
    \item[$A_n$:] $x^2+y^2+z^{n+1}$ for $n \geq 1$,\\
    \dynkin[root radius=.05cm,edge length=0.75cm] A{}
    \item[$D_n$:] $x^2+y^2z+z^{n-1}$ for $n \geq 4$,\\
     \dynkin[root radius=.05cm,edge length=0.75cm] D{}
    \item[$E_6$:] $x^2+y^3+z^4$,\\
     \dynkin[root radius=.05cm,edge length=0.75cm] E6
    \item[$E_7$:] $x^2+y^3+yz^3$,\\
    \dynkin[root radius=.05cm,edge length=0.75cm] E7
    \item[$E_8$:] $x^2+y^3+z^5$,\\
    \dynkin[root radius=.05cm,edge length=0.75cm] E8
\end{itemize}

In \cite{in,um1,um2} minimal desingularizations of a normal quartic surface are classified. The case of canonical singularities leads us to K3 surfaces. See \cite[Lemma 5.1.8]{alv}.

Regarding the configuration of canonical singularities on quartic surfaces, it is not practical to deal with the problem by direct analysis of its homogeneous equation. There exist many papers on this question in the last decades and some answers were given by Kato \& Naruki and Urabe \cite{kn,ur1,ur2}. Some years after their works, an explicit enumeration of configurations of such singularities was obtained by Yang \cite{ya}. We warn the reader that the list in \cite{ya} is enormous. All those authors based their works on the geometry of K3 surfaces. 

\begin{ex}
    Let $S \subset \p^3$ be a quartic surface having a single singularity of type $A_n$ or $D_n$. It follows that $n$ is at most $19$ in both cases by comparing the Picard number of $S$ with 20, the maximum possible Picard number for a K3 surface. We have the existence, respectively nonexistence, of such a surface $S$ with an $A_{19}$ singularity (therefore unique), respectively $D_{19}$, due to the classification of sublattices of type $A_{19}$ and $D_{19}$ in the K3 lattice \cite[Propositions 1 \& 2]{kn}, \cite[Corollary 0.3 \& Proposition 3.5]{ur3}. Kato \& Naruki \cite{kn} show that there exists a unique $S \subset \p^3$ having an $A_{19}$ singularity, up to automorphisms of $\p^3$. In fact, taking $P = (1:0:0:0) \in \p^3_{(w:x:y:z)}$ to be the singular point of $S$, its equation in affine coordinates in $\{w \neq 0\} \simeq \Aa^3_{(x,y,z)}$ is given by
    \begin{multline*}
    16(x^2 +y^2) +32xz^2- 16y^3 +16z^4 -32yz^3 +8(2x^2- 2xy+5y^2)z^2 \\ 
+ 8(2x^3-5x^2y-6xy^2 -7y^3)z+20x^4 +44x^3y+65x^2y^2 +40xy^3 +41y^4 =0.
\end{multline*}

In particular, we must have $\rho(S)=1$. It is also possible to justify the uniqueness of $S$ in such a case by analyzing the corresponding parameter space of quartics. 
\end{ex}

\subsection{Explicit resolution of Du Val singularities}\label{exp res Du Val}

In the case of canonical surface singularities, it is known that the intersection graph of the exceptional divisor of a minimal resolution is given by a simple-laced Dynkin diagram of type ADE.

In this subsection, we collect some results about the resolution of the Du Val singularities whose proofs can be found in \cite[Section 5.1]{alv}. They tell us how many blowups at singular points we need to do and indicate the configuration of the singularities in the intermediate steps of the desingularization process.

\begin{lem}\label{res A_n} The resolution of a Du Val surface singularity of type $A_n$ can be reached by performing $\left\lceil \dfrac{n}{2} \right\rceil$ blowups at a sequence of infinitely near base points lying over the singular point.
\end{lem}

The previous lemma tells us that each point blowup weakens the singularity from $A_n$ to $A_{n-2}$. From this, we can also infer the behavior of the corresponding exceptional divisor along this resolution process. The dynamic of the exceptional divisor of an $A_n$ singularity is to grow from the middle to the ``outermost'' components. The following pictures illustrate this phenomenon:

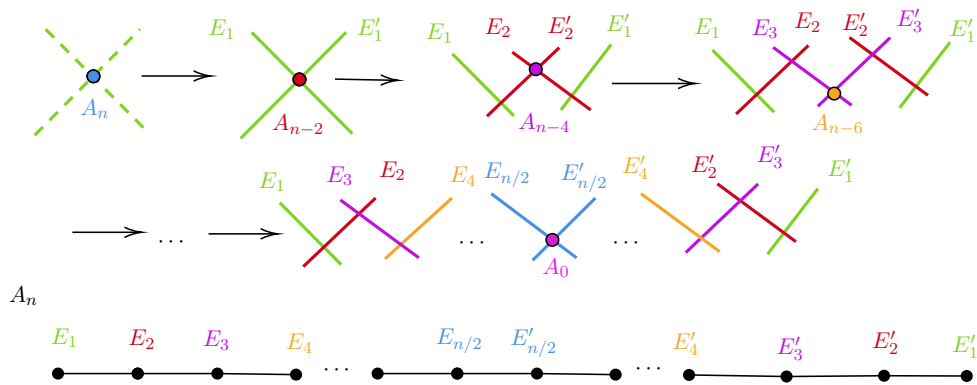
\begin{figure}[htb]
\centering
\resizebox{0.8\textwidth}{!}%
{

\tikzset{every picture/.style={line width=0.75pt}} 

\begin{tikzpicture}[x=0.75pt,y=0.75pt,yscale=-1,xscale=1]

\draw  [fill={rgb, 255:red, 74; green, 144; blue, 226 }  ,fill opacity=1 ] (70.99,46.6) .. controls (70.99,44.39) and (72.78,42.6) .. (74.99,42.6) .. controls (77.19,42.6) and (78.99,44.39) .. (78.99,46.6) .. controls (78.99,48.81) and (77.19,50.6) .. (74.99,50.6) .. controls (72.78,50.6) and (70.99,48.81) .. (70.99,46.6) -- cycle ;
\draw [color={rgb, 255:red, 126; green, 211; blue, 33 }  ,draw opacity=1 ][line width=1.5]  [dash pattern={on 5.63pt off 4.5pt}]  (44.6,16.4) -- (108.6,81.4) ;
\draw [color={rgb, 255:red, 126; green, 211; blue, 33 }  ,draw opacity=1 ][line width=1.5]  [dash pattern={on 5.63pt off 4.5pt}]  (37.6,84.4) -- (104.6,16.4) ;
\draw  [fill={rgb, 255:red, 74; green, 144; blue, 226 }  ,fill opacity=1 ] (70.99,46.6) .. controls (70.99,44.39) and (72.78,42.6) .. (74.99,42.6) .. controls (77.19,42.6) and (78.99,44.39) .. (78.99,46.6) .. controls (78.99,48.81) and (77.19,50.6) .. (74.99,50.6) .. controls (72.78,50.6) and (70.99,48.81) .. (70.99,46.6) -- cycle ;
\draw    (105.6,46.4) -- (147.6,46.4) ;
\draw [shift={(149.6,46.4)}, rotate = 180] [color={rgb, 255:red, 0; green, 0; blue, 0 }  ][line width=0.75]    (10.93,-3.29) .. controls (6.95,-1.4) and (3.31,-0.3) .. (0,0) .. controls (3.31,0.3) and (6.95,1.4) .. (10.93,3.29)   ;
\draw  [fill={rgb, 255:red, 74; green, 144; blue, 226 }  ,fill opacity=1 ] (201.99,48.6) .. controls (201.99,46.39) and (203.78,44.6) .. (205.99,44.6) .. controls (208.19,44.6) and (209.99,46.39) .. (209.99,48.6) .. controls (209.99,50.81) and (208.19,52.6) .. (205.99,52.6) .. controls (203.78,52.6) and (201.99,50.81) .. (201.99,48.6) -- cycle ;
\draw [color={rgb, 255:red, 126; green, 211; blue, 33 }  ,draw opacity=1 ][line width=1.5]    (175.6,18.4) -- (239.6,83.4) ;
\draw [color={rgb, 255:red, 126; green, 211; blue, 33 }  ,draw opacity=1 ][line width=1.5]    (168.6,86.4) -- (235.6,18.4) ;
\draw  [fill={rgb, 255:red, 209; green, 2; blue, 27 }  ,fill opacity=1 ] (201.99,48.6) .. controls (201.99,46.39) and (203.78,44.6) .. (205.99,44.6) .. controls (208.19,44.6) and (209.99,46.39) .. (209.99,48.6) .. controls (209.99,50.81) and (208.19,52.6) .. (205.99,52.6) .. controls (203.78,52.6) and (201.99,50.81) .. (201.99,48.6) -- cycle ;
\draw    (228.6,48.4) -- (266.6,49.16) ;
\draw [shift={(268.6,49.2)}, rotate = 181.15] [color={rgb, 255:red, 0; green, 0; blue, 0 }  ][line width=0.75]    (10.93,-3.29) .. controls (6.95,-1.4) and (3.31,-0.3) .. (0,0) .. controls (3.31,0.3) and (6.95,1.4) .. (10.93,3.29)   ;
\draw [color={rgb, 255:red, 126; green, 211; blue, 33 }  ,draw opacity=1 ][line width=1.5]    (303.2,32.2) -- (342.39,72) -- (342.39,72) ;
\draw [color={rgb, 255:red, 126; green, 211; blue, 33 }  ,draw opacity=1 ][line width=1.5]    (372.2,68.2) -- (382.62,54.52) -- (404.2,26.2) ;
\draw    (404.6,51.2) -- (455.6,51.2) ;
\draw [shift={(457.6,51.2)}, rotate = 180] [color={rgb, 255:red, 0; green, 0; blue, 0 }  ][line width=0.75]    (10.93,-3.29) .. controls (6.95,-1.4) and (3.31,-0.3) .. (0,0) .. controls (3.31,0.3) and (6.95,1.4) .. (10.93,3.29)   ;
\draw [color={rgb, 255:red, 209; green, 2; blue, 27 }  ,draw opacity=1 ][line width=1.5]    (325.2,72.2) -- (371.2,28.2) ;
\draw [color={rgb, 255:red, 209; green, 2; blue, 27 }  ,draw opacity=1 ][line width=1.5]    (391.2,68.2) -- (341.29,31.1) ;
\draw  [fill={rgb, 255:red, 189; green, 16; blue, 224 }  ,fill opacity=1 ] (352.25,42.07) .. controls (352.25,39.86) and (354.04,38.07) .. (356.25,38.07) .. controls (358.46,38.07) and (360.25,39.86) .. (360.25,42.07) .. controls (360.25,44.28) and (358.46,46.07) .. (356.25,46.07) .. controls (354.04,46.07) and (352.25,44.28) .. (352.25,42.07) -- cycle ;
\draw [color={rgb, 255:red, 126; green, 211; blue, 33 }  ,draw opacity=1 ][line width=1.5]    (468.6,29.8) -- (507.79,69.6) -- (507.79,69.6) ;
\draw [color={rgb, 255:red, 126; green, 211; blue, 33 }  ,draw opacity=1 ][line width=1.5]    (587.6,66.8) -- (598.02,53.12) -- (619.6,24.8) ;
\draw [color={rgb, 255:red, 209; green, 2; blue, 27 }  ,draw opacity=1 ][line width=1.5]    (483.6,70.8) -- (529.6,26.8) ;
\draw [color={rgb, 255:red, 209; green, 2; blue, 27 }  ,draw opacity=1 ][line width=1.5]    (605.6,58.8) -- (555.69,21.7) ;
\draw [color={rgb, 255:red, 189; green, 16; blue, 224 }  ,draw opacity=1 ][line width=1.5]    (535.6,65.8) -- (581.6,21.8) ;
\draw [color={rgb, 255:red, 189; green, 16; blue, 224 }  ,draw opacity=1 ][line width=1.5]    (557.08,65.16) -- (507.17,28.06) ;
\draw  [fill={rgb, 255:red, 245; green, 166; blue, 35 }  ,fill opacity=1 ] (541.99,57.67) .. controls (541.99,55.46) and (543.78,53.67) .. (545.99,53.67) .. controls (548.19,53.67) and (549.99,55.46) .. (549.99,57.67) .. controls (549.99,59.88) and (548.19,61.67) .. (545.99,61.67) .. controls (543.78,61.67) and (541.99,59.88) .. (541.99,57.67) -- cycle ;
\draw    (61.6,146.6) -- (103.6,146.6) ;
\draw [shift={(105.6,146.6)}, rotate = 180] [color={rgb, 255:red, 0; green, 0; blue, 0 }  ][line width=0.75]    (10.93,-3.29) .. controls (6.95,-1.4) and (3.31,-0.3) .. (0,0) .. controls (3.31,0.3) and (6.95,1.4) .. (10.93,3.29)   ;
\draw    (148.2,147.6) -- (190.2,147.6) ;
\draw [shift={(192.2,147.6)}, rotate = 180] [color={rgb, 255:red, 0; green, 0; blue, 0 }  ][line width=0.75]    (10.93,-3.29) .. controls (6.95,-1.4) and (3.31,-0.3) .. (0,0) .. controls (3.31,0.3) and (6.95,1.4) .. (10.93,3.29)   ;
\draw [color={rgb, 255:red, 126; green, 211; blue, 33 }  ,draw opacity=1 ][line width=1.5]    (193.4,127.6) -- (232.59,167.4) -- (232.59,167.4) ;
\draw [color={rgb, 255:red, 209; green, 2; blue, 27 }  ,draw opacity=1 ][line width=1.5]    (208.4,168.6) -- (254.4,124.6) ;
\draw [color={rgb, 255:red, 245; green, 166; blue, 35 }  ,draw opacity=1 ][line width=1.5]    (260.4,163.6) -- (303,124.6) ;
\draw [color={rgb, 255:red, 189; green, 16; blue, 224 }  ,draw opacity=1 ][line width=1.5]    (281.88,162.96) -- (231.97,125.86) ;
\draw [color={rgb, 255:red, 74; green, 144; blue, 226 }  ,draw opacity=1 ][line width=1.5]    (354.2,164.6) -- (394.2,124.6) ;
\draw [color={rgb, 255:red, 74; green, 144; blue, 226 }  ,draw opacity=1 ][line width=1.5]    (382.2,163) -- (327.97,121.86) ;
\draw  [fill={rgb, 255:red, 216; green, 30; blue, 205 }  ,fill opacity=1 ] (362.79,151.47) .. controls (362.79,149.26) and (364.58,147.47) .. (366.79,147.47) .. controls (368.99,147.47) and (370.79,149.26) .. (370.79,151.47) .. controls (370.79,153.68) and (368.99,155.47) .. (366.79,155.47) .. controls (364.58,155.47) and (362.79,153.68) .. (362.79,151.47) -- cycle ;
\draw [color={rgb, 255:red, 126; green, 211; blue, 33 }  ,draw opacity=1 ][line width=1.5]    (503.6,160.6) -- (514.02,146.92) -- (535.6,118.6) ;
\draw [color={rgb, 255:red, 209; green, 2; blue, 27 }  ,draw opacity=1 ][line width=1.5]    (521.6,152.6) -- (471.69,115.5) ;
\draw [color={rgb, 255:red, 189; green, 16; blue, 224 }  ,draw opacity=1 ][line width=1.5]    (451.6,159.6) -- (497.6,115.6) ;
\draw [color={rgb, 255:red, 245; green, 166; blue, 35 }  ,draw opacity=1 ][line width=1.5]    (473.08,158.96) -- (423.17,121.86) ;
\draw    (103.13,237.17) -- (153.8,237.17) ;
\draw [shift={(153.8,237.17)}, rotate = 0] [color={rgb, 255:red, 0; green, 0; blue, 0 }  ][fill={rgb, 255:red, 0; green, 0; blue, 0 }  ][line width=0.75]      (0, 0) circle [x radius= 3.35, y radius= 3.35]   ;
\draw [shift={(103.13,237.17)}, rotate = 0] [color={rgb, 255:red, 0; green, 0; blue, 0 }  ][fill={rgb, 255:red, 0; green, 0; blue, 0 }  ][line width=0.75]      (0, 0) circle [x radius= 3.35, y radius= 3.35]   ;
\draw    (153.8,237.17) -- (203.8,237.84) ;
\draw [shift={(203.8,237.84)}, rotate = 0.76] [color={rgb, 255:red, 0; green, 0; blue, 0 }  ][fill={rgb, 255:red, 0; green, 0; blue, 0 }  ][line width=0.75]      (0, 0) circle [x radius= 3.35, y radius= 3.35]   ;
\draw    (577.7,238.56) -- (630.24,238.72) ;
\draw [shift={(630.24,238.72)}, rotate = 0.17] [color={rgb, 255:red, 0; green, 0; blue, 0 }  ][fill={rgb, 255:red, 0; green, 0; blue, 0 }  ][line width=0.75]      (0, 0) circle [x radius= 3.35, y radius= 3.35]   ;
\draw    (52.47,237.17) -- (103.13,237.17) ;
\draw [shift={(103.13,237.17)}, rotate = 0] [color={rgb, 255:red, 0; green, 0; blue, 0 }  ][fill={rgb, 255:red, 0; green, 0; blue, 0 }  ][line width=0.75]      (0, 0) circle [x radius= 3.35, y radius= 3.35]   ;
\draw [shift={(52.47,237.17)}, rotate = 0] [color={rgb, 255:red, 0; green, 0; blue, 0 }  ][fill={rgb, 255:red, 0; green, 0; blue, 0 }  ][line width=0.75]      (0, 0) circle [x radius= 3.35, y radius= 3.35]   ;
\draw    (453.7,238.06) -- (515.7,239.06) ;
\draw [shift={(515.7,239.06)}, rotate = 0.92] [color={rgb, 255:red, 0; green, 0; blue, 0 }  ][fill={rgb, 255:red, 0; green, 0; blue, 0 }  ][line width=0.75]      (0, 0) circle [x radius= 3.35, y radius= 3.35]   ;
\draw [shift={(453.7,238.06)}, rotate = 0.92] [color={rgb, 255:red, 0; green, 0; blue, 0 }  ][fill={rgb, 255:red, 0; green, 0; blue, 0 }  ][line width=0.75]      (0, 0) circle [x radius= 3.35, y radius= 3.35]   ;
\draw    (515.7,239.06) -- (577.7,238.56) ;
\draw [shift={(577.7,238.56)}, rotate = 359.54] [color={rgb, 255:red, 0; green, 0; blue, 0 }  ][fill={rgb, 255:red, 0; green, 0; blue, 0 }  ][line width=0.75]      (0, 0) circle [x radius= 3.35, y radius= 3.35]   ;
\draw    (306.33,237.17) -- (357,237.17) ;
\draw [shift={(357,237.17)}, rotate = 0] [color={rgb, 255:red, 0; green, 0; blue, 0 }  ][fill={rgb, 255:red, 0; green, 0; blue, 0 }  ][line width=0.75]      (0, 0) circle [x radius= 3.35, y radius= 3.35]   ;
\draw [shift={(306.33,237.17)}, rotate = 0] [color={rgb, 255:red, 0; green, 0; blue, 0 }  ][fill={rgb, 255:red, 0; green, 0; blue, 0 }  ][line width=0.75]      (0, 0) circle [x radius= 3.35, y radius= 3.35]   ;
\draw    (357,237.17) -- (407,237.84) ;
\draw [shift={(407,237.84)}, rotate = 0.76] [color={rgb, 255:red, 0; green, 0; blue, 0 }  ][fill={rgb, 255:red, 0; green, 0; blue, 0 }  ][line width=0.75]      (0, 0) circle [x radius= 3.35, y radius= 3.35]   ;
\draw    (255.67,237.17) -- (306.33,237.17) ;
\draw [shift={(306.33,237.17)}, rotate = 0] [color={rgb, 255:red, 0; green, 0; blue, 0 }  ][fill={rgb, 255:red, 0; green, 0; blue, 0 }  ][line width=0.75]      (0, 0) circle [x radius= 3.35, y radius= 3.35]   ;
\draw [shift={(255.67,237.17)}, rotate = 0] [color={rgb, 255:red, 0; green, 0; blue, 0 }  ][fill={rgb, 255:red, 0; green, 0; blue, 0 }  ][line width=0.75]      (0, 0) circle [x radius= 3.35, y radius= 3.35]   ;

\draw (65.2,59.8) node [anchor=north west][inner sep=0.75pt]  [color={rgb, 255:red, 74; green, 144; blue, 226 }  ,opacity=1 ]  {$A_{n}$};
\draw (187.2,69.8) node [anchor=north west][inner sep=0.75pt]  [color={rgb, 255:red, 209; green, 2; blue, 27 }  ,opacity=1 ]  {$A_{n-2}$};
\draw (343.8,68.6) node [anchor=north west][inner sep=0.75pt]  [color={rgb, 255:red, 189; green, 16; blue, 224 }  ,opacity=1 ]  {$A_{n-4}$};
\draw (530.2,69.53) node [anchor=north west][inner sep=0.75pt]  [color={rgb, 255:red, 245; green, 166; blue, 35 }  ,opacity=1 ]  {$A_{n-6}$};
\draw (114.4,148.2) node [anchor=north west][inner sep=0.75pt]    {$\cdots $};
\draw (305.6,148.8) node [anchor=north west][inner sep=0.75pt]    {$\cdots $};
\draw (359,161.73) node [anchor=north west][inner sep=0.75pt]  [color={rgb, 255:red, 229; green, 29; blue, 235 }  ,opacity=1 ]  {$A_{0}$};
\draw (403.6,148.8) node [anchor=north west][inner sep=0.75pt]    {$\cdots $};
\draw (362.27,6.13) node [anchor=north west][inner sep=0.75pt]  [color={rgb, 255:red, 209; green, 2; blue, 27 }  ,opacity=1 ]  {$E'_{2}$};
\draw (492,8.4) node [anchor=north west][inner sep=0.75pt]  [color={rgb, 255:red, 189; green, 16; blue, 224 }  ,opacity=1 ]  {$E_{3}$};
\draw (300,103) node [anchor=north west][inner sep=0.75pt]  [color={rgb, 255:red, 245; green, 166; blue, 35 }  ,opacity=1 ]  {$E_{4}$};
\draw (146.67,8.93) node [anchor=north west][inner sep=0.75pt]  [color={rgb, 255:red, 126; green, 211; blue, 33 }  ,opacity=1 ]  {$E_{1}$};
\draw (241.07,7.4) node [anchor=north west][inner sep=0.75pt]  [color={rgb, 255:red, 126; green, 211; blue, 33 }  ,opacity=1 ]  {$E'_{1}$};
\draw (322.27,7.73) node [anchor=north west][inner sep=0.75pt]  [color={rgb, 255:red, 209; green, 2; blue, 27 }  ,opacity=1 ]  {$E_{2}$};
\draw (285.07,11.4) node [anchor=north west][inner sep=0.75pt]  [color={rgb, 255:red, 126; green, 211; blue, 33 }  ,opacity=1 ]  {$E_{1}$};
\draw (399.47,5) node [anchor=north west][inner sep=0.75pt]  [color={rgb, 255:red, 126; green, 211; blue, 33 }  ,opacity=1 ]  {$E'_{1}$};
\draw (455.47,8.2) node [anchor=north west][inner sep=0.75pt]  [color={rgb, 255:red, 126; green, 211; blue, 33 }  ,opacity=1 ]  {$E_{1}$};
\draw (617.87,8.2) node [anchor=north west][inner sep=0.75pt]  [color={rgb, 255:red, 126; green, 211; blue, 33 }  ,opacity=1 ]  {$E'_{1}$};
\draw (178.67,106.6) node [anchor=north west][inner sep=0.75pt]  [color={rgb, 255:red, 126; green, 211; blue, 33 }  ,opacity=1 ]  {$E_{1}$};
\draw (540.27,97.8) node [anchor=north west][inner sep=0.75pt]  [color={rgb, 255:red, 126; green, 211; blue, 33 }  ,opacity=1 ]  {$E'_{1}$};
\draw (519.07,4.2) node [anchor=north west][inner sep=0.75pt]  [color={rgb, 255:red, 209; green, 2; blue, 27 }  ,opacity=1 ]  {$E_{2}$};
\draw (549.47,4.2) node [anchor=north west][inner sep=0.75pt]  [color={rgb, 255:red, 209; green, 2; blue, 27 }  ,opacity=1 ]  {$E'_{2}$};
\draw (255.07,99.4) node [anchor=north west][inner sep=0.75pt]  [color={rgb, 255:red, 209; green, 2; blue, 27 }  ,opacity=1 ]  {$E_{2}$};
\draw (453.47,95.4) node [anchor=north west][inner sep=0.75pt]  [color={rgb, 255:red, 209; green, 2; blue, 27 }  ,opacity=1 ]  {$E'_{2}$};
\draw (95.87,207.4) node [anchor=north west][inner sep=0.75pt]  [color={rgb, 255:red, 209; green, 2; blue, 27 }  ,opacity=1 ]  {$E_{2}$};
\draw (568.67,208.2) node [anchor=north west][inner sep=0.75pt]  [color={rgb, 255:red, 209; green, 2; blue, 27 }  ,opacity=1 ]  {$E'_{2}$};
\draw (46.67,206.6) node [anchor=north west][inner sep=0.75pt]  [color={rgb, 255:red, 126; green, 211; blue, 33 }  ,opacity=1 ]  {$E_{1}$};
\draw (619.47,210.6) node [anchor=north west][inner sep=0.75pt]  [color={rgb, 255:red, 126; green, 211; blue, 33 }  ,opacity=1 ]  {$E'_{1}$};
\draw (196,210.2) node [anchor=north west][inner sep=0.75pt]  [color={rgb, 255:red, 245; green, 166; blue, 35 }  ,opacity=1 ]  {$E_{4}$};
\draw (409.6,95.8) node [anchor=north west][inner sep=0.75pt]  [color={rgb, 255:red, 245; green, 166; blue, 35 }  ,opacity=1 ]  {$E'_{4}$};
\draw (441.6,208.6) node [anchor=north west][inner sep=0.75pt]  [color={rgb, 255:red, 245; green, 166; blue, 35 }  ,opacity=1 ]  {$E'_{4}$};
\draw (584,1.8) node [anchor=north west][inner sep=0.75pt]  [color={rgb, 255:red, 189; green, 16; blue, 224 }  ,opacity=1 ]  {$E'_{3}$};
\draw (221.6,102.6) node [anchor=north west][inner sep=0.75pt]  [color={rgb, 255:red, 189; green, 16; blue, 224 }  ,opacity=1 ]  {$E_{3}$};
\draw (495.2,89) node [anchor=north west][inner sep=0.75pt]  [color={rgb, 255:red, 189; green, 16; blue, 224 }  ,opacity=1 ]  {$E'_{3}$};
\draw (143.2,208.2) node [anchor=north west][inner sep=0.75pt]  [color={rgb, 255:red, 189; green, 16; blue, 224 }  ,opacity=1 ]  {$E_{3}$};
\draw (508,212.2) node [anchor=north west][inner sep=0.75pt]  [color={rgb, 255:red, 189; green, 16; blue, 224 }  ,opacity=1 ]  {$E'_{3}$};
\draw (219.6,230) node [anchor=north west][inner sep=0.75pt]    {$\cdots $};
\draw (418,229.2) node [anchor=north west][inner sep=0.75pt]    {$\cdots $};
\draw (321.8,100.6) node [anchor=north west][inner sep=0.75pt]  [color={rgb, 255:red, 74; green, 144; blue, 226 }  ,opacity=1 ]  {$E_{n/2}$};
\draw (370.6,100.6) node [anchor=north west][inner sep=0.75pt]  [color={rgb, 255:red, 74; green, 144; blue, 226 }  ,opacity=1 ]  {$E'_{n/2}$};
\draw (291.4,207) node [anchor=north west][inner sep=0.75pt]  [color={rgb, 255:red, 74; green, 144; blue, 226 }  ,opacity=1 ]  {$E_{n/2}$};
\draw (338.6,206.2) node [anchor=north west][inner sep=0.75pt]  [color={rgb, 255:red, 74; green, 144; blue, 226 }  ,opacity=1 ]  {$E'_{n/2}$};
\draw (19.87,179.8) node [anchor=north west][inner sep=0.75pt]  [color={rgb, 255:red, 0; green, 0; blue, 0 }  ,opacity=1 ]  {$A_{n}$};

\end{tikzpicture}
}  
\caption{Resolution of the singularity $A_n$ for $n$ even.}
\label{fig:res A_n}
\end{figure}

The next lemma says that we can resolve a Du Val singularity of type $A_n$ by blowing up rational curves $C \simeq \p^1$ passing through the singularity. In this process, each curve blowup weakens the singularity from $A_n$ to $A_{n-1}$

\begin{lem}\label{res A_n lines}
    A general Du Val singularity of type $A_n$ can be resolved by blowing up the singular point followed by $n-2$ blowups along components of the exceptional divisor of the previous blowups.
\end{lem}

\begin{lem}\label{res D_n} The resolution of a Du Val singularity of type $D_n$ can be reached by performing $2 \cdot \left\lceil \dfrac{n-1}{2} \right\rceil$ blowups at a cluster\footnote{By \textit{cluster} we mean a sequence of infinitely near points of the ambient space, which may not be themselves totally ordered regarding the partial order $\succcurlyeq$: $Q \succcurlyeq P$ if $P$ is the image of $Q$ under some finite sucession of blowdowns.} of infinitely near base points lying over the singular point.
\end{lem}

The last three lemmas together with the explicit resolution of the canonical singularities of type $E$ \cite[Appendix A]{alv} allow us to have the following table:

\newpage


\renewcommand{\arraystretch}{1.7} 
\begin{table}[htp]
\begin{center}
\begin{tabular}{|c|c|}
\hline
type of singularity & number of point blowups to resolve                 \\ 
\hline
$A_n$               & $\left\lceil \dfrac{n}{2} \right\rceil$            \\ \hline
$D_n$               & $ 2 \cdot \left\lceil \dfrac{n-1}{2} \right\rceil$ \\ \hline
$E_6$               & 4                                                  \\ \hline
$E_7$               & 7                                                  \\ \hline
$E_8$               & 8                                                  \\ \hline
\end{tabular}
\caption{Number of point blowups to resolve the Du Val singularities.} 
\label{number of pt blowup to resolve the du val sing}
\end{center}
\end{table}

\section{Volume-preserving Sarkisov links} \label{vp Sark links}

According to Definition \ref{defn vp Sark links}, volume-preserving Sarkisov links are Sarkisov links endowed with additional data and the volume-preserving property. 

This section aims to study volume-preserving Sarkisov links whose source is a Calabi-Yau pair $(\p^3,D)$ with $\coreg(\p^3,D)=2$. Recall that in this situation, by Lemma \ref{coreg 2}, $D$ is an irreducible normal quartic surface with canonical singularities. When $D$ is nonsingular or has a singularity of type $A_1$ or $A_2$, results were obtained in \cite{acm}. Our main result, which extends part of the work in \cite{acm}, is the following:

\begin{thm}\label{thm vp weights}
     Let $(\p^3,D)$ be a log Calabi-Yau pair of coregularity 2 and $\pi \colon (X,D_X) \rightarrow (\p^3,D)$ be a volume-preserving toric $(1,a,b)$-weighted blowup of a torus invariant point. Then this point is necessarily a singularity of $D$ and,  up to permutation, the only possibilities for the weights, depending on the type of singularities, are listed in the following Table 
    \ref{table vp weights III}.

\end{thm}


\renewcommand{\arraystretch}{1.0}
\begin{table}[htp]
\begin{center}
\begin{tabular}{|c|c|}
\hline
type of singularity & volume-preserving weights \\
\hline
$A_1$ & (1,1,1) \\
$A_2$ & (1,1,1), (1,1,2), \\ 
$A_3$ & (1,1,1), (1,1,2), \textcolor{red}{(1,1,3)} \\
$A_4$ & (1,1,1), (1,1,2), \textcolor{red}{(1,1,3)}, \textcolor{purple}{(1,2,3)}\\
$A_5$ & (1,1,1), (1,1,2), \textcolor{red}{(1,1,3)}, \textcolor{purple}{(1,2,3)}\\
$A_6 $ & (1,1,1), (1,1,2), \textcolor{red}{(1,1,3)}, \textcolor{purple}{(1,2,3)}, \textcolor{cyan}{(1,2,5)}, \textcolor{green}{(1,3,4)}\\
$A_7 $ & (1,1,1), (1,1,2), \textcolor{red}{(1,1,3)}, \textcolor{purple}{(1,2,3)}, \textcolor{cyan}{(1,2,5)}, \textcolor{green}{(1,3,4)}, \textcolor{brown}{(1,3,5)}\\
\hline
$D_4$ & (1,1,1), (1,1,2)\\
$D_5$ & (1,1,1), (1,1,2), \textcolor{purple}{(1,2,3)}\\
$D_6$ & (1,1,1), (1,1,2), \textcolor{purple}{(1,2,3)}\\
$D_7$ & (1,1,1), (1,1,2), \textcolor{purple}{(1,2,3)}, \textcolor{green}{(1,3,4)}, \textcolor{brown}{(1,3,5)}\\
$D_8$ & (1,1,1), (1,1,2), \textcolor{purple}{(1,2,3)}, \textcolor{green}{(1,3,4)}, \textcolor{brown}{(1,3,5)}\\
$D_9$ & (1,1,1), (1,1,2), \textcolor{purple}{(1,2,3)}, \textcolor{green}{(1,3,4)}, \textcolor{brown}{(1,3,5)}, \textcolor{orange}{(1,4,5)}\\
$D_{10}$ & (1,1,1), (1,1,2), \textcolor{purple}{(1,2,3)}, \textcolor{green}{(1,3,4)}, \textcolor{brown}{(1,3,5)}, \textcolor{orange}{(1,4,5)}\\
\hline
$E_6$ & (1,1,1), (1,1,2), \textcolor{purple}{(1,2,3)}, \textcolor{green}{(1,3,4)}, \textcolor{brown}{(1,3,5)} \\
$E_7$ & (1,1,1), (1,1,2), \textcolor{purple}{(1,2,3)}, \textcolor{green}{(1,3,4)}, \textcolor{brown}{(1,3,5)} \\
$E_8$ & (1,1,1), (1,1,2), \textcolor{purple}{(1,2,3)}, \textcolor{green}{(1,3,4)}, \textcolor{brown}{(1,3,5)}, \textcolor{orange}{(1,4,5)}\\
\hline
\end{tabular}
\caption{Table summarizing volume-preserving weights, up to permutation.} 
\label{table vp weights III}
\end{center}
\end{table}

Given a Mori fibered space $X/\Spec(\C)$, where $X$ is a Fano variety of Picard rank 1, the first step in the Sarkisov decomposition of any birational map with source $X/\Spec(\C)$ is a divisorial contraction $\pi \colon Y \rightarrow X$ which will initiate the first link. This first link is necessarily of type I or II because $\rho(X)=1$. 

A divisorial contraction from a terminal variety to a nonsingular Fano threefold falls into two cases depending on the dimension of the image of the prime exceptional divisor: it contracts such a divisor to either a curve or a point.

Both cases have been studied, and there are some classification results depending on the singularities involved and assumptions made. For the former, see \cite{ccc, cut,duc div ext, mori, tzi}. For the latter, a result by Kawakita \cite[Theorem 1.1]{kaw1} says that for suitable analytic coordinates at the point, this divisorial contraction can be described as the weighted blowup with weights $(1,a,b)$, where $\text{GCD}(a,b)=1$. We call this map a \textit{Kawakita blowup} of the point.

We point out that if one disregards the nonsingularity assumption on the threefold in the range, additional classes of divisorial contractions beyond these last two may arise. See \cite{tzi,zik}.

Now consider a reduced Weil divisor $D$ on $X$ such that $(X,D) \rightarrow \Spec(\C)$ has the structure of a Mori fibered Calabi-Yau pair. By the previous considerations and \cite[Proposition 3.1]{acm}, the possible volume-preserving divisorial contractions are 
\begin{center}
 $\pi=$
$\begin{cases} \hfil
 \text{divisorial contraction mapping the exceptional divisor onto a curve }C \subset D,\\
\text{Kawakita blowup of a point $P \in D$}
\end{cases}.$
\end{center}

We will see at a glance that such a blown up point must be a singularity of $D$ when $(X,D)$ is canonical. The following lemma deduced from \cite[Proposition 3.1]{acm} can be readily verified for the case where the divisorial contraction is an ordinary blowup by computing some discrepancies and comparing divisors.

\begin{lem}\label{vp center sing pt}
Let $\pi \colon (Y,D_Y) \rightarrow (X,D_X)$ be a volume-preserving terminal divisorial contraction between threefold Calabi-Yau pairs contracting a divisor $E \subset Y$ to a closed point $P \in X$. Assume $(X,D_X)$ has canonical singularities. Then $P$ is a singularity of $D_X$, and $D_Y$ is the strict transform of $D_X$ in $Y$.
\end{lem}

\begin{proof}
    By the first part of \cite[Proposition 3.1]{acm}, $P \in D_X$. Suppose $P$ is not a singularity of $D_X$. For surfaces, being terminal at a point is equivalent to being nonsingular at a point. So $D_X$ is terminal at $P$, and hence by \cite[Proposition 3.1]{acm} we must have $\codim_{X}P = 2$. But this is absurd, since $\codim_{X}P = 3$. Therefore $P \in \Sing(D_X)$. By \cite[Proposition 2.6]{acm} we get that $D_Y = \pi_*^{-1}D_X$.
\end{proof}

Our problem is the following:
\begin{center}
    \textit{Given the Du Val singularity at a point $P \in D$, to determine for which weights $(1,a,b)$ the Kawakita blowup $\pi \colon (X,\Tilde{D}) \rightarrow (\p^3,D)$ of $P$ with weights $(1,a,b)$ is volume-preserving.}
\end{center}

The following result by Duarte Guerreiro will allow us to restrict our possibilities for the weights if we are interested in those ones inducing Sarkisov links.

\begin{thm}[cf. \cite{gue} Proposition 1.2]\label{gue thm}
    Let $\varphi \colon X \rightarrow \mathbb{P}^3$ be the toric $(1,a,b)$-weighted blowup of a point. Then $\varphi$ initiates a (toric) Sarkisov link from $\mathbb{P}^3$ if and only if, up to permutation of $a$ and $b$, 
\[
(a,b) \in \{(1,1), (1,2), (2,3), (2,5) \}.
\]
\end{thm}

The following result is the volume-preserving version of Theorem \ref{gue thm} for the case where $\coreg(\p^3,D)=2$. Amending Table \ref{table vp weights III} accordingly to Theorem \ref{gue thm}, we have the list of volume-preserving Sarkisov links.

\begin{cor}\label{cor vp weights sark}
     Let $(\p^3,D)$ be a log Calabi-Yau pair of coregularity 2 and $\pi \colon (X,D_X) \rightarrow (\p^3,D)$ be a volume-preserving toric $(1,a,b)$-weighted blowup of a torus invariant point. Then this point is necessarily a singularity of $D$ and,  up to permutation, the only possibilities for the weights initiating a volume-preserving Sarkisov link, depending on the type of singularities, are listed in the following Table \ref{table vp weights IV}.

\end{cor}


\begin{table}[htp]
\begin{center}
\begin{tabular}{|c|c|}
\hline
type of singularity & volume-preserving weights \\
\hline
$A_1$ & (1,1,1) \\
$A_2$ & (1,1,1), (1,1,2) \\ 
$A_3$ & (1,1,1), (1,1,2) \\
$A_4$ & (1,1,1), (1,1,2),  \textcolor{purple}{(1,2,3)}\\
$A_5$ & (1,1,1), (1,1,2), \textcolor{purple}{(1,2,3)}\\
$A_{\ge 6} $ & (1,1,1), (1,1,2), \textcolor{purple}{(1,2,3)}, \textcolor{cyan}{(1,2,5)}\\
\hline
$D_4$ & (1,1,1), (1,1,2)\\
$D_{\ge 5}$ & (1,1,1), (1,1,2), \textcolor{purple}{(1,2,3)}\\
\hline
$E_6$ & (1,1,1), (1,1,2), \textcolor{purple}{(1,2,3)} \\
$E_7$ & (1,1,1), (1,1,2), \textcolor{purple}{(1,2,3)} \\
$E_8$ & (1,1,1), (1,1,2), \textcolor{purple}{(1,2,3)} \\
\hline
\end{tabular}
\caption{Table summarizing volume-preserving weights initiating Sarkisov links, up to permutation.} 
\label{table vp weights IV}
\end{center}
\end{table}


Since we will deal with (t,c) Calabi-Yau pairs, by \cite[Proposition 2.6]{acm} and \cite[Lemma 2.8]{acm}, the boundary divisors with respect to which we want the Sarkisov link to be volume-preserving are precisely the strict transforms of the initial one. 

Before continuing to deal with this problem, we check that our hypotheses are nonempty, i.e., that there exists an irreducible normal quartic surface having a canonical singularity of each possible type. The answer is given by Kato \& Naruki, Urabe and Yang \cite{kn,ur1,ur2,ur3,ya}, who analyzed all the possibilities of combinations of singularities in our context, as explained in Subsection \ref{exp res Du Val}. For the cases $A_n$ and $D_n$, for instance, we know that $n$ is at most $19$ and $18$, respectively. For the remaining case $E_n$, there also exists a quartic surface with this type of singularity.

By Proposition \ref{prop crep bir morp}, the volume-preserving property of the toric $(1,a,b)$-weighted blowup can be detected by the vanishing of $a(E,\p^3,D)$, the discrepancy of $E = \Exc(\pi)$ with respect to $(\p^3,D)$.

Concerning the definition and main ideas involving weighted projective spaces and weighted blowups, we refer the reader to \cite{cks,cls,ful}. They can be realized as toric varieties and geometric quotients. This realization will be crucial in our work.

We usually ask that the $(n+1)$-tuple of positive integers defining the weights consists of coprime elements.

\paragraph{First approach to the problem:} Let $\sigma=(a_1,\ldots,a_n) \in \ZZ_{>0}^n$ and consider $\pi \colon X \rightarrow \p^n$ the $\sigma$-weighted blowup of $\p^n$ at the point $P = (1 : 0 : \cdots : 0)$. Set $E \coloneqq \Exc(\pi)$.

We define the $\sigma$-weight of a monomial $M = x_1^{p_1}\cdots x_n^{p_n}$ as $\wt_{\sigma}(M)=p_1a_1+ \cdots + p_na_n$. Given a nonzero homogeneous polynomial $f \in \C[x_0,\ldots,x_n]$, write its dehomogenization with respect to the variable $x_0$ as
\begin{center}
    $f_* = \displaystyle \sum_I \alpha_IM_I \in \C[x_1,\ldots,x_n]$,
\end{center}
where $M_I$ runs through the monomials in $x_1,\ldots,x_n$.

The $\sigma$-weight of $f$ is defined by 
\begin{center}
    $\wt_{\sigma}(f) \coloneqq \min \{\wt_{\sigma}(M_I) | \alpha_I \neq 0\}$.
\end{center}

Let $D = V(f) \subset \p^n$ be a hypersurface. Set $\wt_{\sigma}(D) \coloneqq \wt_{\sigma}(f)$. 
By \cite[Lemma 3.5]{hay} we have \begin{align*}
   \Tilde{D} & = \pi^*D - (\wt_{\sigma}(D))E,\ \text{and}\\
    K_X & = \pi^*K_{\p^n}+(a_1 + \cdots + a_n - 1) E \ .
\end{align*}

Assume $1 \leq a < b$ in the context of our problem. We point out the $a$ and $ b$ are coprime for the weighted blowup to be in the Mori category of terminal varieties.  Thus one has $a(E,\p^3,D)=a+b - \text{wt}_{\sigma}(D)$. The $(1,a,b)$-weighted blowup is volume-preserving if and only $a(E,\p^3,D)=0$, that is, we must have $a+b = \text{wt}_{\sigma}(D)$. 

Thus, one can proceed with a straightforward analysis of the homogeneous equation $f$ that defines $D$. This can be done by verifying a lot of conditions on $f$ which imply that certain monomials will not appear. 
However, this analysis will be very tedious if we make the singularity worse. Moreover, this strategy depends on the coordinates chosen. \\

Recalling that our problem is to determine ``volume-preserving weights'', our task is to check when $a(E,\p^3,D)=0$. Realizing this discrepancy by a divisorial contraction $\pi$ is relatively feasible because weighted blowups can be described in terms of charts and they locally look like affine spaces up to an unramified cover. However, this would require a lot of discrepancy computations.

Using the property that discrepancies only depend on the valuations associated to divisors, we can have an alternative way to perform the task of finding the volume-preserving weights. There exists a realization of the divisorial valuation $\nu_E$ associated to $E$ based on a toric description of the weighted blowup. This was the same strategy adopted by Araujo, Corti \& Massarenti and Kawakita \cite{acm,kaw1}, which has the advantage of being coordinate-free.

\paragraph{Toric description of the weighted blowup:} As a toric variety, $X$ is determined by the fan $\Sigma$ in $\RR^3$ given by the union of all possible 3-dimensional cones and its subcones generated by the vectors in $\{v_0,v_i,v_j\}$ and in $\{v_i,v_j,v\}$ for $i,j \in \{1,2,3\}$, where $v_0=(-1,-1,-1)$, $v_i = e_i$ for $i \in \{1,2,3\}$, and $v=(1,a,b)$.


The ``toric description of the weighted blowup'' consists of a finite sequence of ordinary blowups with nonsingular centers such that the valuation on $\p^3$ corresponding to the last exceptional divisor coincides with the valuation associated to $E$ on $\p^3$. 

This corresponds to the toric realization of a resolution of singularities of the weighted blowup by adding rays in order to make its corresponding fan $\Sigma$ smooth (or regular), that is, in such a way that the minimal generators of all its subcones form part of a $\mathbb{Z}$-basis of the corresponding lattice $N$. See \cite[Theorem 11.1.9]{cls}. 
It is straightforward to check that the set of minimal generators $\{v_1,v_2,v\}$ is not a $\ZZ$-basis for the lattice $\ZZ^3$ unless $(a,b)=(1,1)$. So the first blowup corresponds to adding the ray $(1,1,1)$, then $(1,2,2)$, and so on until adding the ray $(1,a,a)$. Then we add the rays $(1,a,a+1),\ldots,(1,a,b)$. This process of adding rays ends with a regular fan $\Sigma'$ that is a refinement of $\Sigma$, obtained by a star subdivision of $\Sigma$ such that the toric morphism $X_{\Sigma'} \rightarrow X_{\Sigma}$ is a projective resolution of singularities. See \cite[Definition 3.3.13 \& Section 11.1]{cls}. 

Let $\Sigma_i$ be the fan in $\RR^3$ corresponding to the toric variety $X_i$ at the $i$-th step of this process. According to the toric description, $\Sigma_i$ is obtained by a star subdivision of $\Sigma_{i-1}$ along a cone $\sigma_{i-1}$, that is,
\begin{center}
    $\Sigma_i \coloneqq \Sigma_{i-1}^*(\sigma_{i-1})$.
\end{center}

The induced toric morphism $X_i \rightarrow X_{i-1}$ corresponds to blowing up the orbit of the torus action on $X_{i-1}$ associated to the cone $\sigma_{i-1}$ by means of the orbit-cone correspondence theorem. See \cite[Theorem 3.2.6]{cls}.

Inductively we have 
\begin{center}
$X_b \xrightarrow{\pi_b} X_{b-1} \xrightarrow{\pi_{b-1}} \cdots \xrightarrow{\pi_2} X_1 \xrightarrow{\pi_1} X_0 \coloneqq \p^3,$
\end{center} 

\nin where $X_i \xrightarrow{\pi_i} X_{i-1}$ is the blowup of the center $z_{i-1} = z_E X_{i-1}$ of the valuation $E$ on $X_{i-1}$. Notice that $X_1 \to X_0$ is the blow-up of $z_0=P$. For every $i$, we denote by $E_i\subset X_i$ the exceptional divisor of $\pi_i$, and for $j>i$ we denote by $E^j_i\subset X_j$ the strict transform of $E_i$ in $X_j$. The
following key properties follow directly from the toric description of the weighted blowup:
\begin{enumerate}
\item\label{key prop 1 tor desc} For all $0\leq j< a$, the center $z_{j}$ is a closed point of $X_{j}$. If $j\geq 1$, then 
$z_{j}\in E_{j}\subset X_{j}$, and if $j\geq 2$, then  
\[
z_{j}\in E_{j}\setminus E^{j}_{j-1}.
\]
\item\label{key prop 2 tor desc} The center $z_a\in X_a$ is the generic point of a line
  $L_a\subset E_a \simeq \p^2$. If $a\geq 2$, then
\[
L_a\not \subset E^a_{a-1}.
\]
\item\label{key prop 3 tor desc} For all $a+1\leq j < b$, the center $z_{j}\in X_{j}$ is a
  section 
\[
L_{j}\subset E_{j}\setminus E^{j}_{j-1}
\] 
of the projection $E_{j}\to L_{j-1}$.
\item $E_b=E$ (by this we mean that the exceptional divisors $E_b$ and $E$ induce the same valuation on $X_0 = \p^3$).
\end{enumerate}

Set $\pi' \coloneqq \pi_1 \circ \cdots \circ \pi_b$. We have the following 
\[\begin{tikzcd}
	X && {X_b} \\
	{\mathbb{P}^3}
	\arrow["\pi"', from=1-1, to=2-1]
	\arrow["{\pi_1 \circ \cdots \circ \pi_b}", from=1-3, to=2-1]
\end{tikzcd}. \]

Thus $E$ and $E_b = \Exc(\pi_b)$ induce the same valuation on $K(\p^3) \simeq K(X) \simeq K(X_b)$. Roughly speaking, $a = \sharp \{\text{blown up points}\}$  and $b-a = \sharp \{\text{blown up curves}\}$ in the sequence. 

We remark that the exceptional divisors $\Exc(\pi_i)$, for $i \geq a+1$ are Hirzebruch surfaces, since $\Exc(\pi_i) \simeq \p(\mathcal{N}_{L_{i-1}/X_{i-1}}^{\vee})$ and $\mathcal{N}_{L_{i-1}/X_{i-1}}^{\vee}$ is a rank 2 vector bundle over $L_{i-1} \simeq \p^1$. In fact, $z_a \simeq \p^1$ and so are all the remaining centers given by sections $z_j$. Up to isomorphism, a rank 2 vector bundle over $\p^1$ is of the form $\oo_{\p^1} \oplus \oo_{\p^1}(n)$ for some $n \in \ZZ_{\geq 0}$. See \cite[Proposition III.15 (i)]{bea} or face this as a particular case of the famous Grothendieck Theorem about finite rank vector bundles over $\p^1$.\\

The following key lemma will be essential because it will give us bounds for the possibilities for $a$ and $b$, besides the relation $a+b=\text{wt}_{\sigma}(D)$. 

\begin{lem}[Key lemma]\label{key lemma} Let $D \subset \p^3$ be an irreducible normal quartic surface with an isolated canonical singularity at $P \in \p^3$. The weighted blowup at $P$ with coprime weights $(1,a,b)$ is volume-preserving if and only if each blowup in its toric description is volume-preserving.
\end{lem}

\begin{proof} Let $\pi \colon X \rightarrow \p^3$ be the $(1,a,b)$-weighted blowup of $P$. Consider the chain of blowups which realizes the valuation associated to $E = \Exc(\pi)$, as explained above:
\begin{center}
$X_b \xrightarrow{\pi_b} X_{b-1} \xrightarrow{\pi_{b-1}} \cdots \xrightarrow{\pi_2} X_1 \xrightarrow{\pi_1} X_0 \eqqcolon \p^3 $.
\end{center} 

Denote by $D_i$ the strict transform of $D$ on $X_i$.\\

\nin $(\Rightarrow)$ Let us show by a contrapositive argument.
Suppose some blowup in the toric description is not volume-preserving, and define
\begin{center}
    $i_0 \coloneqq \max\{ i \mid \pi_1, \ldots,\pi_i \ \text{are all volume-preserving} \}$.
\end{center}

Thus, $\pi_{i_0+1}$ is the first blowup which is not volume-preserving. Since the pair $(\p^3,D)$ is canonical, the corresponding discrepancies are always nonnegative.

By Proposition \ref{prop crep bir morp}, such discrepancies are zero if and only if the corresponding blowup is volume-preserving. This implies that $a(E_{i_0+1},\p^3,D) > 0$. We have two possible cases:\\

\nin \textbf{Case 1:} $z_{i_0} \not\subset D_{i_0} $.\\
It follows from the toric description that the next centers of the blowups will no longer belong to the strict transforms of $D$, that is, $z_j \not\subset D_j$ for all $j \geq i_0 + 1$. Consequently we will have $a(E_j,\p^3,D) > 0$ for all $j \geq i_0 + 1$.

Thus, $a(E_b,\p^3,D) > 0$ and therefore $\pi$ is not volume-preserving.\\

\nin \textbf{Case 2:} $z_{i_0} \subset D_{i_0}$.\\
In this case, $z_{i_0}$ is either a curve or a closed point. If $z_{i_0}$ is a curve, then $\pi_{i_0+1} \colon X_{i_0+1} \rightarrow X_{i_0}$ is volume-preserving, since one has $a(E_{i_0+1},X_{i_0},D_{i_0}) = 0$.

If $z_{i_0}$ is a closed point, then either $z_{i_0} \in \Sing(D_{i_0})$ or $z_{i_0} \notin \Sing(D_{i_0})$.

In the first scenario, since canonicity is preserved under volume-preserving maps by \cite[Lemma 2.8]{acm}, it follows that $z_{i_0}$ is a strict canonical singularity of $D_{i_0}$. Hence, one also has $a(E_{i_0+1},X_{i_0},D_{i_0}) = 0$.

So we are inevitably in the second scenario $z_{i_0} \notin \Sing(D_{i_0})$. By the toric description, we observe that necessarily $i_0 < a$. It follows that $z_i=P_i \in D_i$ is a nonsingular point for every $i$ such that $i_0 \leq i \leq a-1$. Indeed, if $z_i \not\subset D_i$ for some $i$, arguing analogously as in Case 1, we deduce that $\pi$ is not volume-preserving. 

Therefore, $a(E_a,X_{a-1},D_{a-1}) > 0$.

Next, we analyze what happens for $i \geq a$. 
Recall that we have $z_a = D_a \cap E_a \subset D_a$ and $\pi_{a+1} \colon X_{a+1} \rightarrow X_a$ is volume-preserving.

We have the following:
\[\begin{tikzcd}
	{X_{a-1}} & {X_a} & {X_{a+1}} & \cdots & {X_b}
	\arrow["{\text{Bl}_{P_{a-1}}}"', from=1-2, to=1-1]
	\arrow["{\text{Bl}_{L_a}}"', from=1-3, to=1-2]
	\arrow["{\text{Bl}_{L_{a+1}}}"', from=1-4, to=1-3]
	\arrow["{\text{Bl}_{L_{b-1}}}"', from=1-5, to=1-4]
\end{tikzcd}.\]

We have that
\[\begin{tikzcd}
	{D_{a-1}} & {D_a+E_a} & {D_{a+1}+E_a^{a+1}+2E_{a+1}} \\
	&& {D_{a+2}+E_a^{a+2}+2E_{a+1}^{a+2}+3E_{a+2}} \\
	&& \vdots \\
	&& {D_b+E_a^b+2E_{a+1}^b+3E_{a+2}^b+\cdots+(b-a-1)E_{b-1}^b+(b-a)E_b}
	\arrow["{\pi_a^*}", from=1-1, to=1-2]
	\arrow["{\pi_{a+1}^*}", from=1-2, to=1-3]
	\arrow["{\pi_{a+2}^*}", from=1-3, to=2-3]
	\arrow["{\pi_{a+3}^*}", from=2-3, to=3-3]
	\arrow["{\pi_b^*}", from=3-3, to=4-3]
\end{tikzcd},\]
which implies
\begin{align*}
   D_a & = \pi_a^*D_{a-1}-E_a,\\
    D_{a+1} & = \pi_{a+1}^*D_a - E_{a+1} \\
    & = \pi_{a+1}^*(\pi_a^*D_{a-1}-E_a) - E_{a+1} \\
    & = (\pi_a \circ \pi_{a+1})^*D_{a-1} - E_a^{a+1} - 2E_{a+1}, \\
 D_{a+2} & = (\pi_a \circ \pi_{a+1} \circ \pi_{a+2})^*D_{a-1} - E_a^{a+2} - 2E_{a+1}^{a+2} - 3E_{a+2}, \\
 &\  \vdots \\
 D_b & = (\pi_a \circ \cdots \circ \pi_b)^*D_{a-1}-E_a^b-2E_{a+1}^b-3E_{a+2}^b-\cdots-(b-a-1)E_{b-1}^b-(b-a)E_b .
\end{align*}

By the adjunction formula, we have that
\begin{align*}
   K_{X_a} & = \pi_a^*K_{X_{a-1}}+2E_a,\\
    K_{X_{a+1}} & = \pi_{a+1}^*(K_{X_a}) + E_{a+1} \\
    & = \pi_{a+1}^*(\pi_a^*K_{X_{a-1}}+2E_a) + E_{a+1} \\
    & = (\pi_a \circ \pi_{a+1})^*K_{X_{a-1}}+2E_a^{a+1} + 3E_{a+1},\\
 K_{X_{a+2}} & = (\pi_a \circ \pi_{a+1} \circ \pi_{a+2})^*K_{X_{a-1}} + 2E_a^{a+2} + 3E_{a+1}^{a+2} + 4E_{a+2}^{a+2}, \\
 &\  \vdots \\
 K_{X_b} & = (\pi_a \circ \cdots \circ \pi_b)^*K_{X_{a-1}}+2E_a^b+3E_{a+1}^b+4E_{a+2}^b+\cdots+(b-a)E_{b-1}^b+(b-a+1)E_b.
\end{align*}

Therefore,
\begin{center}
    $K_{X_b} + D_b = (\pi_a \circ \cdots \circ \pi_b)^*(K_{X_{a-1}}+D_{a-1})+E_a^b+E_{a+1}^b+E_{a+2}^b+\cdots+E_{b-1}^b+E_b $.
\end{center}

It follows that $a(E_b,X_{a-1},D_{a-1}) = 1$, whence we deduce $a(E_b,\p^3,D) = 1 > 0$. Consequently $\pi$ is not volume-preserving.\\

\nin $(\Rightarrow)$ Since all the blowups are volume-preserving, the corresponding discrepancies 
\begin{center}
    $a(E_1,\p^3,D), \ldots, a(E_b,X_{b-1},D_{b-1})$
\end{center}
will be zero as a consequence of the fact that composition of volume-preserving maps is also volume-preserving. See Remark \ref{comp of vp is vp}. 

Thus, $a(E,\p^3,D)= a(E_b,\p^3,D) = 0$ and therefore $\pi$ is volume-preserving.

\end{proof}

Our Lemma \ref{key lemma} indicates that after realizing the divisorial valuation $\nu_E$ associated to $E$ by means of a sequence of ordinary blowups at points or curves through a toric description, we must verify that all these intermediate blowups are volume-preserving.

Concerning the $A_n$ case of the problem, Lemma \ref{key lemma} restricts the possibilities for the weights. 

\begin{cor}\label{bound a for A_n}
For $D \subset \p^3$ a quartic surface with an isolated canonical singularity of type $A_n$, a sequence of $1 \leq m \leq \left\lceil \dfrac{n}{2} \right\rceil$ ordinary blowups at the singular points of the strict transforms of $D$ is volume-preserving. On the other hand, the weights of the form $(1, a, b)$ with $a > \left\lceil \dfrac{n}{2} \right\rceil$ are not volume-preserving.
\end{cor}

\begin{proof}
The first part is direct from the Lemma \ref{res A_n}, which tells us that $D$ is resolved after $\left\lceil \dfrac{n}{2} \right\rceil$ blowups. If $a > \left\lceil \dfrac{n}{2} \right\rceil$, by the toric description this means that we must blowup points outside the singular locus of the strict transform of $D$. These blowups are not volume-preserving, and hence by Lemma \ref{key lemma} the assertion follows.
\end{proof}

For the $D$-$E$ case, analogous results are harder to be stated in a precise way. This is because the graph of the exceptional divisor of a minimal resolution is more complicated. It has a vertex with the degree $\geq 2$ and more than 2 leaves. The \textit{degree of a vertex} of a graph is the number of edges leaving from it and a \textit{leaf} is a vertex of degree 1. In any case, $a$ is also certainly bounded by the number of point blowups needed to resolve such singularities. Therefore, by Lemma \ref{res D_n} and the discussion in the Subsection \ref{exp res Du Val} we must have
\begin{center}
    $a \leq 2 \cdot \left\lceil \dfrac{n-1}{2} \right\rceil, 4, 7, 8 $,
\end{center}
for the $D_n,E_6,E_7$ and $E_8$ cases, respectively.

\section{Determination of the volume-preserving weights} \label{deter vp weights}

With the results of the previous sections, we are ready to determine the desired volume-preserving weights. 
The content of this section is a classification of all possible toric weighted blowups $\pi \colon X \rightarrow \p^3$ that are volume-preserving for a Calabi-Yau pair $(\p^3,D)$ with $\coreg(\p^3,D)=2$, depending on the singularities on $D$. This classification will culminate, thanks to Theorem \ref{gue thm}, in a further classification, whose weights will drive volume-preserving Sarkisov links originating from these Calabi-Yau pairs.

In what follows we will determine which weights will be volume-preserving and not necessarily the candidates to initiate a volume-preserving Sarkisov link. The answer will depend on the type of canonical singularity of the quartic surface.

Even if some volume-preserving weights do not yield a volume-preserving link, they will be relevant to producing a model of the quartic surface embedded in a toric variety or links of a log version of the Sarkisov program. See Remark \ref{rmk (1,1,3)}. 

Due to the difficulty level of this problem in terms of the number of computations needed, we will restrict ourselves to singularities $A_n$ with $1 \leq n \leq 7$ and $D_n$ with $4 \leq n \leq 10$. All the singularities of type E will be contemplated. This will be enough for the second classification to cover all types of strict canonical singularities. 

\begin{rmk}
    Throughout the proof of Theorem \ref{thm vp weights}, we only consider global toric weighted blowups. In this context, we do not describe all divisorial contractions to $\p^3$ which are in local analytic coordinates a $(1,a,b)$-weighted blowup according to \cite[Theorem 1]{kaw1}, much less all volume-preserving ones to Calabi-Yau pairs $(\p^3,D)$. Indeed, the weighted blowup does depend on the local coordinates chosen (algebraic x analytic). 
\end{rmk}

\begin{rmk}\label{rmk vp weights for A_m, m geq n}
By the toric description of the weighted blowup, if the weights $(1,a,b)$ are volume-preserving for an $A_n$ or $D_n$ singularity, then $(1,a,b)$ are volume-preserving for an $A_m$ or $D_m$ singularity with $m \geq n$.  
\end{rmk}

\subsection{Notation for $D \subset \p^3$ quartic surface and some generalities}\label{notation quartic}

Let $D \subset \p^3$ be an irreducible normal quartic surface having a strict canonical singularity at $P$ so that we have $m_P(D)=2$. Fix homogeneous coordinates such that $P=(1:0:0:0)$ and the equation of $D$ has the form 
\begin{center}
    $x_0^2A+x_0B+C = 0$,
\end{center}

\nin where $A,B,C \in \C[x_1,x_2,x_3]$ are homogeneous polynomials of degrees 2, 3 and 4, respectively. Notice that the tangent cone of $D$ at $P$ is given by $TC_PD = \{A = 0\}$.

Given a positive integer $n$, we will say that \textit{$P$ is of type $A_{\geq n}$} to express that $P$ is a singularity of type $A_m$ with $m \geq n$. Analogously for the $D_n$ case. Sometimes, instead of saying ``$P$ is a singularity of type $A_n$'', we will abbreviate it to ``$P$ is $A_n$''. 

From the local description of the Du Val singularities, we see that:
\begin{equation*}
    \left.\begin{aligned}
  P\  \text{is}\  A_1 \Leftrightarrow \rk(A) = 3 \\
  P\  \text{is}\  A_{\geq 2} \Leftrightarrow \rk(A) = 2\\
  \end{aligned}\right\} \Rightarrow P\  \text{is}\  A_n \Leftrightarrow \rk(A) > 1 \ \text{and}
\end{equation*}
\begin{center}
    $P\  \text{is}\  D\text{-}E \Leftrightarrow \rk(A) = 1$.
\end{center}

Notice that the projectivization of $TC_PD$ is a plane conic, not necessarily irreducible. The ranks 3, 2 and 1 correspond to an irreducible conic, two concurrent lines and a double line (nonreduced), respectively.  

Therefore, we have an immediate criterion to detect an $A_1$ singularity and necessary conditions for the other types.

Possibly after changing homogeneous coordinates, if $\rk(A) =2$ or 1, we may assume that $A=x_2x_3$ or $A=x_3^2$, respectively. 

We will adopt the following notation which will ease our computations. Write 
\begin{align*}
    B(x_1,x_2,x_3) & = \displaystyle\sum_{i=0}^3 b_ix_1^{3-i},\ \text{where}\ b_i \in \C[x_2,x_3]_i\ \text{and} \\
    C(x_1,x_2,x_3) & = \displaystyle\sum_{i=0}^4 c_ix_1^{4-i},\ \text{where}\ c_i \in \C[x_2,x_3]_i .
\end{align*}

Also, write the following, where all the Greek letters indicating coefficients are complex numbers:
\begin{itemize}
    \item $b_1 = \beta_2x_2+\beta_3x_3$;
    \item $b_2 = \rho_2x_2^2 + \rho_{23}x_2x_3 + \rho_3x_3^2$;
    \item $b_3 = \displaystyle \sum_{i=0}^3 \sigma_ix_2^{3-i}x_3^i = \sigma_0x_2^3 + \sigma_1x_2^2x_3 + \sigma_2x_2x_3^2 + \sigma_3x_3^3$;
     \item $c_1 = \delta_2x_2+\delta_3x_3$;
    \item $c_2 = \varepsilon_2x_2^2 + \varepsilon_{23}x_2x_3 + \varepsilon_3x_3^2$;
    \item $c_3 = \displaystyle \sum_{i=0}^3 \tau_ix_2^{3-i}x_3^i = \tau_0x_2^3 + \tau_1x_2^2x_3 + \tau_2x_2x_3^2 + \tau_3x_3^3$;
    \item $c_4 = \displaystyle \sum_{i=0}^4 \lambda_ix_2^{4-i}x_3^i = \lambda_0x_2^4 + \lambda_1x_2^3x_3 + \lambda_2x_2^2x_3^2 + \lambda_3x_2x_3^3 + \lambda_4x_3^4$.
\end{itemize}

\subsection{The $A_n$ case}

Let $D \subset \p^3$ be an irreducible normal quartic surface having a canonical singularity of type $A_n$ at $P=(1:0:0:0)$. Abusing notation, extend it by writing $A_i, i \leq 0$ to mean that $D$ is nonsingular at $P$, that is, $D$ is terminal at $P$. 

Let $\pi \colon X \rightarrow \p^3$ be the toric $(1,a,b)$-weighted blowup at $P$ and $E \coloneqq \Exc(\pi)$.

According to the toric description of the $(1,a,b)$-weighted blowup, the valuation associated to $E$ can be realized by an exceptional divisor on a sequence of $a$ blowups at points followed by $b-a$ blowups along curves. We will follow the notation adopted so far.

By Lemma \ref{res A_n} the effect of $a$ blowups at singular points is to weaken the singularity from $A_n$ to $A_{n-2a}$, whereas by Lemma \ref{res A_n lines} the remaining $b-a$ blowups along curves will weaken it from $A_{n-2a}$ to $A_{n-2a-(b-a)}=A_{n-a-b}$. 

\begin{lem}\label{bound a+b leq n+1}
Let $X,D,P \in D$ and $\pi \colon X \rightarrow \p^3$ be as above and suppose that $D$ has an $A_n$ singular point at $P$. If $\pi$ is volume-preserving, then $a \leq \left\lceil \dfrac{n}{2} \right\rceil$ and
\begin{center}
    $b-a \leq n-2a+1 \Rightarrow a+b \leq n+1$.
\end{center}
\end{lem}

\begin{proof}
    The first inequality is simply Lemma \ref{bound a for A_n}.
    If $b-a > n-2a+1$, then the types of singularity of the strict transforms of $D$ in the first $n-a+1$ steps of the toric description of $\pi$ are as follows:
\begin{equation}\label{tor desc A_n}
        \underbrace{A_n \overset{\Bl_{P_0}} \longrightarrow A_{n-2} \overset{\Bl_{P_1}} \longrightarrow \ldots \overset{\Bl_{P_{a-1}}} \longrightarrow A_{n-2a}}_{a\ \text{blowups at singular points}} \underbrace{\overset{\Bl_{L_a}} \longrightarrow A_{n-2a-1} \overset{\Bl_{L_{a+1}}} \longrightarrow \ldots 
         \overset{\Bl_{L_{n-a-2}}} \longrightarrow A_1 \overset{\Bl_{L_{n-a-1}}} \longrightarrow A_0 \overset{\Bl_{L_{n-a}}} \longrightarrow A_{-1} }_{n-2a+1\ \text{blowups along curves}}.
\end{equation}

In terms of the ambient varieties we have:

\begin{equation}
        \underbrace{\p^3=X_0 \overset{\pi_1} \longleftarrow X_1 \overset{\pi_2} \longleftarrow \ldots \overset{\pi_a} \longleftarrow X_a}_{a\ \text{blowups at singular points}} \underbrace{\overset{\pi_{a+1}} \longleftarrow X_{a+1} \overset{\pi_{a+2}} \longleftarrow \ldots \overset{\pi_{n-a-1}} \longleftarrow X_{n-a-1} \overset{\pi_{n-a}} \longleftarrow X_{n-a} \overset{\pi_{n-a+1}} \longleftarrow X_{n-a+1} }_{n-2a+1\ \text{blowups along curves}}.
\end{equation}

We observe that $n-2a \geq -1$. Otherwise,
\begin{center}
    $n-2a < -1 \Rightarrow a > \dfrac{n+1}{2} \geq \left\lceil \dfrac{n}{2} \right\rceil$,
\end{center}
contradicting Corollary \ref{bound a for A_n}. 

By the Key Property \ref{key prop 3 tor desc} of the toric description, the next center $z_{n-a+1} \in X_{n-a+1}$ in \ref{tor desc A_n} is the generic point of a section $L_{n-a+1} \subset E_{n-a+1}$ of the projection $E_{n-a+1} \rightarrow L_{n-a}$ disjoint from $E_{n-a}^{n-a+1}$. 

On the other hand, since $D_{n-a+1} \subset X_{n-a+1}$ is a nonsingular surface, $z_{n-a+1}$ is the generic point of a curve in $D_{n-a+1} \cap E_{n-a+1}$ by \cite[Proposition 3.1]{acm}. Since 
\begin{center}
    $L_{n-a+1} \coloneqq D_{n-a+1} \cap E_{n-a+1} \simeq \p^1$
\end{center}
\nin is irreducible, we must have $z_{n-a+1} = L_{n-a+1}$. 

Let us show that the curves $L_{n-a+1}$ and $E_{n-a}^{n-a+1} \cap E_{n-a+1}$ intersect. This contradicts the Key Property \ref{key prop 3 tor desc} of the toric description, and therefore we must have $b-a \leq n-2a+1$.

This nonempty intersection follows from the subsequent discussion.

By description \ref{tor desc A_n}, $D_j$ is singular at a point $P_j$ if $j \leq n-a-1$, and is nonsingular otherwise.

Suppose $a+1 \leq j \leq n-a  < b$. Then $E_j \simeq \p(\mathcal{N}_{L_{j-1}/X_{j-1}}^{\vee})\simeq \F_m$ for some $m \geq 0$, parametrizing all the normal directions to $L_{j-1}$. The projection $p_j \colon E_j \rightarrow L_{j-1}$ corresponds to the $\p^1$-bundle structure of $\F_m \rightarrow \p^1$ and each fiber represents all the normal directions to $L_{j-1}$ at its respective image.

Notice that $\pi_j|_{D_j} \colon D_j \rightarrow D_{j-1}$ is the blowup of $D_{j-1}$ along $L_{j-1}$ which a is non Cartier divisor of $D_{j-1}$ passing through its singularity. One has $\Exc(\pi_j|_{D_j}) = E_j \cap D_j$ is a Cartier divisor of $D_j$ and $D_j \not\simeq D_{j-1}$.

By the geometric properties of the blowup in this setting, $E_j \cap D_j$ is reducible, consisting of the union of a fiber $L_j'$ of $p_j$ corresponding to all the normal directions to $L_{j-1}$ at $P_{j-1}$ which form a $\p^1$, and a horizontal section $L_j$ of $p_j$. The setting is depicted in Figure \ref{fig:case 1 - a+1 leq j leq n-a}.

Furthermore, $L_j' \cap E_{j-1}^j=\{Q_j\}$ and $Q_j$ corresponds to the normal direction $L_{j-1}'$ to $L_{j-1}$ determined by $T_{P_{j-1}}E_{j-1}$. 
One has $L_j \cap L_j' = \{P_j\}$ and $P_j$ corresponds to the ``limit'' normal direction to $L_{j-1}$ at $P_{j-1}$ coming from the normal directions that define $T_QD_{j-1}$ for $Q \neq P_{j-1}$. 


\begin{center}


 
\tikzset{
pattern size/.store in=\mcSize, 
pattern size = 5pt,
pattern thickness/.store in=\mcThickness, 
pattern thickness = 0.3pt,
pattern radius/.store in=\mcRadius, 
pattern radius = 1pt}
\makeatletter
\pgfutil@ifundefined{pgf@pattern@name@_hyq4ior04}{
\pgfdeclarepatternformonly[\mcThickness,\mcSize]{_hyq4ior04}
{\pgfqpoint{-\mcThickness}{-\mcThickness}}
{\pgfpoint{\mcSize}{\mcSize}}
{\pgfpoint{\mcSize}{\mcSize}}
{
\pgfsetcolor{\tikz@pattern@color}
\pgfsetlinewidth{\mcThickness}
\pgfpathmoveto{\pgfpointorigin}
\pgfpathlineto{\pgfpoint{0}{\mcSize}}
\pgfusepath{stroke}
}}
\makeatother
\tikzset{every picture/.style={line width=0.75pt}} 

\begin{tikzpicture}[x=0.5pt,y=0.5pt,yscale=-1,xscale=1]

\draw  [color={rgb, 255:red, 189; green, 16; blue, 224 }  ,draw opacity=1 ][pattern=_hyq4ior04,pattern size=7.5pt,pattern thickness=0.75pt,pattern radius=0pt, pattern color={rgb, 255:red, 189; green, 16; blue, 224}] (189.71,43.51) -- (477.39,113.48) -- (478.66,309.5) -- (190.97,239.52) -- cycle ;
\draw  [color={rgb, 255:red, 209; green, 2; blue, 27 }  ,draw opacity=1 ][line width=1.5]  (259.37,50.74) -- (544.86,122.45) -- (398.94,224.24) -- (113.45,152.53) -- cycle ;
\draw [color={rgb, 255:red, 248; green, 231; blue, 28 }  ,draw opacity=1 ][line width=1.5]    (409,97.6) -- (410,293.6) ;
\draw [color={rgb, 255:red, 248; green, 231; blue, 28 }  ,draw opacity=1 ][line width=1.5]    (188,100.6) -- (476,171) ;
\draw [line width=1.5]    (310,330) -- (507,189.6) ;
\draw  [fill={rgb, 255:red, 126; green, 211; blue, 33 }  ,fill opacity=1 ] (404.6,154.98) .. controls (404.6,152.77) and (406.39,150.98) .. (408.6,150.98) .. controls (410.81,150.98) and (412.6,152.77) .. (412.6,154.98) .. controls (412.6,157.19) and (410.81,158.98) .. (408.6,158.98) .. controls (406.39,158.98) and (404.6,157.19) .. (404.6,154.98) -- cycle ;
\draw    (321.44,377.12) -- (322.58,502.2) ;
\draw [shift={(322.6,504.2)}, rotate = 269.48] [color={rgb, 255:red, 0; green, 0; blue, 0 }  ][line width=0.75]    (10.93,-3.29) .. controls (6.95,-1.4) and (3.31,-0.3) .. (0,0) .. controls (3.31,0.3) and (6.95,1.4) .. (10.93,3.29)   ;
\draw   (268,741.3) .. controls (268,731.75) and (291.51,724) .. (320.5,724) .. controls (349.49,724) and (373,731.75) .. (373,741.3) .. controls (373,750.85) and (349.49,758.6) .. (320.5,758.6) .. controls (291.51,758.6) and (268,750.85) .. (268,741.3) -- cycle ;
\draw   (268,547.3) .. controls (268,537.75) and (291.51,530) .. (320.5,530) .. controls (349.49,530) and (373,537.75) .. (373,547.3) .. controls (373,556.85) and (349.49,564.6) .. (320.5,564.6) .. controls (291.51,564.6) and (268,556.85) .. (268,547.3) -- cycle ;
\draw [color={rgb, 255:red, 189; green, 16; blue, 224 }  ,draw opacity=1 ][line width=1.5]    (268,547.3) -- (373,741.3) ;
\draw [line width=1.5]    (373,547.3) -- (268,741.3) ;
\draw  [color={rgb, 255:red, 209; green, 2; blue, 27 }  ,draw opacity=1 ][line width=1.5]  (188.79,541.8) -- (427.55,541.8) -- (475.84,767.8) -- (237.08,767.8) -- cycle ;
\draw [color={rgb, 255:red, 248; green, 231; blue, 28 }  ,draw opacity=1 ][line width=1.5]    (281,571) -- (304.02,573.48) ;
\draw [shift={(307,573.8)}, rotate = 186.15] [color={rgb, 255:red, 248; green, 231; blue, 28 }  ,draw opacity=1 ][line width=1.5]    (14.21,-4.28) .. controls (9.04,-1.82) and (4.3,-0.39) .. (0,0) .. controls (4.3,0.39) and (9.04,1.82) .. (14.21,4.28)   ;
\draw [color={rgb, 255:red, 248; green, 231; blue, 28 }  ,draw opacity=1 ][line width=1.5]    (288,585) -- (311.02,587.48) ;
\draw [shift={(314,587.8)}, rotate = 186.15] [color={rgb, 255:red, 248; green, 231; blue, 28 }  ,draw opacity=1 ][line width=1.5]    (14.21,-4.28) .. controls (9.04,-1.82) and (4.3,-0.39) .. (0,0) .. controls (4.3,0.39) and (9.04,1.82) .. (14.21,4.28)   ;
\draw [color={rgb, 255:red, 248; green, 231; blue, 28 }  ,draw opacity=1 ][line width=1.5]    (294,599) -- (317.02,601.48) ;
\draw [shift={(320,601.8)}, rotate = 186.15] [color={rgb, 255:red, 248; green, 231; blue, 28 }  ,draw opacity=1 ][line width=1.5]    (14.21,-4.28) .. controls (9.04,-1.82) and (4.3,-0.39) .. (0,0) .. controls (4.3,0.39) and (9.04,1.82) .. (14.21,4.28)   ;
\draw [color={rgb, 255:red, 248; green, 231; blue, 28 }  ,draw opacity=1 ][line width=1.5]    (303.33,611) -- (326.35,613.48) ;
\draw [shift={(329.33,613.8)}, rotate = 186.15] [color={rgb, 255:red, 248; green, 231; blue, 28 }  ,draw opacity=1 ][line width=1.5]    (14.21,-4.28) .. controls (9.04,-1.82) and (4.3,-0.39) .. (0,0) .. controls (4.3,0.39) and (9.04,1.82) .. (14.21,4.28)   ;
\draw [color={rgb, 255:red, 248; green, 231; blue, 28 }  ,draw opacity=1 ][line width=1.5]    (309.33,625) -- (332.35,627.48) ;
\draw [shift={(335.33,627.8)}, rotate = 186.15] [color={rgb, 255:red, 248; green, 231; blue, 28 }  ,draw opacity=1 ][line width=1.5]    (14.21,-4.28) .. controls (9.04,-1.82) and (4.3,-0.39) .. (0,0) .. controls (4.3,0.39) and (9.04,1.82) .. (14.21,4.28)   ;
\draw [color={rgb, 255:red, 248; green, 231; blue, 28 }  ,draw opacity=1 ][line width=1.5]    (321.33,645) -- (344.35,647.48) ;
\draw [shift={(347.33,647.8)}, rotate = 186.15] [color={rgb, 255:red, 248; green, 231; blue, 28 }  ,draw opacity=1 ][line width=1.5]    (14.21,-4.28) .. controls (9.04,-1.82) and (4.3,-0.39) .. (0,0) .. controls (4.3,0.39) and (9.04,1.82) .. (14.21,4.28)   ;
\draw [color={rgb, 255:red, 248; green, 231; blue, 28 }  ,draw opacity=1 ][line width=1.5]    (327.33,659) -- (350.35,661.48) ;
\draw [shift={(353.33,661.8)}, rotate = 186.15] [color={rgb, 255:red, 248; green, 231; blue, 28 }  ,draw opacity=1 ][line width=1.5]    (14.21,-4.28) .. controls (9.04,-1.82) and (4.3,-0.39) .. (0,0) .. controls (4.3,0.39) and (9.04,1.82) .. (14.21,4.28)   ;
\draw [color={rgb, 255:red, 248; green, 231; blue, 28 }  ,draw opacity=1 ][line width=1.5]    (336,671) -- (359.02,673.48) ;
\draw [shift={(362,673.8)}, rotate = 186.15] [color={rgb, 255:red, 248; green, 231; blue, 28 }  ,draw opacity=1 ][line width=1.5]    (14.21,-4.28) .. controls (9.04,-1.82) and (4.3,-0.39) .. (0,0) .. controls (4.3,0.39) and (9.04,1.82) .. (14.21,4.28)   ;
\draw [color={rgb, 255:red, 248; green, 231; blue, 28 }  ,draw opacity=1 ][line width=1.5]    (342,685) -- (365.02,687.48) ;
\draw [shift={(368,687.8)}, rotate = 186.15] [color={rgb, 255:red, 248; green, 231; blue, 28 }  ,draw opacity=1 ][line width=1.5]    (14.21,-4.28) .. controls (9.04,-1.82) and (4.3,-0.39) .. (0,0) .. controls (4.3,0.39) and (9.04,1.82) .. (14.21,4.28)   ;
\draw [color={rgb, 255:red, 248; green, 231; blue, 28 }  ,draw opacity=1 ][line width=1.5]    (350.67,697.67) -- (373.68,700.15) ;
\draw [shift={(376.67,700.47)}, rotate = 186.15] [color={rgb, 255:red, 248; green, 231; blue, 28 }  ,draw opacity=1 ][line width=1.5]    (14.21,-4.28) .. controls (9.04,-1.82) and (4.3,-0.39) .. (0,0) .. controls (4.3,0.39) and (9.04,1.82) .. (14.21,4.28)   ;
\draw [color={rgb, 255:red, 248; green, 231; blue, 28 }  ,draw opacity=1 ][line width=1.5]    (356.67,711.67) -- (379.68,714.15) ;
\draw [shift={(382.67,714.47)}, rotate = 186.15] [color={rgb, 255:red, 248; green, 231; blue, 28 }  ,draw opacity=1 ][line width=1.5]    (14.21,-4.28) .. controls (9.04,-1.82) and (4.3,-0.39) .. (0,0) .. controls (4.3,0.39) and (9.04,1.82) .. (14.21,4.28)   ;
\draw [color={rgb, 255:red, 248; green, 231; blue, 28 }  ,draw opacity=1 ][line width=1.5]    (364,724.33) -- (387.02,726.81) ;
\draw [shift={(390,727.13)}, rotate = 186.15] [color={rgb, 255:red, 248; green, 231; blue, 28 }  ,draw opacity=1 ][line width=1.5]    (14.21,-4.28) .. controls (9.04,-1.82) and (4.3,-0.39) .. (0,0) .. controls (4.3,0.39) and (9.04,1.82) .. (14.21,4.28)   ;
\draw [color={rgb, 255:red, 248; green, 231; blue, 28 }  ,draw opacity=1 ][line width=1.5]    (370,738.33) -- (393.02,740.81) ;
\draw [shift={(396,741.13)}, rotate = 186.15] [color={rgb, 255:red, 248; green, 231; blue, 28 }  ,draw opacity=1 ][line width=1.5]    (14.21,-4.28) .. controls (9.04,-1.82) and (4.3,-0.39) .. (0,0) .. controls (4.3,0.39) and (9.04,1.82) .. (14.21,4.28)   ;
\draw [color={rgb, 255:red, 248; green, 231; blue, 28 }  ,draw opacity=1 ][line width=1.5]    (321.33,645) -- (340.08,629.64) ;
\draw [shift={(342.4,627.73)}, rotate = 140.66] [color={rgb, 255:red, 248; green, 231; blue, 28 }  ,draw opacity=1 ][line width=1.5]    (14.21,-4.28) .. controls (9.04,-1.82) and (4.3,-0.39) .. (0,0) .. controls (4.3,0.39) and (9.04,1.82) .. (14.21,4.28)   ;
\draw [color={rgb, 255:red, 248; green, 231; blue, 28 }  ,draw opacity=1 ][line width=1.5]    (321.33,645) -- (301.8,653.83) ;
\draw [shift={(299.07,655.07)}, rotate = 335.67] [color={rgb, 255:red, 248; green, 231; blue, 28 }  ,draw opacity=1 ][line width=1.5]    (14.21,-4.28) .. controls (9.04,-1.82) and (4.3,-0.39) .. (0,0) .. controls (4.3,0.39) and (9.04,1.82) .. (14.21,4.28)   ;
\draw [color={rgb, 255:red, 248; green, 231; blue, 28 }  ,draw opacity=1 ][line width=1.5]    (320.5,644.3) -- (296.5,639.01) ;
\draw [shift={(293.57,638.37)}, rotate = 12.42] [color={rgb, 255:red, 248; green, 231; blue, 28 }  ,draw opacity=1 ][line width=1.5]    (14.21,-4.28) .. controls (9.04,-1.82) and (4.3,-0.39) .. (0,0) .. controls (4.3,0.39) and (9.04,1.82) .. (14.21,4.28)   ;
\draw [color={rgb, 255:red, 248; green, 231; blue, 28 }  ,draw opacity=1 ][line width=1.5]    (321.33,645) -- (311.54,668.96) ;
\draw [shift={(310.4,671.73)}, rotate = 292.24] [color={rgb, 255:red, 248; green, 231; blue, 28 }  ,draw opacity=1 ][line width=1.5]    (14.21,-4.28) .. controls (9.04,-1.82) and (4.3,-0.39) .. (0,0) .. controls (4.3,0.39) and (9.04,1.82) .. (14.21,4.28)   ;
\draw [color={rgb, 255:red, 248; green, 231; blue, 28 }  ,draw opacity=1 ][line width=1.5]    (321.33,645) -- (324.69,671.42) ;
\draw [shift={(325.07,674.4)}, rotate = 262.76] [color={rgb, 255:red, 248; green, 231; blue, 28 }  ,draw opacity=1 ][line width=1.5]    (14.21,-4.28) .. controls (9.04,-1.82) and (4.3,-0.39) .. (0,0) .. controls (4.3,0.39) and (9.04,1.82) .. (14.21,4.28)   ;
\draw [color={rgb, 255:red, 248; green, 231; blue, 28 }  ,draw opacity=1 ][line width=1.5]    (321.33,645) -- (323.55,621.95) ;
\draw [shift={(323.83,618.97)}, rotate = 95.49] [color={rgb, 255:red, 248; green, 231; blue, 28 }  ,draw opacity=1 ][line width=1.5]    (14.21,-4.28) .. controls (9.04,-1.82) and (4.3,-0.39) .. (0,0) .. controls (4.3,0.39) and (9.04,1.82) .. (14.21,4.28)   ;
\draw [color={rgb, 255:red, 248; green, 231; blue, 28 }  ,draw opacity=1 ][line width=1.5]    (320.5,644.3) -- (305.12,627.92) ;
\draw [shift={(303.07,625.73)}, rotate = 46.8] [color={rgb, 255:red, 248; green, 231; blue, 28 }  ,draw opacity=1 ][line width=1.5]    (14.21,-4.28) .. controls (9.04,-1.82) and (4.3,-0.39) .. (0,0) .. controls (4.3,0.39) and (9.04,1.82) .. (14.21,4.28)   ;
\draw [color={rgb, 255:red, 248; green, 231; blue, 28 }  ,draw opacity=1 ][line width=1.5]    (322.78,645.77) -- (339.67,661.77) ;
\draw [shift={(341.85,663.83)}, rotate = 223.46] [color={rgb, 255:red, 248; green, 231; blue, 28 }  ,draw opacity=1 ][line width=1.5]    (14.21,-4.28) .. controls (9.04,-1.82) and (4.3,-0.39) .. (0,0) .. controls (4.3,0.39) and (9.04,1.82) .. (14.21,4.28)   ;
\draw  [fill={rgb, 255:red, 126; green, 211; blue, 33 }  ,fill opacity=1 ] (316.5,644.3) .. controls (316.5,642.09) and (318.29,640.3) .. (320.5,640.3) .. controls (322.71,640.3) and (324.5,642.09) .. (324.5,644.3) .. controls (324.5,646.51) and (322.71,648.3) .. (320.5,648.3) .. controls (318.29,648.3) and (316.5,646.51) .. (316.5,644.3) -- cycle ;
\draw [color={rgb, 255:red, 74; green, 144; blue, 226 }  ,draw opacity=1 ][line width=1.5]    (188.64,201.76) -- (477.44,276.96) ;
\draw  [fill={rgb, 255:red, 80; green, 227; blue, 194 }  ,fill opacity=1 ] (404.5,259.8) .. controls (404.5,257.59) and (406.29,255.8) .. (408.5,255.8) .. controls (410.71,255.8) and (412.5,257.59) .. (412.5,259.8) .. controls (412.5,262.01) and (410.71,263.8) .. (408.5,263.8) .. controls (406.29,263.8) and (404.5,262.01) .. (404.5,259.8) -- cycle ;

\draw (328,340.4) node [anchor=north west][inner sep=0.75pt]    {$L'_{j-1}$};
\draw (392,303.4) node [anchor=north west][inner sep=0.75pt]  [color={rgb, 255:red, 248; green, 231; blue, 28 }  ,opacity=1 ]  {$L'_{j}$};
\draw (496,160.4) node [anchor=north west][inner sep=0.75pt]  [color={rgb, 255:red, 248; green, 231; blue, 28 }  ,opacity=1 ]  {$E_{j} \ \cap D_{j}$};
\draw (159,79.4) node [anchor=north west][inner sep=0.75pt]  [color={rgb, 255:red, 248; green, 231; blue, 28 }  ,opacity=1 ]  {$L_{j}$};
\draw (181,12.4) node [anchor=north west][inner sep=0.75pt]  [color={rgb, 255:red, 189; green, 16; blue, 224 }  ,opacity=1 ]  {$E_{j} \ \simeq \mathbb{F}_{m}$};
\draw (75,163.4) node [anchor=north west][inner sep=0.75pt]  [color={rgb, 255:red, 209; green, 2; blue, 27 }  ,opacity=1 ]  {$\widetilde{\{x_{3} =0\}}$};
\draw (209.8,418) node [anchor=north west][inner sep=0.75pt]    {$\pi _{j} =\text{Bl}_{L_{j-1}}$};
\draw (491.6,264.24) node [anchor=north west][inner sep=0.75pt]  [color={rgb, 255:red, 74; green, 144; blue, 226 }  ,opacity=1 ]  {$E_{j-1}^{j} \cap E_{j}$};
\draw (420.4,273) node [anchor=north west][inner sep=0.75pt]  [color={rgb, 255:red, 80; green, 227; blue, 194 }  ,opacity=1 ]  {$Q_{j}$};
\draw (442.2,529.8) node [anchor=north west][inner sep=0.75pt]  [color={rgb, 255:red, 209; green, 2; blue, 7 }  ,opacity=1 ]  {$\{x_{3} =0\}$};
\draw (377.8,119.2) node [anchor=north west][inner sep=0.75pt]  [color={rgb, 255:red, 126; green, 211; blue, 33 }  ,opacity=1 ]  {$P_{j}$};
\draw (260.4,622.12) node [anchor=north west][inner sep=0.75pt]  [color={rgb, 255:red, 126; green, 211; blue, 33 }  ,opacity=1 ]  {$P_{j-1}$};
\draw (371.2,559.72) node [anchor=north west][inner sep=0.75pt]    {$L'_{j-1}$};
\draw (229.2,555.12) node [anchor=north west][inner sep=0.75pt]  [color={rgb, 255:red, 189; green, 16; blue, 224 }  ,opacity=1 ]  {$L_{j-1}$};
\draw (258,504.72) node [anchor=north west][inner sep=0.75pt]    {$D_{j-1}$};

\end{tikzpicture}
\end{center}
\begin{figure}[htb]
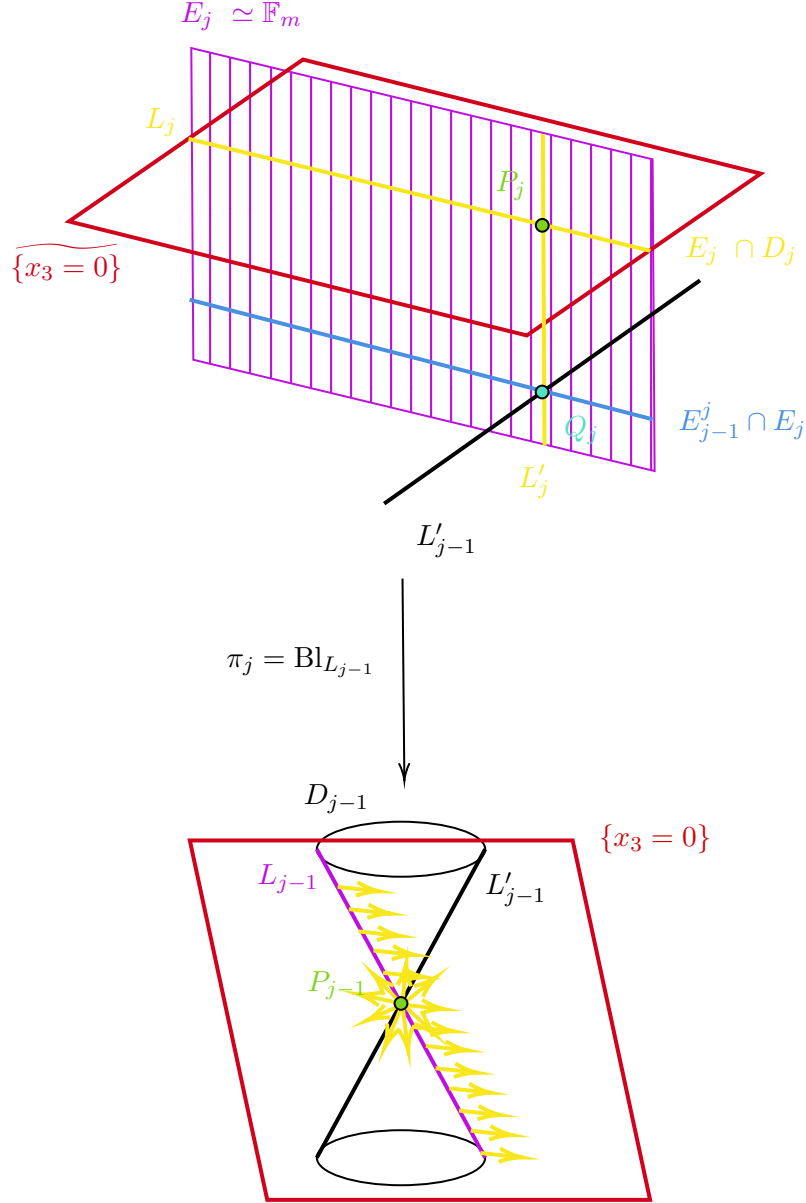

\caption{Case $a+1 \leq j \leq n-a  < b$.}
\label{fig:case 1 - a+1 leq j leq n-a}
\end{figure}

Note that $D_{n-a}$ is nonsingular and $\pi_{n-a+1} \colon X_{n-a+1} \rightarrow X_{n-a}$, the blowup of $L_{n-a} \subset D_{n-a}$, is such that $D_{n-a+1} = \widetilde{D_{n-a}} \simeq D_{n-a}$. So, for $j=n-a+1$, we have that $E_{n-a+1} \cap D_{n-a+1}$ is a section of the projection $p_{n-a+1} \colon E_{n-a+1} \rightarrow L_{n-a}$.

Indeed, since $D_{n-a}$ is nonsingular, for each $Q' \in L_{n-a}$, there exists only one normal direction to $L_{n-a}$ at $Q'$ that is tangent to $D_{n-a}$.

Furthermore, $E_{n-a+1} \cap D_{n-a+1} \cap E_{n-a}^{n-a+1} = \{Q_{n-a+1}\} \neq \emptyset$, where $Q_{n-a+1}$ corresponds to the normal direction to $L_{n-a}$ at $P_{n-a}$ determined by $T_{P_{n-a}}E_{n-a} = T_{P_{n-a}}D_{n-a}$, and this concludes the proof. The setting is depicted in Figure \ref{fig:case j=n-a+1}.

\begin{center}

 
\tikzset{
pattern size/.store in=\mcSize, 
pattern size = 5pt,
pattern thickness/.store in=\mcThickness, 
pattern thickness = 0.3pt,
pattern radius/.store in=\mcRadius, 
pattern radius = 1pt}
\makeatletter
\pgfutil@ifundefined{pgf@pattern@name@_rk56uwo6q}{
\pgfdeclarepatternformonly[\mcThickness,\mcSize]{_rk56uwo6q}
{\pgfqpoint{-\mcThickness}{-\mcThickness}}
{\pgfpoint{\mcSize}{\mcSize}}
{\pgfpoint{\mcSize}{\mcSize}}
{
\pgfsetcolor{\tikz@pattern@color}
\pgfsetlinewidth{\mcThickness}
\pgfpathmoveto{\pgfpointorigin}
\pgfpathlineto{\pgfpoint{0}{\mcSize}}
\pgfusepath{stroke}
}}
\makeatother
\tikzset{every picture/.style={line width=0.75pt}} 

\begin{tikzpicture}[x=0.5pt,y=0.5pt,yscale=-1,xscale=1]

\draw  [color={rgb, 255:red, 189; green, 16; blue, 224 }  ,draw opacity=1 ][pattern=_rk56uwo6q,pattern size=7.5pt,pattern thickness=0.75pt,pattern radius=0pt, pattern color={rgb, 255:red, 189; green, 16; blue, 224}] (520.88,87.45) -- (808.57,157.42) -- (809.83,353.44) -- (522.15,283.47) -- cycle ;
\draw [line width=1.5]    (604.6,347.2) -- (832,184.6) ;
\draw [color={rgb, 255:red, 74; green, 144; blue, 226 }  ,draw opacity=1 ][line width=1.5]    (520.78,188.9) -- (809.58,264.1) ;
\draw [color={rgb, 255:red, 248; green, 231; blue, 28 }  ,draw opacity=1 ][line width=1.5]    (522.33,143) .. controls (615.37,70.6) and (769.33,310.8) .. (809.33,280.8) ;
\draw  [fill={rgb, 255:red, 80; green, 227; blue, 194 }  ,fill opacity=1 ] (740.17,248.13) .. controls (740.17,245.92) and (741.96,244.13) .. (744.17,244.13) .. controls (746.38,244.13) and (748.17,245.92) .. (748.17,248.13) .. controls (748.17,250.34) and (746.38,252.13) .. (744.17,252.13) .. controls (741.96,252.13) and (740.17,250.34) .. (740.17,248.13) -- cycle ;
\draw    (476.2,193) -- (352.6,193.2) ;
\draw [shift={(350.6,193.2)}, rotate = 359.91] [color={rgb, 255:red, 0; green, 0; blue, 0 }  ][line width=0.75]    (10.93,-3.29) .. controls (6.95,-1.4) and (3.31,-0.3) .. (0,0) .. controls (3.31,0.3) and (6.95,1.4) .. (10.93,3.29)   ;
\draw  [color={rgb, 255:red, 74; green, 144; blue, 226 }  ,draw opacity=1 ][line width=1.5]  (62.5,104.2) -- (263.2,104.2) -- (263.2,299.2) -- (62.5,299.2) -- cycle ;
\draw [color={rgb, 255:red, 189; green, 16; blue, 224 }  ,draw opacity=1 ][line width=1.5]    (62,196.2) -- (262,195.2) ;
\draw [color={rgb, 255:red, 0; green, 0; blue, 0 }  ,draw opacity=1 ][line width=1.5]    (117,301.2) -- (119,102.2) ;
\draw  [fill={rgb, 255:red, 126; green, 211; blue, 33 }  ,fill opacity=1 ] (114,197.7) .. controls (114,195.49) and (115.79,193.7) .. (118,193.7) .. controls (120.21,193.7) and (122,195.49) .. (122,197.7) .. controls (122,199.91) and (120.21,201.7) .. (118,201.7) .. controls (115.79,201.7) and (114,199.91) .. (114,197.7) -- cycle ;
\draw [color={rgb, 255:red, 189; green, 16; blue, 224 }  ,draw opacity=1 ][line width=1.5]    (59,380.2) -- (259,379.2) ;
\draw    (150,316.2) -- (150.57,361.2) ;
\draw [shift={(150.6,363.2)}, rotate = 269.27] [color={rgb, 255:red, 0; green, 0; blue, 0 }  ][line width=0.75]    (10.93,-3.29) .. controls (6.95,-1.4) and (3.31,-0.3) .. (0,0) .. controls (3.31,0.3) and (6.95,1.4) .. (10.93,3.29)   ;

\draw (512,61.4) node [anchor=north west][inner sep=0.75pt]  [color={rgb, 255:red, 189; green, 16; blue, 224 }  ,opacity=1 ]  {$E_{n-a+1} \simeq \mathbb{F}_{m}$};
\draw (818,252.84) node [anchor=north west][inner sep=0.75pt]  [color={rgb, 255:red, 74; green, 144; blue, 226 }  ,opacity=1 ]  {$E_{n-a}^{n-a+1} \cap E_{n-a+1}$};
\draw (722.67,276.73) node [anchor=north west][inner sep=0.75pt]  [color={rgb, 255:red, 80; green, 227; blue, 194 }  ,opacity=1 ]  {$Q_{n-a+1}$};
\draw (619.33,352.07) node [anchor=north west][inner sep=0.75pt]    {$L'_{n-a}$};
\draw (819.4,285.2) node [anchor=north west][inner sep=0.75pt]  [color={rgb, 255:red, 248; green, 231; blue, 28 }  ,opacity=1 ]  {$E_{n-a+1} \ \cap D_{n-a+1}$};
\draw (345.8,163) node [anchor=north west][inner sep=0.75pt]    {$\pi _{n-a+1} =\text{Bl}_{L_{n-a}}$};
\draw (251,75.4) node [anchor=north west][inner sep=0.75pt]  [color={rgb, 255:red, 74; green, 144; blue, 226 }  ,opacity=1 ]  {$E_{n-a}$};
\draw (128.4,208.12) node [anchor=north west][inner sep=0.75pt]  [color={rgb, 255:red, 126; green, 211; blue, 33 }  ,opacity=1 ]  {$P_{n-a}$};
\draw (271.2,192.12) node [anchor=north west][inner sep=0.75pt]  [color={rgb, 255:red, 189; green, 16; blue, 224 }  ,opacity=1 ]  {$L_{n-a}$};
\draw (106.2,76.4) node [anchor=north west][inner sep=0.75pt]    {$L'_{n-a}$};
\draw (275.2,371.12) node [anchor=north west][inner sep=0.75pt]  [color={rgb, 255:red, 189; green, 16; blue, 224 }  ,opacity=1 ]  {$L_{n-a-1}$};
\draw (89,336.4) node [anchor=north west][inner sep=0.75pt]    {$p_{n-a}$};

\end{tikzpicture}
\end{center}

\begin{figure}[htb]
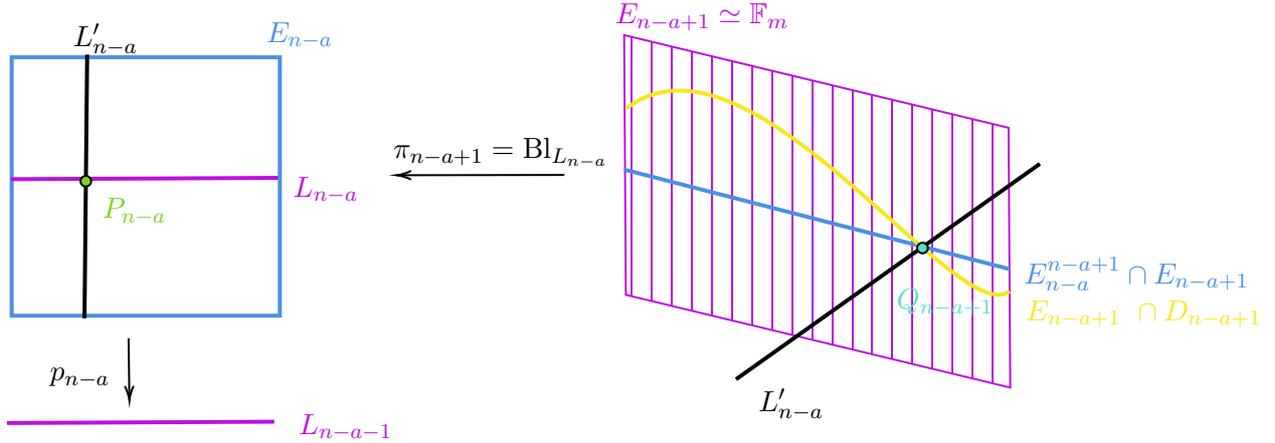

\caption{Case $j=n-a+1$.}
\label{fig:case j=n-a+1}
\end{figure}

\end{proof}

\paragraph{Notation:} Let $\varphi \colon Y \rightarrow \mathbb{P}^3$ be the blowup of $\p^3$ at the coordinate point $P=(1:0:0:0)$. One can describe as $$Y = V(\{x_iy_j-x_jy_i\}_{1\leq i<j \leq 3}) \subset \p^3_{(x_0:x_1:x_2:x_3)} \times\p^2_{(y_1:y_2:y_3)}$$ 
and one has the identification $U \coloneqq \varphi^{-1}(\{x_0 \neq 0\}) \simeq \Bl_{(0,0,0)}(\Aa^3)$. For almost all the purposes of this work involving such a blowup, the relevant computations will be done in the ``first affine chart'' $\{ y_1 \neq 0\}$, namely, 

$$U \cap \{y_1 \neq 0\} \simeq \Aa^3_{(x_1,y_2,y_3)},$$
\nin where $x_j = x_1y_j $ for $j \neq 1$. By abuse of notation, let $U \cap \{y_1 \neq 0\}$ always have affine coordinates $(x_1,x_2,x_3)$. We observe that the exceptional divisor is given by $\{x_1=0\}$ in this affine chart. 

The same notation will be employed analogously for the case of a blowup along a coordinate line.

\subsubsection{Criteria for $A_n$ singularities on a quartic surface}\label{criteria A_n}

Let $D = \{x_0^2A+x_0B+C = 0\} \subset \p^3$ be a quartic surface having a canonical singularity of type $A_n$ at $P = (1: 0 : 0: 0)$, $n \in \{1,\ldots,19\}$. Our task now is to establish conditions on $A,B$ and $C$ so that $P$ is of type $A_n$.

This was done in a different setting in \cite{kn} and included the $D_n$ case, where Kato \& Naruki also give a description of the parameter space of some of these quartics.

However, we are working in different coordinates than \cite{kn}. We will obtain our conditions by analyzing how the singularity is resolved by a sequence of blowups at singular points, and checking the geometry of the exceptional divisor along the process.  
The bigger $n$ is, the more complicated these explicit criteria will be in terms of the coefficients of the quartic.

In view of Theorem \ref{gue thm} and Lemma \ref{bound a+b leq n+1}, we will do it for $n \leq 7$.

In fact, $n \leq 6$ would be enough for the problem of determining which weighted blowups initiate a volume-preserving Sarkisov link, but we will also deal with the case $n=7$, since they are still manageable. 

From what was explained in Subsection \ref{notation quartic} concerning the notation for the coefficients of the homogeneous equation that defines $D$, we have that 
\begin{center}
    $P$ is $A_1 \Leftrightarrow \rk(A) = 3$.
\end{center}

Let us move to the $A_{\geq 2}$ case. Since a blowup at the singular point weakens the singularity from $A_n$ to $A_{n-2}$, to obtain the criteria for $n=7$, by Lemma \ref{res A_n} we will need to blow up 4 times.

We will denote by $D_i$ the strict transform of $D$ in the $i$-th blowup and $E_i$ the exceptional divisor. All the following equations are with respect to the ``first affine chart'' $\{ y_1 \neq 0\}$ in the corresponding blowup. We observe that $E_i = \{x_1 = 0\}$ in this affine chart.

At Step 0 we have our initial situation, and at Step $i$ the situation after the $i$-th blowup. After performing many computations and applying the Jacobian Criterion, one obtains the following:

\paragraph{Step 0:} $D_0 = \{x_0^2A+x_0B+C = 0\}$. In the affine open $\{x_0 \neq 0\}$ we have
\begin{center}
    $D_0 = \{A + B + C = 0 \} \subset \Aa_{(x_1,x_2,x_3)}^3$.
\end{center}

Set $P_0 \coloneqq (0,0,0) \in \Sing( D_{0})$. One has
\begin{center}
    $P_0$ is $A_{\geq 2} \Leftrightarrow \rk(A) = 2 \Leftrightarrow A = x_2x_3$, without loss of generality.
\end{center}

\paragraph{Step 1:} $D_1 = \{x_2x_3+x_1B(1,x_2,x_3)+x_1^2C(1,x_2,x_3) = 0\} \subset \Aa_{(x_1,x_2,x_3)}^3$. We have 
\begin{center}
    $E_1 \cap D_1 = \{x_1 = x_2x_3 = 0\} = L_1 \cup L_1'$, 
\end{center}

\nin where $L_1$ and $L_1'$ are lines on $D_1 \subset \Aa^3$. Observe that $P_1 \coloneqq (0,0,0) \in \Sing(E_1 \cap D_1)$. One has
\begin{align*}
P_0 \ \text{is} \ A_2 & \Leftrightarrow P_1 \ \text{is} \ A_0 \Leftrightarrow b_0 \neq 0,\\
P_0 \ \text{is} \ A_{\geq 3} & \Leftrightarrow P_1 \ \text{is} \ A_{\geq 1} \Leftrightarrow b_0 = 0 \ (*_1),
\end{align*}

\nin where $b_0$ is the coefficient of the monomial $x_1^3$ in the defining equation of $D_0$ in the afine open $\{x_0 \neq 0\}$.

\paragraph{Step 2:} 

\begin{center}
    $D_2 = \Biggl\{x_2x_3+b_1(x_2,x_3)+x_1b_2(x_2,x_3)+x_1^2b_3(x_2,x_3) + \displaystyle \sum_{i=0}^4 x_1^ic_i(x_2,x_3) = 0\Bigg\} \subset \Aa_{(x_1,x_2,x_3)}^3$.
\end{center}

We have that
\begin{center}
    $E_2 \cap D_2 = \{x_1 = x_2x_3 + \beta_2x_2 + \beta_3x_3 +c_0 = 0\}$, 
\end{center}

\nin is a conic (not necessarily irreducible) on $D_1 \subset \Aa^3$. One has
\begin{align*}
    P_0 \ \text{is} \ A_3 & \Leftrightarrow P_1 \ \text{is} \ A_1 \Leftrightarrow (*_1) \ \text{and}\  c_0 \neq \beta_2\beta_3,\\
    P_0 \ \text{is} \ A_{\geq 4} & \Leftrightarrow P_1 \ \text{is} \ A_{\geq 2} \Leftrightarrow (*_1) \ \text{and}\  c_0 = \beta_2\beta_3 \ (*_2).
\end{align*}

In the latter case, we have 
\begin{center}
    $E_2 \cap D_2 = \{x_1 = (x_2+\beta_3)(x_3+\beta_2)= 0\} = L_2 \cup L_2'$,
\end{center}

\nin where $L_2$ and $L_2'$ are lines on $D_2$. Observe that $P_2 \coloneqq (0,-\beta_3,-\beta_2) \in \Sing(E_2 \cap D_2)$. One has
\begin{align*}
    P_0 \ \text{is} \ A_4 & \Leftrightarrow (*_2)\ \text{and}\ b_2(\beta_3,\beta_2)-c_1(\beta_3,\beta_2) = \rho_2\beta_3^2 + \rho_{23}\beta_2\beta_3 + \rho_3\beta_2^2 - \delta_2\beta_3 - \delta_3\beta_2 \eqqcolon \zeta \neq 0, \\
    P_0 \ \text{is} \ A_{\geq 5} & \Leftrightarrow P_2 \ \text{is} \ A_{\geq 1} \Leftrightarrow (*_2) \ \text{and}\  \zeta = 0 \ (*_3).
\end{align*}

\paragraph{Step 3:} \begin{multline*}
D_3 = \Biggl\{x_2x_3 + \rho_2x_1x_2 - 2\rho_2\beta_3x_2 + \rho_{23}(x_1x_2x_3-\beta_2x_2-\beta_3x_3) + \rho_3x_1x_3 -2\rho_3\beta_2x_3  \\+ b_3(x_1x_2-\beta_3,x_1x_3-\beta_2) + \delta_2x_2 + \delta_3x_3 + \displaystyle \sum_{i=2}^4 x_1^{i-2}c_i(x_1x_2-\beta_3,x_1x_3-\beta_2) = 0\Biggl\} \\
\subset \Aa_{(x_1,x_2,x_3)}^3.
\end{multline*}

Set 
\begin{align*}
    \xi_2 & \coloneqq -2\rho_2\beta_3 - \rho_{23}\beta_2 + \delta_2, \\
    \xi_3 & \coloneqq -2\rho_3\beta_2 - \rho_{23}\beta_3 + \delta_3 \ \text{and} \\
    \alpha & \coloneqq -b_3(\beta_3,\beta_2)+c_2(\beta_3,\beta_2).
\end{align*}

We have that

\begin{center}
    $E_3 \cap D_3 = \{x_1 = x_2x_3 + \xi_2x_2 + \xi_3x_3 +\alpha = 0\}$, 
\end{center}

\nin is a conic (not necessarily irreducible) on $D_3 \subset \Aa^3$. One has
\begin{align*}
    P_0 \ \text{is} \ A_5 & \Leftrightarrow P_1 \ \text{is} \ A_3 \Leftrightarrow (*_3) \ \text{and}\  \xi_2\xi_3 \neq \alpha,\\
    P_0 \ \text{is} \ A_{\geq 6} & \Leftrightarrow P_1 \ \text{is} \ A_{\geq 4} \Leftrightarrow (*_3) \ \text{and}\  \xi_2\xi_3 = \alpha \ (*_4).
\end{align*}

In the latter case, we have 
\begin{center}
    $E_3 \cap D_3 = \{x_1 = (x_2+\xi_3)(x_3+\xi_2)= 0\} = L_3 \cup L_3'$,
\end{center}

\nin where $L_3$ and $L_3'$ are lines on $D_3 \subset \Aa^3$. Observe that $P_3 \coloneqq (0,-\xi_3,-\xi_2) \in \Sing(E_3 \cap D_3)$.

Set 
\begin{align*}
    \omega & \coloneqq -3\sigma_0\xi_3\beta_3^2+\sigma_1(-2\xi_3\beta_2\beta_3-\xi_2\beta_3^2) +\sigma_2(-\xi_3\beta_2^2-2\xi_2\beta_2\beta_3) - 3\sigma_3\xi_2\beta_2^2 \ \text{and} \\
    \eta & \coloneqq 2\varepsilon_2\xi_3\beta_3+\varepsilon_{23}\xi_3\beta_2+\varepsilon_{23}\xi_2\beta_3 + 2\varepsilon_3\xi_2\beta_2.    
\end{align*}

One has
\begin{align*}
    P_0 \ \text{is} \ A_6 & \Leftrightarrow (*_4)\ \text{and}\ b_2(\xi_3,\xi_2) + \omega + \eta -c_3(\beta_3,\beta_2) \eqqcolon \theta \neq 0, \\
    P_0 \ \text{is} \ A_{\geq 7} & \Leftrightarrow P_1 \ \text{is} \ A_{\geq 5} \Leftrightarrow (*_4) \ \text{and}\  \theta = 0 \ (*_5).
\end{align*}

\paragraph{Step 4:} Set \begin{align*}
    \gamma_2 & \coloneqq -2\rho_2\xi_3 - \rho_{23}\xi_2+3\sigma_0\beta_3^2+2\sigma_1\beta_2\beta_3+\sigma_2\beta_2^2-2\varepsilon_2\beta_3-\varepsilon_{23}\beta_2 , \\
    \gamma_3 & \coloneqq - \rho_{23}\xi_3 -2\rho_3\xi_2 +3\sigma_3\beta_2^2+2\sigma_2\beta_2\beta_3+\sigma_1\beta_3^2-2\varepsilon_3\beta_2-\varepsilon_{23}\beta_3 \ \text{and}
\end{align*}
\begin{multline*}
    \mu \coloneqq -3\sigma_0\beta_3\xi_3^2-3\sigma_3\beta_2\xi_2^2-2\sigma_1\beta_3\xi_2\xi_3 - 2\sigma_2\beta_2\xi_2\xi_3 -\sigma_1\xi_3^2\beta_2 - \sigma_2\xi_2^2\beta_3 + \varepsilon_2\xi_3^2 + \varepsilon_{23}\xi_2\xi_3 +\varepsilon_3\xi_2^2 \\
    -3\tau_0\beta_3^2\xi_3-3\tau_3\beta_2^2\xi_2+\tau_1(-2\beta_2\beta_3\xi_3 - \beta_3^2\xi_2)+\tau_2(-\xi_3\beta_2^2-2\beta_2\beta_3\xi_2) + c_4(\beta_3,\beta_2).
\end{multline*} 

In this step 
\begin{center}
    $D_4 = \{x_2x_3+\gamma_2x_2+\gamma_3x_3+\mu + (\text{higher order terms in}\ x_1) = 0   \} \subset \Aa_{(x_1,x_2,x_3)}^3$.
\end{center}

We have that

\begin{center}
    $E_4 \cap D_4 = \{x_1 =x_2x_3+\gamma_2x_2+\gamma_3x_3+\mu = 0\}$, 
\end{center}

\nin is a conic (not necessarily irreducible) on $D_4 \subset \Aa^3$. One has
\begin{align*}
    P_0 \ \text{is} \ A_7 & \Leftrightarrow P_1 \ \text{is} \ A_5 \Leftrightarrow (*_5) \ \text{and}\  \gamma_2\gamma_3 \neq \mu,\\
    P_0 \ \text{is} \ A_{\geq 8} & \Leftrightarrow P_1 \ \text{is} \ A_{\geq 6} \Leftrightarrow (*_5) \ \text{and}\  \gamma_2\gamma_3 = \mu.
\end{align*}

In the latter case, we have 
\begin{center}
    $E_4 \cap D_4 = \{x_1 = (x_2+\gamma_3)(x_3+\gamma_2)= 0\} = L_4 \cup L_4'$,
\end{center}

\nin where $L_4$ and $L_4'$ are lines on $D_4 \subset \Aa^3$. Observe that $P_4 \coloneqq (0,-\gamma_3,-\gamma_2) \in \Sing(E_4 \cap D_4)$.\\

We will stop at this step. 

We observe that all the criteria obtained present symmetries. The degree of the conditions with respect to the coefficients of the quartic increases as $n$ increases.

In \cite{kn}, Kato \& Naruki found an equation of a quartic surface with a single $A_{19}$ singularity. We may use it to double-check all our criteria. 

In \cite{kn}, the equation of this quartic is in affine coordinates such that the tangent cone at $P$ is given by $\{x_1^2 + x_2^2 = (x_1 + ix_2)(x_1 - ix_2)=0\}$, namely,
\begin{multline*}
    16(x_1^2 +x_2^2) +32x_1x_3^2- 16x_2^3 +16x_3^4 -32x_2x_3^3+8(2x_1^2- 2x_1x_2+5x_2^2)x_3^2 \\ 
    +8(2x_1^3-5x_1^2x_2-6x_1x_2^2 -7x_2^3)x_3+20x_1^4 +44x_1^3x_2+65x_1^2x_2^2 +40x_1x_2^3 +41x_2^4 =0.
\end{multline*}

Performing the affine change of coordinates induced by the projective change of coordinates
\begin{center}
    $(x_0:x_1:x_2:x_3) \mapsto \left( x_0: \dfrac{x_2+x_3}{8} : \dfrac{i(x_2-x_3)}{8} :x_1 \right) $,
\end{center}

\nin we can bring this quartic to have $TC_P D = \{x_2x_3=0\}$. The equation becomes
\begin{multline}\label{quartic A_19}
    16 x_1^4 - 4 i x_1^3 x_2 + 4 i x_1^3 x_3 - \left(\dfrac{3}{8} + \dfrac{i}{4}\right) x_1^2 x_2^2 + \dfrac{7}{4} x_1^2 x_2 x_3 + 4 x_1^2 x_2 \\
    - \left(\dfrac{3}{8} - \dfrac{i}{4}\right) x_1^2 x_3^2 + 4 x_1^2 x_3 + \dfrac{i x_2^3}{32} - \dfrac{3}{32} i x_2^2 x_3 + \dfrac{3}{32} i x_2 x_3^2 + x_2 x_3 - \dfrac{i x_3^3}{32}\\
    +\left(\dfrac{1}{8} + \dfrac{i}{32} \right) x_1 x_2^3 - \dfrac{13}{32} i x_1 x_2^2 x_3 + \dfrac{13}{32} i x_1 x_2 x_3^2 + \left(\dfrac{1}{8} - \dfrac{i}{32} \right) x_1 x_3^3 \\
    + \left(- \dfrac{1}{1024} + \dfrac{i}{1024} \right) x_2^4 -\left(\dfrac{21}{1024} - \dfrac{21i}{512} \right) x_2^3 x_3 \\
    + \dfrac{31 x_2^2 x_3^2}{256} - \left(\dfrac{21}{1024} + \dfrac{21i}{512} \right) x_2 x_3^3 - \left(\dfrac{1}{1024} + \dfrac{i}{1024} \right) x_3^4 = 0.
\end{multline}

One can easily check that $P$ is indeed $A_{\geq 8}$ according to all our criteria.

\subsection{Toric description of the weights $(1,a,b)$}

In this part, we will analyze the toric description of the weighted blowup with weights $(1,a,b)$ and conditions imposed on its center so that it is volume-preserving. We will follow the notation in the corresponding paragraph section on Section \ref{vp Sark links}.

By \cite[Proposition 3.1]{acm} and  Lemma \ref{vp center sing pt}, the centers are:
\begin{center}
 $z_{i}=$
$\begin{cases} \hfil
\text{singular points in $D_i$ and its strict transforms}
,& \text{for $0 \leq i \leq a-1$},\\
\hfil \text{curves on $D_i$}, & \text{for $a \leq i \leq b-1$}.
\end{cases}$
\end{center}

Recalling our notation, consider $\pi \colon X \rightarrow \p^3$ the $(1,a,b)$-weighted blowup at $P \in \Sing(D)$.

Identify $\p^{34}$ with the space of quartics in $\p^3$ in the following way:
\begin{align*}
\p^{34} & \longrightarrow \text{\{quartics in $\p^3$\}}\\
(a_0:\ldots:a_{34})&\longmapsto \Biggl\{ \displaystyle \sum_{i=0}^{34} a_iM_i = 0 \Biggl\},
\end{align*}

\nin where $\{M_0,\ldots,M_{34}\}$ are all the monomials of degree 4 in $\C[x_0,x_1,x_2,x_3]$ with a fixed order.

Let $\mathcal{A}_{\geq n} \subset \p^{34}$ be the parameter space of irreducible quartic surfaces passing through $P=(1:0:0:0)$ and having an $A_{\geq n}$ singularity at this point.

\subsubsection{Toric description of the weights $(1,1,b)$}

\paragraph{Weights $(1,1,1)$:} We need to insert in $\Sigma_0$ the ray
\begin{center}
    $(1,1,1) \in \cone(v_1,v_2,v_3)$.
\end{center}

This insertion corresponds to blowing up the orbit of the torus action on $X_0 = \p^3$ associated to the 3-dimensional cone $\cone(v_1,v_2,v_3)$. This orbit is precisely 
\begin{center}
    $\{x_1 =0\} \cap \{x_2 =0\} \cap \{x_3 =0\} =  P = (1:0:0:0)$,
\end{center}
that is, $z_0 = P \eqqcolon P_0$.

By Proposition \ref{prop crep bir morp}, $\pi_1$ is volume-preserving.

Therefore the weights $(1,1,1)$ are volume-preserving for any quartic with an $A_n$ singularity at $P$, $n \geq 1$. 

\paragraph{Weights $(1,1,2)$:} We need to insert in $\Sigma_1$ the ray
\begin{center}
    $(1,1,2) \in \cone((1,1,1),(0,0,1))$.
\end{center}

This insertion corresponds to blowing up the orbit of the torus action on $X_1$ associated to the 2-dimensional cone $\cone((1,1,1),v_3)$. This orbit is precisely the line $L_1 = E_1 \cap \widetilde{\{x_3=0\}}$, that is, 
\begin{center}
    $z_1 = L_1 \subset E_1 \simeq \p^2$.
\end{center}

In order for this blowup to be volume-preserving, by \cite[Proposition 3.1]{acm} we need that $L_1 \subset E_1 \cap D_1$. 

If $P$ is a singularity of type $A_1$ or $A_2$, by \cite[Lemma 5.18]{alv} and Proposition \ref{res A_n}, $D_1 \subset X_1$ is a nonsingular K3 surface. For the $A_1$ case, we have that $E_1 \cap D_1$ is an irreducible conic. So the weights $(1,1,2)$ are not volume-preserving in this situation because $L_1$ is a line and $E_1 \cap D_1$ is an irreducible conic. Notice that this also follows from Lemma \ref{bound a+b leq n+1}.

Consider the $A_{\geq 2}$ case from now on. The following equations are with respect to the ``first affine chart'' $\{ y_1 \neq 0\}$ in $X_1$. The computations are analogous in the other ones.

One has 
\begin{align*}
        D_1 & = \{x_2x_3+x_1B(1,x_2,x_3)+x_1^2C(1,x_2,x_3) = 0\}\\
        & = \Biggl\{ x_2x_3+\displaystyle x_1\sum_{i=0}^3 b_i(x_2,x_3) + \displaystyle x_1^2\sum_{i=0}^4 c_i(x_2,x_3) = 0 \Biggl\} \subset \Aa_{(x_1,x_2,x_3)}^3
    \end{align*}
and $E_1 = \{x_1=0\}$.

Therefore
\begin{center}
    $E_1 \cap D_1 = \{x_1 = x_2x_3 = 0\} = \{x_1 = x_2 = 0\} \cup \{x_1 = x_3 = 0\}$.
\end{center}

Observe that $L_1=\{x_1 = x_3 = 0\}$ is indeed contained in $E_1 \cap D_1$.

By Proposition \ref{prop crep bir morp}, $\pi_2$ is
volume-preserving.

Therefore the weights $(1,1,2)$ are volume-preserving for any quartic in $\mathcal{A}_{\geq 2}$.

\paragraph{Weights $(1,1,3)$:} By Lemma \ref{bound a+b leq n+1} we are led to consider the $A_{\geq 3}$ case.

We need to insert in $\Sigma_2$ the ray
\begin{center}
    $(1,1,3) \in \cone((1,1,2),(0,0,1))$.
\end{center}

This insertion corresponds to blowing up the orbit of the torus action on $X_2$ associated to the 2-dimensional cone $\cone((1,1,2),v_3)$. This orbit is precisely the line $L_2 = E_2 \cap \widetilde{\{x_3=0\}}$, that is, 
\begin{center}
    $z_2 = L_2 \subset E_2 \simeq \p(\mathcal{N}_{L_1/X_1}^{\vee}) \simeq \p(\oo_{\p^1}(1) \oplus \oo_{\p^1}(-1)) \simeq \F_2$.
\end{center}

In order for this blowup to be volume-preserving, by \cite[Proposition 3.1]{acm} we need that $L_2 \subset E_2 \cap D_2$. 

In the ``first affine chart'' $\{ y_1 \neq 0\}$ in $X_2$, one has 
\begin{center}
    $D_2 = \Biggl\{ x_2x_3+\displaystyle \sum_{i=0}^3 b_i(x_2,x_1x_3) + \displaystyle x_1\sum_{i=0}^4 c_i(x_2,x_1x_3) = 0 \Biggl\} \subset \Aa_{(x_1,x_2,x_3)}^3$
\end{center}
and $E_2 = \{x_1=0\}$.

Therefore
\begin{center}
    $E_2 \cap D_2 = \{x_1 = x_2x_3 + b_0 + \beta_2x_2 + \rho_2x_2^2 + \sigma_0x_2^3= 0\}$
\end{center}
and it is a conic or cubic depending on whether $\sigma_0=0$.

One can check
\begin{align*}
    L_2 = \{ x_1 = x_3 = 0 \} \subset E_2 \cap D_2 & \Leftrightarrow b_0 = \beta_2 = \rho_2 = \sigma_0 = 0 \\
    & \Leftrightarrow x_3 \mid B \\
    & \Leftrightarrow D_1 \ \text{is tangent to $\{x_3=0\}$ along}\ L_1 .
\end{align*}

The last condition is a consequence of the geometric properties of the blowup. 

Note that $L_2 \subset E_2 \setminus E_1^2$. Indeed $E_1^2$ does not appear in the affine chart $\{y_1 \neq 0\}$, that is, 
\begin{center}
    $E_1^2 \cap (\text{affine chart}\ \{y_1 \neq 0\}) = \emptyset$,
\end{center}

and $E_1^2 \subset (\text{affine chart}\ \{y_3 \neq 0\})$, while
\begin{center}
    $L_2 \cap (\text{affine chart}\ \{y_3 \neq 0\}) = \emptyset$.
\end{center}

If $L_2 \subset E_2 \cap D_2$, then Proposition \ref{prop crep bir morp} implies that $\pi_3$ is
volume-preserving.

Therefore the weights $(1,1,3)$ are volume-preserving for any quartic in $\mathcal{A}_{\geq 3}$ such that $b_0 = \beta_2 = \rho_2 = \sigma_0 = 0$. We point out that the condition $b_0=0$ on an element of $\cA_{\geq 2}$ implies that it belongs to $\cA_{\geq 3}$. Conversely,
\begin{center}
    $\cA_{\geq 2} \cap \{b_0=0\} = \cA_{\geq 3} \subset \p^{34}$.
\end{center}

So the weights $(1,1,3)$ are not volume-preserving for a generic $D$ in the corresponding parameter space.

Besides the relation $b_0=0$, we need the extra closed conditions $\beta_2 = \rho_2 = \sigma_0 = 0$ on an element of $\cA_{\geq 3}$ so that the weights $(1,1,3)$ are volume-preserving. 

\begin{center}
    $(1,1,3)$ are volume-preserving weights for $D \in \cA_{\geq 3} \Leftrightarrow x_3 \mid B$. 
\end{center}

Observe that if $x_3 \mid B$, then $x_3 \nmid C$. Otherwise, $x_3$ divides the equation that defines $D$, and therefore $D$ would be reducible.
Write
\begin{center}
    $C = \displaystyle \sum_{i=0}^4 c_i'x_3^{4-i}$,
\end{center}
where $c_i' \in \C[x_1,x_2]_i$. The condition $x_3 \nmid C$ implies that $c_4' \neq 0$. Since $c_4'$ is a homogeneous polynomial in 2 variables over an algebraically closed field, we can factorize it into linear factors
\begin{align*}
    c_4' & = \upsilon_0x_1^4 + \upsilon_1x_1^3x_2 + \upsilon_2x_1^2x_2^2 + \upsilon_3x_1x_2^3 + \upsilon_4x_2^4 \\
    & = \displaystyle \prod_{i=1}^4 (\alpha_ix_1+\varrho_ix_2),
\end{align*}
where $\upsilon_i,\alpha_i,\varrho_i \in \C$ for all $i$. Take $\ell_i \coloneqq \{ x_3 = \alpha_ix_1+\varrho_ix_2 = 0\}$ and notice that $D \supset \ell_i$. 

Thus, an element of $\cA_{\geq 3}$ for which the weights $(1,1,3)$ are volume-preserving necessarily contains lines through the singular point $P$. The union of these lines $\ell_i$ constitutes a hyperplane section of $D$.

\paragraph{Weights $(1,1,4)$:} By Lemma \ref{bound a+b leq n+1} we are led to consider the $A_{\geq 4}$ case. From the previous case, we must have $b_0 = \beta_2 = \rho_2 = \sigma_0 = 0$. We will consider elements in 
\begin{center}
    $\cA_{\geq 4} \cap \{\beta_2 = \rho_2 = \sigma_0 = 0\} \subset \p^{34}$.
\end{center}

We suppressed the condition $\{b_0 = 0\}$ because it is already satisfied for elements in $\cA_{\geq 4}$.

We need to insert in $\Sigma_3$ the vector
\begin{center}
    $(1,1,4) \in \cone((1,1,3),(0,0,1))$.
\end{center}

This insertion corresponds to blowing up the orbit of the torus action on $X_3$ associated to the 2-dimensional cone $\cone((1,1,3),v_3)$. This orbit is precisely the line $L_3 = E_3 \cap \widetilde{\{x_3=0\}}$, that is, 
\begin{center}
    $z_3 = L_3 \subset E_3 \simeq \p(\mathcal{N}_{L_2/X_2}^{\vee}) \simeq \F_{m_2}$.
\end{center}

For this blowup to be volume-preserving, by \cite[Proposition 3.1]{acm} we need that $L_3 \subset E_3 \cap D_3$. 

In the ``first affine chart'' $\{ y_1 \neq 0\}$ in $X_3$, one has 
\begin{center}
    $D_3 = \Biggl\{ x_2x_3+\displaystyle x_1\sum_{i=1}^3 b_i'(x_2,x_3) + \displaystyle \sum_{i=0}^4 c_i(x_2,x_1^2x_3) = 0 \Biggl\} \subset \Aa_{(x_1,x_2,x_3)}$,
\end{center}
where $b_i'(x_2,x_3) \in \C[x_2,x_3]_i$ is such that $x_1^2b_i'(x_2,x_3)=b_i(x_2,x_1^2x_3)$ for $i \in \{1,2,3\}$, and $E_3 = \{x_1=0\}$. Recall that here we have $b_0 = \beta_2 = \rho_2 = \sigma_0 = 0$.

Therefore
\begin{center}
    $E_3 \cap D_3 = \{x_1 = x_2x_3 + c_0 + \delta_2x_2 + \varepsilon_2x_2^2 + \tau_0x_2^3 + \lambda_0x_2^4= 0\}$
\end{center}
and
\begin{align*}
    L_3 = \{ x_1 = x_3 = 0 \} \subset E_3 \cap D_3 & \Leftrightarrow c_0 = \delta_2 = \varepsilon_2 = \tau_0 = \lambda_0 = 0 \\
    & \Leftrightarrow x_3 \mid C .
\end{align*}

We have the following:

\begin{center}
    $(1,1,4)$ are volume-preserving weights for $D \in \cA_{\geq 4}$ \\
    $\Leftrightarrow (1,1,3)$ are volume-preserving weights for $D \in \cA_{\geq 4}\ \text{and}\ x_3 \mid C$\\
    $  \Leftrightarrow x_3 \mid B,C $. 
\end{center}

Observe that if $x_3 \mid B,C$, so $x_3$ divides the equation that defines $D$, and therefore $D$ would be reducible.

Therefore the weights $(1,1,4)$ are not volume-preserving for any element in $\cA_{\geq 4}$.

\subsubsection{Toric description of the weights $(1,2,b)$}

The toric description of the weighted blowup with weights $(1,2,b)$ says that we must start by inserting in $\Sigma_0$ the ray 
\begin{center}
    $(1,1,1) \in \cone(v_1,v_2,v_3)$.
\end{center}

We already discussed this step, which corresponds to the blowup of $P_0 = P$.

Then we need to insert in $\Sigma_1$ the ray
\begin{center}
    $(1,2,2) \in \cone((1,1,1),v_2,v_3)$.
\end{center}

This insertion corresponds to blowing up the orbit of the torus action on $X_1$ associated to the 3-dimensional cone $\cone((1,1,1),v_2,v_3)$. This orbit is precisely
\begin{center}
    $E_1 \cap \widetilde{\{x_2 = 0\}}  \cap \widetilde{\{x_3 = 0\}} = (0,0,0) \coloneqq P_1$
\end{center}
in the ``first affine chart'' $\{y_1 \neq 0\}$.

In order for this blowup to be volume-preserving, by Lemma \ref{vp center sing pt} we need that $P_1 \in \Sing(D_1)$.

In the ``first affine chart'' $\{ y_1 \neq 0\}$ in $X_1$, one has 
\begin{align*}
        D_1 & = \{x_2x_3+x_1B(1,x_2,x_3)+x_1^2C(1,x_2,x_3) = 0\}\\
        & = \Biggl\{ x_2x_3+\displaystyle x_1\sum_{i=0}^3 b_i(x_2,x_3) + \displaystyle x_1^2\sum_{i=0}^4 c_i(x_2,x_3) = 0 \Biggl\} \subset \Aa_{(x_1,x_2,x_3)}^3
    \end{align*}
and $E_1 = \{x_1=0\}$.

Notice that $P_1 \in D_1$ and by the Jacobian Criterion,
\begin{center}
    $P_1 \in \Sing(D_1) \Leftrightarrow b_0 = 0$.
\end{center}

In this case, Proposition \ref{prop crep bir morp} implies that $\pi_2$ is
volume-preserving.

\paragraph{Weights $(1,2,3)$:} By Lemma \ref{bound a+b leq n+1} we are led to consider the $A_{\geq 4}$ case.

After inserting the rays $(1,1,1)$ and $(1,1,2)$, we need to insert in $\Sigma_2$ the ray
\begin{center}
    $(1,2,3) \in \cone((1,2,2),(0,0,1))$.
\end{center}

This insertion corresponds to blowing up the orbit of the torus action on $X_2$ associated to the 2-dimensional cone $\cone((1,2,2),v_3)$. This orbit is precisely the line $L_2 = E_2 \cap \widetilde{\{x_3=0\}}$, that is, 
\begin{center}
    $z_2 = L_2 \subset E_2 \simeq \p^2$.
\end{center}

In order for this blowup to be volume-preserving, by \cite[Proposition 3.1]{acm} we need that $L_2 \subset E_2 \cap D_2$. 

If $P$ is a singularity of type $A_3$, by \cite[Lemma 5.18]{alv}  and Proposition \ref{res A_n}, $D_2 \subset X_2$ is nonsingular and $E_2 \cap D_2$ is an irreducible conic. So the weights $(1,2,3)$ are not volume-preserving. Suppose now $P$ is a singularity of type $A_{\geq 4}$.

In the ``first affine chart'' $\{ y_1 \neq 0\}$ in $X_2$, one has 
\begin{align*}
        D_2 & = \Biggl\{x_2x_3+\sum_{i=1}^3 x_1^{i-1}b_i(x_2,x_3) + \displaystyle \sum_{i=0}^4 x_1^ic_i(x_2,x_1x_3) = 0 \Biggl\}\\
        & = \Biggl\{ x_2x_3+ b_1(x_2,x_3) + x_1b_2(x_2,x_3)+x_1^2b_1(x_2,x_3) + \displaystyle \sum_{i=0}^4 x_1^ic_i(x_2,x_3) = 0 \Biggl\} \subset \Aa_{(x_1,x_2,x_3)}^3
    \end{align*}
and $E_2 = \{x_1=0\}$.

Therefore
\begin{center}
    $E_2 \cap D_2 = \{x_1 = x_2x_3 + \beta_2x_2 + \beta_3x_3 + c_0= 0\}$
\end{center}
and
\begin{center}
    $ L_2 = \{ x_1 = x_3 = 0 \} \subset E_2 \cap D_2 \Leftrightarrow  \beta_2 = c_0 = 0$.
\end{center}

If $L_2 \subset E_2 \cap D_2$, then Proposition \ref{prop crep bir morp} implies that $\pi_3$ is volume-preserving.

Therefore the weights $(1,2,3)$ are volume-preserving for any quartic in $\mathcal{A}_{\geq 4}$ such that $\beta_2 = c_0 = 0$. We point out that the condition $c_0=\beta_2\beta_3$ on an element of $\cA_{\geq 3}$ implies that it belongs to $\cA_{\geq 4}$. Conversely,
\begin{center}
    $\cA_{\geq 3} \cap \{c_0=\beta_2\beta_3\} = \cA_{\geq 4} \subset \p^{34}$.
\end{center}

So the weights $(1,2,3)$ are not volume-preserving for a generic $D$ in the corresponding parameter space.

Besides the relation $c_0=\beta_2\beta_3$, we need the extra closed conditions $\beta_2 = c_0 = 0$ on an element of $\cA_{\geq 4}$ so that the weights $(1,2,3)$ are volume-preserving. 
\begin{center}
    $(1,2,3)$ are volume-preserving weights for $D \in \cA_{\geq 4} \Leftrightarrow \beta_2 = c_0 = 0$. 
\end{center}

\paragraph{Weights $(1,2,5)$:} By Lemma \ref{bound a+b leq n+1} we are led to consider the $A_{\geq 6}$ case. 

According to the toric description, we need to insert in $\Sigma_3$ the ray
\begin{center}
    $(1,2,4) \in \cone((1,2,3),(0,0,1))$.
\end{center}

From the previous steps, we must have $b_0 = 0$ and $\beta_2 = c_0 = 0$.

These conditions imply that the quartic is in 
\begin{center}
    $\cA_{\geq 5} \cap \{\beta_2 = 0\} \subset \p^{34}$.
\end{center}

We suppressed the closed conditions $\{b_0 = 0\}$ and $\{c_0 = 0\}$ because the former is already satisfied for elements in $\cA_{\geq 4}$, and the latter follows from the fact in $\cA_{\geq 4}$ we have that 
\begin{center}
    $\beta_2 = 0 \Rightarrow c_0 = 0$.
\end{center}

Inserting the ray $(1,2,4)$ corresponds to blowing up the orbit of the torus action on $X_3$ associated to the 2-dimensional cone $\cone((1,2,3),v_3)$. This orbit is precisely the line $L_3 = E_3 \cap \widetilde{\{x_3=0\}}$, that is, 
\begin{center}
    $z_3 = L_3 \subset E_3 \simeq \p(\mathcal{N}_{L_2/X_2}^{\vee}) \simeq \p(\oo_{\p^1}(1) \oplus \oo_{\p^1}(-1)) \simeq \F_2$.
\end{center}

For this blowup to be volume-preserving, by \cite[Proposition 3.1]{acm} we need that $L_3 \subset E_3 \cap D_3$. 

In the ``first affine chart'' $\{ y_1 \neq 0\}$ in $X_3$, one has 
\begin{center}
    $ D_3 = \Biggl\{ (x_2+\beta_3)x_3 + b_2(x_2,x_1x_3)+x_1b_3(x_2,x_1x_3) + \displaystyle \sum_{i=1}^4 x_1^{i-1}c_i(x_2,x_1x_3) = 0 \Biggl\} \subset \Aa_{(x_1,x_2,x_3)}^3$
\end{center}
\nin and $E_3 = \{x_1=0\}$.

Therefore
\begin{center}
    $E_3 \cap D_3 = \{x_1 = (x_2+\beta_3)x_3 + \rho_2x_2^2 + \delta_2x_2= 0\}$
\end{center}
and
\begin{align*}
    L_3 = \{ x_1 = x_3 = 0 \} \subset E_3 \cap D_3 & \Leftrightarrow \rho_2 = \delta_2 = 0 \\
    & \Leftrightarrow D_2 \ \text{is tangent to $\{x_3=0\}$ along}\ L_2 .
\end{align*}

If $L_3 \subset E_3 \cap D_3$, then Proposition \ref{prop crep bir morp} implies that $\pi_4$ is volume-preserving.

Therefore $\pi_4$ is volume-preserving for any quartic in $\mathcal{A}_{\geq 5}$ such that $\rho_2 = \delta_2 = 0$. We point out that the condition $\zeta = 0$ on an element of $\cA_{\geq 4}$ implies that it belongs to $\cA_{\geq 5}$. Conversely,
\begin{center}
    $\cA_{\geq 4} \cap \{\zeta = b_2(\beta_3,\beta_2)-c_1(\beta_3,\beta_2) = \rho_2\beta_3^2 + \rho_{23}\beta_2\beta_3 + \rho_3\beta_2^2 - \delta_2\beta_3 - \delta_3\beta_2 = 0\} = \cA_{\geq 5} \subset \p^{34}$.
\end{center}

The blowup $\pi_4$ is not volume-preserving for a generic $D$ in $\cA_{\geq 5}$.

Besides the relation $\zeta = 0$, we need the extra closed conditions $\rho_2 = \delta_2 = 0$ on an element of $\cA_{\geq 5}$ so that $\pi_4$ is volume-preserving. 
\begin{center}
    $\pi_4$ is volume-preserving for $D \in \cA_{\geq 5} \cap 
\{\beta_2 = 0\} \Leftrightarrow \rho_2 = \delta_2 = 0$. 
\end{center}

Recall that by Lemma \ref{bound a+b leq n+1} we are led to consider the $A_{\geq 6}$ case so that the weights $(1,2,5)$ are volume-preserving. From the previous steps $(1,2,2), (1,2,3)$ and $(1,2,4)$, we must have $b_0 = 0$, $\beta_2 = c_0 = 0$ and $\rho_2 = \delta_2 = 0$.

Removing the redundant conditions, we must consider elements in 
\begin{center}
    $\cA_{\geq 6} \cap \{\beta_2 = \rho_2 = \delta_2 = 0\} \subset \p^{34}$.
\end{center}

We need to insert in $\Sigma_4$ the ray
\begin{center}
    $(1,2,5) \in \cone((1,2,4),(0,0,1))$.
\end{center}

This insertion corresponds to blowing up the orbit of the torus action on $X_4$ associated to the 2-dimensional cone $\cone((1,2,4),v_3)$. This orbit is precisely the line $L_4 = E_4 \cap \widetilde{\{x_3=0\}}$, that is, 
\begin{center}
    $z_4 = L_4 \subset E_4 \simeq \p(\mathcal{N}_{L_3/X_3}^{\vee}) \simeq \F_{m_3}$.
\end{center}

For this blowup to be volume-preserving, by \cite[Proposition 3.1]{acm} we need that $L_4 \subset E_4 \cap D_4$. 

In the ``first affine chart'' $\{ y_1 \neq 0\}$ in $X_4$, one has 
\begin{center}
    $ D_4 = \Biggl\{ (x_2+\beta_3)x_3 + \rho_{23}x_1x_2x_3 + \rho_3x_1^3x_3^2+b_3(x_2,x_1^2x_3) + \delta_3x_1x_3 + \displaystyle \sum_{i=2}^4 x_1^{i-2}c_i(x_2,x_1x_3) = 0 \Biggl\}$\\
    $\subset \Aa_{(x_1,x_2,x_3)}^3$
\end{center}
\nin and $E_4 = \{x_1=0\}$.

Therefore
\begin{center}
    $E_4 \cap D_4 = \{x_1 = (x_2+\beta_3)x_3 + \sigma_0x_2^3 + \varepsilon_2x_2^2= 0\}$
\end{center}
and
\begin{align*}
    L_4 = \{ x_1 = x_3 = 0 \} \subset E_4 \cap D_4 & \Leftrightarrow \sigma_0 = \varepsilon_2 = 0 \\
    & \Leftrightarrow D_3 \ \text{is tangent to $\{x_3=0\}$ along}\ L_3 .
\end{align*}

If $L_4 \subset E_4 \cap D_4$, then Proposition \ref{prop crep bir morp} implies that $\pi_5$ is volume-preserving.

Therefore the weights $(1,2,5)$ are volume-preserving for any quartic in $\mathcal{A}_{\geq 6}$ such that $\sigma_0 = \varepsilon_2 = 0$. We point out that the condition $\xi_2\xi_3 = \alpha$ on an element of $\cA_{\geq 5}$ implies that it belongs to $\cA_{\geq 6}$. Conversely,
\begin{center}
    $\cA_{\geq 5} \cap \{\xi_2\xi_3 = \alpha\} = \cA_{\geq 6} \subset \p^{34}$.
\end{center}

So the weights $(1,2,5)$ are not volume-preserving for a generic $D$ in $\cA_{\geq 6}$.

Besides the relation $\xi_2\xi_3 = \alpha$, we need the extra closed conditions $\sigma_0 = \varepsilon_2 = 0$ over an element of $\cA_{\geq 6}$ so that the weights $(1,2,5)$ are volume-preserving. 

\begin{center}
    $(1,2,5)$ are volume-preserving weights for $D \in \cA_{\geq 6} \cap \{\beta_2 = \rho_2 = \delta_2 = 0\} \Leftrightarrow \sigma_0 = \varepsilon_2 = 0$. 
\end{center}

At this point, for quartics in 
\begin{center}
     $D \in \cA_{\geq 6} \cap \{\beta_2 = \rho_2 = \delta_2 = \sigma_0 = \varepsilon_2 = 0\}$
\end{center}
the criteria $(*_5)$; $(*_5)$ and $\gamma_2\gamma_3 \neq \mu$; and $(*_5)$ and $\gamma_2\gamma_3 = \mu$ to detect whether $P$ is $A_{\geq 7}$, $A_7$ and $A_{\geq 8}$, respectively, become much simpler because for $\cA_{\geq 6} \cap \{\beta_2 = \rho_2 = \delta_2 = \sigma_0 = \varepsilon_2 = 0\}$ we have
\begin{itemize}
    \item $(*_5) \colon  (*_4)\ \text{and}\  \theta = -\tau_0\beta_3^3 $;
    \item $\gamma_2 = 0$;
    \item $\gamma_3 = \rho_{23}\xi_3 + \sigma_1\beta_3^2-\varepsilon_{23}\beta_3$;
    \item $\mu = \beta_3^2(-3\tau_0\xi_3+\lambda_0)$.
\end{itemize}

Since in this setting $\gamma_2\gamma_3=0$, to detect whether $P$ is $A_{\geq 8}$ it is enough to check if $\mu$ equals $0$. 

\subsubsection{Toric description of the weights $(1,3,b)$}

The toric description of the weighted blowup with weights $(1,3,b)$ says that we must start by inserting in $\Sigma_0$ the ray 
\begin{center}
    $(1,1,1) \in \cone(v_1,v_2,v_3)$,
\end{center}
followed by the insertion in $\Sigma_1$ of the ray 
\begin{center}
    $(1,2,2) \in \cone((1,1,2),v_2,v_3)$.
\end{center}

We already discussed these steps which correspond to the blowup of $P_0 = P$ and $P_1 = E_1 \cap \widetilde{\{x_2=0\}} \cap \widetilde{\{x_3 = 0\}}$, respectively.

Then we need to insert in $\Sigma_2$ the ray
\begin{center}
    $(1,3,3) \in \cone((1,2,2),v_2,v_3)$.
\end{center}

This insertion corresponds to blowing up the orbit of the torus action on $X_2$ associated to the 3-dimensional cone $\cone((1,2,2),v_2,v_3)$. This orbit is
\begin{center}
    $E_2 \cap \widetilde{\{x_2 = 0\}}  \cap \widetilde{\{x_3 = 0\}} = (0,0,0) \coloneqq P_2$.
\end{center}

For this blowup to be volume-preserving, by Lemma \ref{vp center sing pt} we need that $P_2 \in \Sing(D_2)$.

In the ``first affine chart'' $\{ y_1 \neq 0\}$ in $X_2$, one has 
\begin{align*}
        D_2 & = \Biggl\{x_2x_3+\sum_{i=1}^3 x_1^{i-1}b_i(x_2,x_3) + \displaystyle \sum_{i=0}^4 x_1^ic_i(x_2,x_3) = 0 \Biggl\}\\
        & = \Biggl\{ x_2x_3+ b_1(x_2,x_3) + x_1b_2(x_2,x_3)+x_1^2b_1(x_2,x_3) + \displaystyle \sum_{i=0}^4 x_1^ic_i(x_2,x_3) = 0 \Biggl\} \subset \Aa_{(x_1,x_2,x_2)}^3
    \end{align*}
and $E_2 = \{x_1=0\}$.

Notice that $P_2 \in D_2 \Leftrightarrow c_0 = 0$. By the Jacobian Criterion,
\begin{center}
    $P_2 \in \Sing(D_2) \Leftrightarrow \beta_2 = \beta_3 = 0$.
\end{center}

In this case, Proposition \ref{prop crep bir morp} implies that $\pi_3$ is volume-preserving.

\paragraph{Weights $(1,3,4)$:} By Lemma \ref{bound a+b leq n+1} we are led to consider the $A_{\geq 6}$ case. From the previous steps, we must have $b_0 = 0$ and $c_0 = \beta_2 = \beta_3 = 0$, that is, we must consider elements in 
\begin{center}
    $\cA_{\geq 6} \cap \{\beta_2 = \beta_3 = 0\} \subset \p^{34}$.
\end{center}

We suppressed the condition $\{b_0 = 0\}$ because it is already satisfied for elements in $\cA_{\geq 6}$, since we have that 
\begin{center}
    $\cA_{\geq 6} \subset \cA_{\geq 2}$;
\end{center}
as well as the condition $\{c_0 = 0\}$ because we have that
\begin{center}
    $\cA_{\geq 6} \subset \cA_{\geq 4}$,
\end{center}
since  the relation $c_0 = \beta_2\beta_3$ holds for elements in $\cA_{\geq 4}$.

We need to insert in $\Sigma_4$ the ray
\begin{center}
    $(1,3,4) \in \cone((1,3,3),(0,0,1))$.
\end{center}

This insertion corresponds to blowing up the orbit of the torus action on $X_3$ associated to the 2-dimensional cone $\cone((1,3,3),v_3)$. This orbit is precisely the line $L_3 = E_3 \cap \widetilde{\{x_3=0\}}$, that is, 
\begin{center}
    $z_3 = L_3 \subset E_3 \simeq \p^2$.
\end{center}

In order for this blowup to be volume-preserving, by \cite[Proposition 3.1]{acm} we need that $L_3 \subset E_3 \cap D_3$. 

If $P$ is a singularity of type $A_5$, by Proposition \ref{res A_n} and \cite[Lemma 5.18]{alv}, $D_3 \subset X_3$ is nonsingular and $E_3 \cap D_3$ is an irreducible conic. So the weights $(1,3,4)$ are not volume-preserving. 

Suppose now $P$ is a singularity of type $A_{\geq 6}$.

In the ``first affine chart'' $\{ y_1 \neq 0\}$ in $X_3$, one has 
\begin{center}
    $D_3 = \Biggl\{x_2x_3+\displaystyle \sum_{i=2}^3 x_1^{2i-3}b_i(x_2,x_3) + \displaystyle \sum_{i=1}^4 x_1^{2i-2}c_i(x_2,x_3) = 0 \Biggl\} \subset \Aa_{(x_1,x_2,x_3)}^3$
\end{center}
and $E_3 = \{x_1=0\}$.

Therefore
\begin{center}
    $E_3 \cap D_3 = \{x_1 = x_2x_3 + \delta_2x_2 + \delta_3x_3= 0\}$
\end{center}
and
\begin{center}
    $ L_3 = \{ x_1 = x_3 = 0 \} \subset E_3 \cap D_3 \Leftrightarrow  \delta_2 = 0$.
\end{center}

If $L_3 \subset E_3 \cap D_3$, then Proposition \ref{prop crep bir morp} implies that $\pi_4$ is volume-preserving.

Therefore the weights $(1,3,4)$ are volume-preserving for any quartic in $\mathcal{A}_{\geq 6}$ such that $\delta_2 = 0$. We point out that the condition $\theta = 0$ on an element of $\cA_{\geq 5}$ implies that it belongs to $\cA_{\geq 6}$. Conversely,
\begin{center}
    $\cA_{\geq 5} \cap \{\theta = 0\} = \cA_{\geq 6} \subset \p^{34}$.
\end{center}

So the weights $(1,3,4)$ are not volume-preserving for a generic $D$ in $\cA_{\geq 6}$.

Besides the relation $\theta$, we have the extra closed condition $\delta_2 = 0$ on an element of $\cA_{\geq 6}$ so that the weights $(1,3,4)$ are volume-preserving. 

\begin{center}
    $(1,3,4)$ are volume-preserving weights for $D \in \cA_{\geq 6} \cap \{\beta_2 = \beta_3 = 0\} \Leftrightarrow \delta_2 = 0$. 
\end{center}

\paragraph{Weights $(1,3,5)$:} By Lemma \ref{bound a+b leq n+1} we are led to consider the $A_{\geq 7}$ case. From the previous steps, we need to take into account $\beta_2 =\beta_3 = 0$ and $\delta_2 = 0$, that is, we must consider elements in 
\begin{center}
    $\cA_{\geq 6} \cap \{\beta_2 =\beta_3 = \delta_2 =0\} \subset \p^{34}$.
\end{center}

We need to insert in $\Sigma_4$ the ray
\begin{center}
    $(1,3,5) \in \cone((1,3,4),(0,0,1))$.
\end{center}

This insertion corresponds to blowing up the orbit of the torus action on $X_4$ associated to the 2-dimensional cone $\cone((1,3,4),v_3)$. This orbit is precisely the line $L_4 = E_4 \cap \widetilde{\{x_3=0\}}$, that is, 
\begin{center}
    $z_4 = L_4 \subset E_4 \simeq \p(\mathcal{N}_{L_3/X_3}^{\vee}) \simeq \p(\oo_{\p^1}(1) \oplus \oo_{\p^1}(-1)) \simeq \F_2$.
\end{center}

In order for this blowup to be volume-preserving, by \cite[Proposition 3.1]{acm} we need that $L_4 \subset E_4 \cap D_4$. 

In the ``first affine chart'' $\{ y_1 \neq 0\}$ in $X_4$, one has 
\begin{center}
    $ D_4 = \Biggl\{(x_2+\delta_3)x_3+\displaystyle \sum_{i=2}^3 x_1^{2i-4}b_i(x_2,x_1x_3) + \displaystyle \sum_{i=2}^4 x_1^{2i-3}c_i(x_2,x_1x_3) = 0 \Biggl\} \subset \Aa_{(x_1,x_2,x_3)}^3$
\end{center}
\nin and $E_4 = \{x_1=0\}$.

Therefore
\begin{center}
    $E_4 \cap D_4 = \{x_1 = (x_2+\delta_3)x_3 + \rho_2x_2^3= 0\}$
\end{center}
and
\begin{align*}
    L_4 = \{ x_1 = x_3 = 0 \} \subset E_4 \cap D_4 & \Leftrightarrow \rho_2 = 0 \\
    & \Leftrightarrow D_3 \ \text{is tangent to $\{x_3=0\}$ along}\ L_3 .
\end{align*}

If $L_4 \subset E_4 \cap D_4$, then Proposition \ref{prop crep bir morp} implies that $\pi_5$ is volume-preserving.

Therefore the weights $(1,3,5)$ are volume-preserving for any quartic in $\mathcal{A}_{\geq 7}$ such that $\rho_2 = 0$. We point out that the condition $\gamma_2\gamma_3 = \mu$ on an element of $\cA_{\geq 6}$ implies that it belongs to $\cA_{\geq 7}$. Conversely,
\begin{center}
    $\cA_{\geq 6} \cap \{\gamma_2\gamma_3 = \mu\} = \cA_{\geq 7} \subset \p^{34}$.
\end{center}

So the weights $(1,3,5)$ are not volume-preserving for a generic $D$ in $\cA_{\geq 6}$.

Besides the relation $\gamma_2\gamma_3 = \mu$, we have the extra condition $\rho_2 = 0$ on an element of $\cA_{\geq 7}$ so that the weights $(1,3,5)$ are volume-preserving. 

\begin{center}
    $(1,3,5)$ are volume-preserving weights for $D \in \cA_{\geq 7} \cap \{\beta_2 = \beta_3 = \delta_2 = 0\} \Leftrightarrow \rho_2 = 0$. 
\end{center}

At this point, for quartics in 
\begin{center}
     $D \in \cA_{\geq 7} \cap \{\beta_2 = \beta_3 = \delta_2 = \rho_2 = 0\}$
\end{center}
the criteria $(*_5)$ and $\gamma_2\gamma_3 \neq \mu$ and $(*_5)$ and $\gamma_2\gamma_3 = \mu$ to detect whether $P$ is $A_7$ and $A_{\geq 8}$, respectively, become much simpler because for $\cA_{\geq 7} \cap \{\beta_2 = \beta_2 = \delta_2 = \rho_2 = 0\}$ we have
\begin{itemize}
    \item $\gamma_2 = 0$;
    \item $\gamma_3 = -\rho_{23}\xi_3$;
    \item $\mu = \varepsilon_2\xi_3^2$.
\end{itemize}

Since in this setting $\gamma_2\gamma_3=0$, to detect whether $P$ is $A_{\geq 8}$ it is enough to check if $\mu$ equals $0$.\\

In the following Table \ref{table vp weights conditions A_n}, we summarize all the necessary and sufficient conditions so that the insertion of the ray $(1,c,d)$ corresponds to a volume-preserving blowup. If there do not exist any conditions, we will write ``generic'' as the case of the rays $(1,1,1)$ and $(1,1,2).$

\begin{table}[htp]
\begin{center}
\begin{tabular}{|c|c|}
\hline
ray inserted & conditions \\
\hline
(1,1,1) & generic on $\cA_{\geq 1}$ \\
(1,1,2) & generic on $\cA_{\geq 2}$\\
(1,1,3) & $b_0 = \beta_2 = \rho_2 = \sigma_0 = 0 \Leftrightarrow x_3 \mid B $ \\
(1,1,4) & $c_0 = \delta_2 = \varepsilon_2 = \tau_0 = \lambda_0 = 0 \Leftrightarrow x_3 \mid C $ \\
\hline
(1,2,2) & $b_0 = 0$ \\
(1,2,3) & $\beta_2 = c_0 = 0$ \\
(1,2,4) & $\rho_2 = \delta_2 = 0$ \\
(1,2,5) & $\sigma_0 = \varepsilon_2 = 0$ \\
\hline
(1,3,3) & $c_0 = \beta_2 = \beta_3 = 0$ \\
(1,3,4) & $\delta_2 = 0$ \\
(1,3,5) & $\rho_2 = 0$ \\
\hline
\end{tabular}
\caption{Table summarizing necessary and sufficient conditions for the induced blowup to be volume-preserving in the $A_n$ case.} 
\label{table vp weights conditions A_n}
\end{center}
\end{table}

\subsection{The $D$-$E$ case}

Let $D \subset \p^3$ be an irreducible normal quartic surface having a canonical singularity of type $D$-$E$ at $P=(1:0:0:0)$. 

Set $\pi \colon X \rightarrow \p^3$ to be the toric $(1,a,b)$-weighted blowup at $P$ and $E \coloneqq \Exc(\pi)$.

We will basically imitate the same strategy for the $A_n$ case. For the $D$-$E$ case, explicit results in terms of criteria to detect such singularities are more complicated to be stated, since the exceptional divisor of a minimal resolution has a different behavior along the process. Nevertheless, we will still be able to deduce some criteria to distinguish the cases $D_4$, $D_{>4}$ and $E_n$.

\subsubsection{Criteria for $D$-$E$ singularities on a quartic surface}\label{criteria D-E}

Keep the same notation concerning a quartic surface $D \subset \p^3$ having a canonical singularity of type $D$-$E$ at $P = (1: 0 : 0: 0)$. With all the considerations and remarks in this section, our task now is to establish conditions over $B$ and $C$ in such a way that we can guarantee a $D_4$, $D_{>4}$ and $E_n$ singularity.

As we have mentioned before, this was also already settled in a different setting regarding distinct coordinates for the $D_n$ case. Since we are working in different coordinates from \cite{kn}, we will obtain them by analyzing straightforwardly how the singularities are in the strict transform $D_1$ of $D$ with respect to the ordinary blowup $\pi_1 \colon X_1 \rightarrow \p^3$ of $P$. Set $E_1 \coloneqq \Exc(\pi_1)$ and consider the affine charts $W_i = X_1 \cap \{y_i \neq 0\}$.

From \cite[Appendix A]{alv}, we have the following Table \ref{table sing on D_1} which summarizes the behavior of the $D$-$E$ singularities, where, by abuse of notation, we identify the singular points in $D_1$ with their respective types.

\begin{table}[htp]
\begin{center}
\begin{tabular}{|c|c|}
\hline
type of singularity on $D$ & singularities on $D_1$  \\
\hline
$D_4$ & $A_1, A_1, A_1$ \\
\hline
$D_5$ & $A_1, A_3$\\
$D_{\ge 5}$ & $A_1, D_{n-2}$ \\
\hline
$E_6$ & $A_5$ \\
$E_7$ & $D_6$ \\
$E_8$ & $E_7$ \\
\hline
\end{tabular}
\caption{Table summarizing singularities on $D_1$.} 
\label{table sing on D_1}
\end{center}
\end{table}

From what was explained in Subsection \ref{notation quartic} concerning the notation for the coefficients of the homogeneous equation that defines $D$, we have that 
\begin{center}
    $P$ is $D\text{-}E \Leftrightarrow \rk(A) = 1$.
\end{center}

Such a condition is an explicit criterion and does not depend on $B$ and $C$. Without loss of generality, suppose $TC_PD= \{x_3^2=0\}$.

Regarding the three affine charts $W_i$ on $X_1$ and abusing notation, the equation $f_1=0$ of $D_1$ in these charts is given respectively by
\begin{align*}
    x_3^2 + x_1B(1,x_2,x_3) + x_1^2C(1,x_2,x_3)& =0, \\
    x_3^2 + x_2B(x_1,1,x_3) + x_2^2C(x_1,1,x_3) & =0\ \text{and} \\
    1+ x_3B(x_1,x_2,1) + x_3^2C(x_1,x_2,1) & =0.
\end{align*}

One has $E_1 \cap W_i = \{x_i = 0\}$. Notice that $E_1 \cap D_1  = \emptyset$ in $W_3$ and therefore $D_1$ has no singular points in this chart, which can also be justified by the Jacobian criterion.

Observe that $E_1 \cap D_1 = \{x_i = x_3^2 =0 \} $ in $W_i$ for $i \in \{1,2\}$, that is, such intersections are nonreduced and therefore all points belonging to them are candidates to be singular points of $D_1$ in the respective charts.

By the Jacobian criterion one can check that
\begin{align*}
    P_1=(0,x_2,0) \in \Sing(D_1) \cap W_1 & \Leftrightarrow x_2\  \text{is a root of}\ p_1(t) \coloneqq b_0+\beta_2t+\rho_2t^2+\sigma_0t^3 \in \C[t], \\
    P_1'=(x_1,0,0) \in \Sing(D_1) \cap W_2 & \Leftrightarrow x_1\  \text{is a root of}\ p_2(t) \coloneqq b_0t^3+\beta_2t^2+\rho_2t+\sigma_0 \in \C[t].
\end{align*}

We must necessarily have $p_1(t)$ and $p_2(t)$ are nonzero. Otherwise, the last two equivalences imply that $D_1$ is singular along the curve $E_1 \cap D_1$ and this contradicts the normality of $D_1$, which is retained in this case because $P$ is a canonical singularity.

Observe that $p_1(t)$ and $p_2(t)$ are reciprocal polynomials, that is, 
\begin{center}
    $\cf_{p_1}(t^i)=\cf_{p_2}(t^{3-i})$.
\end{center}

Set $d_i \coloneqq \deg(p_i(t))$ for $i \in \{1,2\}$.

Following the behavior described in Table \ref{table sing on D_1}, the degrees and the nature of the roots of $p_1$ and $p_2$ in terms of their multiplicity will determine whether $P$ is $D_4$, $D_{> 4}$ or $E_n$.  The last property is detected by means of the discriminant of such polynomials without computing their roots, which gives us the desired explicit criteria. 

We have four cases:\\

\noindent \textit{Case 1:} $d_1=3 \Rightarrow \sigma_0 \neq 0$

Since we are working in characteristic zero, we can find a change of variable such that it takes $p_1$ to a polynomial \begin{center}
    $q_1(t)=t^3+r_1t+s_1$,
\end{center}
where $r_1,s_1 \in \C$. In terms of the coefficients of $p_1$, they are explicitly 
\begin{align*}
    r_1 & = \dfrac{3\sigma_0\beta_2-\rho_2^2}{3\sigma_0^2}, \\
    s_1 & = \dfrac{2\rho_2^3-9\sigma_0\rho_2\beta_2+27\sigma_0^2b_0}{27\sigma_0^3}.
\end{align*}

It is straightforward that $p_1$ and $q_1$ have the same roots. Define $\Delta_1$ to be the discriminant of $q_1$. The discriminant of $p_1$ is the product $\sigma_0 \cdot \Delta_1$ and therefore it follows that one of these two discriminants is zero if and only if the other is also zero.

One computes $\Delta_1 = -(4r_1^3+27s_1^2)$. 

We have the following scenario and criteria:
\begin{align*}
    P\ \text{is}\ D_4 & \Leftrightarrow p_1\ \text{has 3 distinct roots} \Leftrightarrow \Delta_1 \neq 0,  \\
    P\ \text{is}\ D_{> 4} & \Leftrightarrow p_1\ \text{has 1 double root + 1 simple root, both distinct} \Leftrightarrow \Delta_1 = 0, r_1 \neq 0, \\
     P\ \text{is}\ E_n & \Leftrightarrow p_1\ \text{has 1 triple root} \Leftrightarrow \Delta_1 = 0, r_1 = 0\ (\Rightarrow s_1 = 0).
\end{align*}

\noindent \textit{Case 2:} $d_1=2 \Rightarrow \sigma_0 = 0, \rho_2 \neq 0$

We have the following scenario and criteria:
\begin{align*}
    P\ \text{is}\ D_4 & \Leftrightarrow p_1\ \text{has 2 distinct roots} \Leftrightarrow \beta_2^2-4\rho_2 b_0 \neq 0, \\
    P\ \text{is}\ D_{> 4} & \Leftrightarrow p_1\ \text{has 1 double root} \Leftrightarrow \beta_2^2-4\rho_2 b_0 = 0.
\end{align*}

In the remaining two cases $d_1=1 \Rightarrow \rho_2 = \sigma_0 = 0, \beta_2 \neq 0$ and $d_1 = 0 \Rightarrow \beta_2 = \rho_2 = \sigma_0 = 0, b_0 \neq 0$, we have that $P$ is $D_{>4}$, and $E_n$, respectively.

\subsubsection{Toric description of the weights $(1,a,b)$}

The analysis of the toric description of the weights $(1,a,b)$ in the $D$-$E$ case is analogous to the $A_n$ case. There only exists one difference concerning the weights of the form $(1,1,b)$ and some intermediate steps of the remaining ones when the centers are lines. Depending on some conditions imposed by the previously inserted rays, these corresponding lines coming from the intersection of the exceptional divisor and strict transforms of $D$ are nonreduced.

The toric description is with respect to nonsingular centers. Since these centers can be seen as abstract varieties (schemes over $\Spec(\C)$ in our case), then they are integral (reduced and irreducible) separated schemes of finite type over $\C$.

We will therefore take the \textit{reduced induced closed subscheme structure}, see \cite[Example 3.2.6]{har}, which is the most natural and smallest among all the possible \textit{closed subscheme structures}. Thus we will be blowing up such centers with the reduced induced structure. 

We point out that blowups along nonreduced schemes of a regular scheme can produce singularities and surprisingly they can also resolve singularities, as well as the classical blowups along reduced and regular schemes. We refer the reader to \cite[Section IV.2.3]{eh} and \cite[Section 22.4]{vak} for more details about these phenomena.

Let $\mathcal{D}_{\geq n} \subset \p^{34}$ be the parameter space of irreducible quartic surfaces passing through $P=(1:0:0:0)$ and having an $D_{\geq n}$ singularity at this point. Define $\cE_6,\cE_7$ and $\cE_8$ analogously, and abusing notation, set
\begin{center}
    $\cE_n \coloneqq \cE_6 \cup \cE_7 \cup \cE_8$.
\end{center}

Consider $f_i = 0$ the equation of the strict transform $D_i$ of $D$ at the $i$-th step of the toric description of the weighted blowup. By doing once again a lot of computations with ordinary blowups, we can summarize all the necessary and sufficient conditions for the weights analyzed are volume-preserving in the following Table \ref{table vp weights conditions D-E}. We also analyzed the weights of the form $(1,4,b)$ and restricted ourselves to the weights listed due to the feasibility of the computations. 

If there do not exist any conditions, we will write generic as in the case of the weights $(1,1,1)$, $(1,1,2)$, $(1,4,5)$. 
\vspace{20.0cm}

\begin{table}[htp]
\begin{center}
\begin{tabular}{|c|c|}
\hline
ray inserted & conditions \\
\hline
(1,1,1) & generic on $\mathcal{D}_{\geq 4}$ and $\cE_n$ \\
(1,1,2) & generic on $\mathcal{D}_{\geq 4}$ and $\cE_n$\\
(1,1,3) & $b_0 = \beta_2 = \rho_2 = \sigma_0 = 0 \Leftrightarrow x_3 \mid B \Leftrightarrow x_1 \mid f_2 $ \\
\hline
(1,2,2) & $b_0 = 0$ \\
(1,2,3) & $\beta_2 = c_0 = 0$ \\
(1,2,4) & $\rho_2 = \delta_2 = 0$ and $\beta_3 \neq 0$ 
\\
(1,2,5) & $\sigma_0 = \varepsilon_2 = 0$ \\
\hline
(1,3,3) & $c_0 = \beta_2 = \beta_3 = 0$ \\
(1,3,4) & $\delta_2 = 0$ \\
(1,3,5) & $\rho_2 = 0$ and $\delta_3 \neq 0$ \\
(1,3,6) & $\varepsilon_2 = 0$ \\
(1,3,7) & $\sigma_0 = 0$\\
\hline
(1,4,4) & $\delta_2 = \delta_3 = 0$ \\
(1,4,5) & generic on $\mathcal{D}_{\geq 9} \cap \{\text{previous conditions}\}$ and $\cE_8 \cap \{\text{previous conditions}\}$\\
(1,4,6) & $\rho_2 = 0 \Leftrightarrow x_1 \mid f_5$ \\
\hline
\end{tabular}
\caption{Table summarizing necessary and sufficient conditions for the induced blowup to be volume
preserving in the $D$-$E$ case.} 
\label{table vp weights conditions D-E}
\end{center}
\end{table}

In the insertion of the rays $(1,2,4)$ and $(1,3,5)$ we have the presence of the open conditions $\beta_3 \neq 0$ and $\delta_3 \neq 0$, respectively. Otherwise, we have that
\begin{center}
    $\rho_2 = \delta_2 = \beta_3 = 0 \Leftrightarrow x_1 \mid f_3$
\end{center}
and
\begin{center}
    $\rho_2 = \delta_3 = 0 \Leftrightarrow x_1 \mid f_4$,
\end{center}
respectively, which contradicts the irreducibility of $D_3$ and $D_4$, respectively.

The previous weights will exactly be the possible ones occurring in the $D_n$ case for $n \in \{4,\ldots,10\}$ and in the $E_n$ case. 

This can be justified in terms of how the singularities stand along the process of the toric description. Let us illustrate with some examples. The singularities appearing in the process are a consequence of the explicit resolution of the Du Val singularities in Section \ref{exp res Du Val}.
\begin{itemize}
    \item $D_5 \overset{\Bl_{P_0}} \longrightarrow A_3 \overset{\Bl_{P_1}} \longrightarrow A_1 \overset{\Bl_{P_2} \text{or}\ \Bl_{L_2}} \longrightarrow A_0 $
    \item $D_6 \overset{\Bl_{P_0}} \longrightarrow D_4 \overset{\Bl_{P_1}} \longrightarrow A_1 \overset{\Bl_{P_2} \text{or}\ \Bl_{L_2}} \longrightarrow A_0 $
\end{itemize}

In these examples, we can see directly that the weights of the form $(1,3,b)$ are not realizable as volume-preserving due to the toric description. In fact we should have that $z_3$ is a line $L_3$ such that  $L_3 = E_3 \cap D_3$ and $L_3 \not\subset E_2^3$. But $E_3 \cap D_3$ is a conic. 

In the second step in both previous examples, $P_1$ could be the remaining $A_1$ point on $D_1$, but the conditions over the coefficients so that the weights are volume-preserving will allow us to say that $P_1$ is indeed a singularity of type $A_3$ and $D_4$, respectively.

The following examples also allow us to restrict the possibilities for weights based on analogous arguments.

\begin{itemize}
    \item $E_6 \overset{\Bl_{P_0}} \longrightarrow A_5 \overset{\Bl_{P_1}} \longrightarrow A_3 \overset{\Bl_{P_2}} \longrightarrow A_1 \overset{\Bl_{P_3} \text{or}\ \Bl_{L_3}} \longrightarrow A_0 $
    \item $E_7 \overset{\Bl_{P_0}} \longrightarrow D_6 \overset{\Bl_{P_1}} \longrightarrow D_4 \overset{\Bl_{P_2}} \longrightarrow A_1 \overset{\Bl_{P_3} \text{or}\ \Bl_{L_3}} \longrightarrow A_0 $
    \item $E_8 \overset{\Bl_{P_0}} \longrightarrow E_7 \overset{\Bl_{P_1}} \longrightarrow D_6 \overset{\Bl_{P_2}} \longrightarrow D_4 \overset{\Bl_{P_3}} \longrightarrow A_1 \overset{\Bl_{P_4} \text{or}\ \Bl_{L_4}} \longrightarrow A_0$
\end{itemize}

\subsection{Proof of Theorem \ref{thm vp weights} \& Corollary \ref{cor vp weights sark}}

So far we have obtained explicit criteria to recognize certain canonical singularities and at the same time criteria so that the weights analyzed are volume-preserving. By comparing both criteria and plugging in Theorem \ref{gue thm}, we can obtain the desired weights listed in Tables \ref{table vp weights III} \& \ref{table vp weights IV}.

\begin{rmk}
For the unique quartic surface $D$ in $\p^3$, up to automorphisms of ambient space, with a singularity of type $A_{19}$, the only volume-preserving weights are $(1,1,1)$ and $(1,1,2)$ according to the criteria exposed in Table \ref{table vp weights conditions A_n}. Indeed, from its equation \ref{quartic A_19} we have that $\beta_2 = -4i \neq 0$ and therefore the weights $(1,1,3)$ and the weights of the form $(1,2,b)$ and $(1,3,b)$ are not volume-preserving.  

Notice that the Calabi-Yau pair $(\p^3,D)$ satisfies the assumptions of \cite[Proposition 3.3]{acm}. As a consequence of this result and Theorem \ref{cor vp weights sark}, the first volume-preserving Sarkisov link of any Sarkisov factorization of a volume-preserving map from this Calabi-Yau pair must start with a divisorial contraction which in local analytic coordinates is either the $(1,1,1)$ or the $(1,1,2)$-weighted blowup at the singular point of $D$. 

The same holds for a Calabi-Yau pair $(\p^3,D)$ of coregularity 2 such that $\rho_D = 1$ with the possibilities for weights restricted according to Table \ref{table vp weights IV}. 
\end{rmk}

\let\oldbibliography\thebibliography
\renewcommand{\thebibliography}[1]{\oldbibliography{#1}
\setlength{\itemsep}{0pt}} 

\noindent Universität Basel, Departement Mathematik und Informatik, Spiegelgasse 1, CH-4051
Basel, Switzerland \\
\textit{Email address}: \texttt{eduardo.alvesdasilva@unibas.ch}

\end{document}